 \documentclass[graybox]{svmult}
\usepackage{mathptmx}       
\usepackage{helvet}         
\usepackage{courier}        
\usepackage{type1cm}        
\usepackage{makeidx}         
\usepackage{graphicx}        
\usepackage{multicol}        
\usepackage[bottom]{footmisc}
 \usepackage{comment}
 \usepackage{latexsym}
 \usepackage{amsfonts}
 \usepackage{amssymb,amsmath}

\makeindex               
 \usepackage[latin1]{inputenc}
 \usepackage[T1]{fontenc}
 \usepackage[normalem]{ulem}
\usepackage{verbatim}
 \usepackage{graphicx}
\usepackage[latin1]{inputenc}
\usepackage{bbm}
\newcommand{\p}{\mathfrak{p}}

\renewcommand{\u}{\mathfrak{u}}
\renewcommand{\v}{\mathfrak{v}}
\newcommand{\C}{\mathbb{C}}
\newcommand{\N}{\mathbb{N}}
\newcommand{\R}{\mathbb{R}}

\newcommand{\T}{\mathbb{T}}

\newcommand{\Z}{\mathbb{Z}}

\newcommand{\boC}{\mathcal{C}}
\newcommand{\boE}{\mathcal{E}}

\newcommand{\boO}{\mathcal{O}}
\newcommand{\boP}{\mathcal{P}}

\newcommand{\boN}{\mathcal{N}}
\newcommand{\boM}{\mathcal{M}}

\newcommand{\eps}{\varepsilon}

\newcommand{\ch}{{\rm ch}}

\renewcommand{\th}{{\rm th}}

\begin{document}

\title*{IST versus PDE, a comparative study}
\titlerunning{IST-PDE}
\author{Christian Klein\inst{1}\and
Jean-Claude Saut\inst{2}}
\institute{Institut de Math\'ematiques de Bourgogne, UMR 5584\\
                Universit\'e de Bourgogne, 9 avenue Alain Savary, 21078 Dijon
                Cedex, France\\
   \texttt{Christian.Klein@u-bourgogne.fr}
\and
Laboratoire de Math\' ematiques, UMR 8628\\
Universit\' e Paris-Sud et CNRS\\ 91405 Orsay, France\\
\texttt{jean-claude.saut@math.u-psud.fr}}
\date{\today}
\maketitle
\large

\abstract{
We survey and compare, mainly in the two-dimensional case, various results obtained by IST and PDE techniques for integrable equations. We also  comment on what can be predicted from integrable equations on non integrable ones.
}

\vspace{0.3cm}
\begin{center}
{\it To Walter Craig with friendship and admiration}
\end{center}

\section{Introduction}
The theory of nonlinear dispersive equations has been flourishing 
during the last thirty years. Partial differential equations (PDE) techniques (in the large) have led to striking results concerning the resolution of the Cauchy problem, blow-up issues, stability analysis of various "localized" solutions.
On the other hand, a few nonlinear dispersive equations or systems are integrable by  Inverse Scattering Transform (IST) techniques. This allows a deep understanding of the equation dynamic and also to make relevant conjectures on close, non integrable equations. The best example is the Korteweg-de Vries (KdV) equation for which IST allows to prove that any solution to the Cauchy problem with  sufficiently smooth and decaying initial data decomposes into a finite train of solitons traveling to the right and a dispersive tail traveling to the left (see \cite{Schu} and the references therein  and section 1).

The aim of the present paper is to survey and compare, on specific examples, the advantages and shortcomings of PDE and IST techniques and also how they can benefit from each other. We will restrain to the Cauchy problem posed on the whole space $\R$ or $\R^2$ since the corresponding {\it periodic} problems lead to rather different issues, both by the two approaches.

\vspace{0.3cm}
The paper will be organized as follows. The first section is devoted 
to one-dimensional (spatial) problems. After recalling the KdV case  
for which the IST techniques yield the more complete results and 
which can be seen as a paradigm of what is expected for "close", 
possibly not integrable equations, we then consider two nonlocal 
integrable equations, the Benjamin-Ono (BO) and Intermediate Long 
Wave (ILW) equations that have a much less complete IST theory than 
the KdV equation. We will in this section present some numerical 
simulations from \cite{KS} showing that the long time dynamics of KdV 
solutions seems to be inherited by those of some non integrable 
equations, such as the fractional KdV or BBM equations. We close this 
Section by the one-dimensional Gross-Pitaevskii equation, a 
defocusing nonlinear Schr\"{o}dinger equation for which non trivial 
boundary conditions at infinity provide some {\it focusing} behavior. 
At this example, one can compare, for the specific problem of the 
stability of the {\it black soliton}, the differences between the two 
methods. We will provide here some details since this example might 
be less known than for instance the KdV equation.
 
We then turn to two-dimensional equations. Section two is devoted to the Kadomtsev-Petviashvili equations (KP).
Finally, the last section deals with the family of Davey-Stewartson systems, two members  of which are integrable (the so-called DSI and DS II system). An important issue is whether or not some of the remarkable properties of the integrable DS systems persist in the non integrable case.

We conclude by a short mention of two other integrable two-dimensional systems, the Ishimori and the Novikov-Veselov systems.


\section{Notations}

The following notations will be used throughout this article. The partial derivative will be denoted by $u_x,...$ or $\partial_x \phi,..$. For any $s\in \R,$ $D^s=(-\Delta)^{\frac{s}2}$ and $J^s=(I-\Delta)^{\frac{s}2}$ denote the Riesz and Bessel potentials of order $-s$, respectively. 

The Fourier transform of a function $f$ is denoted by $\hat{f} $ or 
$\mathcal{F}(f)$ and the dual variable of $x\in \R^d$ is denoted $\xi.$ For $1 \le p \le \infty$, $L^p(\mathbb R)$ is the 
usual Lebesgue space with the norm $|\cdot |_p$, and for $s \in 
\mathbb R$, the Sobolev spaces $H^s(\mathbb R^2)$ are defined via the usual norm $\|\phi \|_s= |J^s \phi|_2$. 

 $\mathcal S(\R^d)$ will denote the Schwartz space of smooth rapidly decaying functions in $\R^d$, and $\mathcal S'(\R^d)$ the space of tempered distributions.

\section{The one-dimensional case}

\subsection{The KdV equation}

The KdV equation is historically the first nonlinear dispersive equation which has been written down. It was  in fact derived  formally by Boussinesq (1877). We refer to \cite{Da} for a complete historical account. The full rigorous derivation from the water wave system is due to W. Craig \cite{Craig}. We refer to \cite{Da} for historical aspects and to \cite{La} for the systematic rigorous derivation of water waves models in various regimes. The KdV equation is (as in fact most of the classical  nonlinear dispersive equations or systems) a "universal" asymptotic equation describing a specific dynamic (in the long wave, weakly nonlinear regime) of a large class of complex nonlinear dispersive systems. \footnote{Note however that the KdV equation is a one-way model. For waves traveling in both directions, the same asymptotic regime would lead to the class of Boussineq systems (see {\it eg} \cite{BCS1} in the context of water waves) which are not integrable.}

The KdV equation

\begin{equation}\label{KdV}\begin{aligned}
&u_t-6uu_x+u_{xxx}=0,\\
\end{aligned}
\end{equation}

is also the first nonlinear  PDE for which the Inverse scattering technique was successively applied (see for instance \cite{GGKM, Lax}).

It is associated to the spectral problem for the Schr\"{o}dinger operator

$$L(t)=-\frac{d^2\psi}{dx^2} +u(\cdot,t)\psi$$

considered as an unbounded operator in $L^2(\R).$

We thus consider the spectral problem

$$\psi_{xx}+(k^2-u(x,t))\psi =0, \quad -\infty<x<+\infty.$$

Given $u_0=u(\cdot,0)$ sufficiently smooth and decaying at $\pm \infty,$ say in the Schwartz space $\mathcal S(\R),$ one associates to $L(0)$ its spectral data, that is a finite (possibly empty) set of negative eigenvalues $-\kappa_1^2<-\kappa^2_2\cdot \cdot \cdot <-\kappa_N^2,$ together with right normalization coefficients $c_j^r$ and right reflection coefficients $b_r(k)$ (see \cite{Schu} for precise definitions and properties of those objects). 

The spectral data consists thus in the collection of $\lbrace b_r(k),\kappa_j, c_j^r\rbrace.$ It turns out that if $u(x,t)$ evolves according to the KdV equation, the scattering data evolves in a very simple way:

$$\kappa_j(t)=\kappa_j,$$

$$c_j^r(t)=c_j^r\exp (4\kappa_j^3 t),\quad j=1,2,\cdot \cdot \cdot,N,$$

$$b_r(k,t)=b_r(k)\exp(8ik^3t), \quad -\infty<k<+\infty.$$

The potential $u(x,t)$ is recovered as follows. Let

$$\Omega(\xi;t)=2\sum_{j=1}^N [c^r_j(t)]^2e^{-2\kappa_j\xi}+\frac{1}{\pi}\int_{-\infty}^{\infty} b_r(k,t)e^{2ik\xi} dk.$$

One then solves the linear integral equation (Gel'fand-Levitan- Marchenko equation):

\begin{equation}\label{GL}
\beta(y;x,t)+\Omega(x+y;t)+\int_0^\infty \Omega(x+y+z;t)\beta(z;x,t)dz=0, \quad y>0,\;x\in \R,\;t>0.
\end{equation}

The solution of the Cauchy problem \eqref{KdV} is then given by 

$$u(x,t)=-\frac{\partial}{\partial x}\beta(0^+;x,t),\quad x\in \R,\quad t>0.$$

One obtains explicit solutions when $b_r=0.$ A striking case is 
obtained when the scattering data are $\lbrace 
0,\kappa_j,c_j^r(t)\rbrace.$ This corresponds to the so-called $N-soliton$ solution $u_d(x,t)$ according to its asymptotic behavior obtained by Tanaka \cite{Tana}:

\begin{equation}\label{Tanaka}
\lim_{t\to \infty}\;\sup_{x\in \R}\left|u_d(x,t)-\sum_{p=1}^N\left(-2\kappa_p^2\text{sech}^2[\kappa_p(x-x_p^+-4\kappa_p^2t)]\right)\right|=0,
\end{equation}

where

$$x_p^+=\frac{1}{2\kappa_p}\log \left\lbrace \frac{\left[c_p^r\right]^2}{2\kappa_p}\prod_{l=1}^{p-1}\left(\frac{\kappa_l-\kappa_p}{\kappa_l+\kappa_p}\right)^2\right\rbrace.$$

In other words, $u_d(x,t)$ appears for large positive time as a sequence of $N$ solitons, with the largest one in the front, uniformly with respect to $x\in \R.$ 

For $u_0\in \mathcal S(\R)$\footnote{This condition can be weakened, but a decay property is always needed.}, the solution of \eqref{KdV} has the following asymptotics

\begin{equation}\label{nosol}
\sup_{x\geq -t^{1/3}}|u(x,t)|=O(t^{-2/3}), \quad \text{as}\; t\to 
\infty,
\end{equation}

 in the absence of solitons (that is when $L(0)$ has no negative eigenvalues) and
 
 \begin{equation}\label{avecso}
 \sup_{x\geq -t^{1/3}}|u(x,t)-u_d(x,t)|=O(t^{-1/3}),\quad \text{as} \; t\to \infty
 \end{equation}
 
 in the general case, the $N$ in $u_d$ being the number of negative eigenvalues of $L(0).$  One has moreover the convergence result
 
 \begin{equation}\label{super}
 \lim_{t\to +\infty}\;\sup_{x\geq -t^{1/3}}\left|u(x,t)-\sum_{p=1}^N\left(-2\kappa_p^2\text{sech}^2[\kappa_p(x-x_p^+-4\kappa_p^2t)]\right)\right|=0.
 \end {equation}
 
 In both cases, a "dispersive tail" propagates to the left.
 
 \begin{remark}
 The shortcoming of those remarkable results is of course that they apply only to the integrable KdV equation and also to the modified KdV equation
 
 $$  u_t+6u^2u_x+u_{xxx}=0. $$
 
 However, though they are out of reach of "classical" PDE methods, they give  hints on the behavior of other, non integrable, equations whose dynamics could be in some sense similar.
 \end{remark}
 
 \begin{remark}
 As previously noticed, the results obtained by IST methods necessitate a decay property of the initial data, the minimal condition being
 
 \begin{equation}\label{min}
I(0)= \int_{-\infty}^\infty (1+|x|)|u_0(x)|dx<\infty.
 \end{equation}
 
 This condition ensures in particular (\cite{Marc}) that $L(0)$ has a 
 finite number of discrete eigenvalues, more 
 precisely (\cite{CaDe}), the number $N$ of eigenvalues of $L(0)$ is 
 bounded  by $1+I(0)$.
 
 This excludes for instance initial data in the energy space $H^1(\R)$ in which the Cauchy problem is globally well-posed (\cite{KPV2}). The global behavior of the flow might  thus be different from the aforementioned results for such initial data. 
 \end{remark}
 
 \begin{remark}
 1. Stemming from the seminal work of Bourgain \cite{Bourgain}, PDE 
 techniques allow to prove the well-posedness of the Cauchy problem 
 for the KdV equation with initial data in very large spaces, namely 
 $H^s(\R), s>-\frac{3}{4},$ (see \cite{KPV3})  which includes in 
 particular measures, for instance the Dirac distribution.

2. PDE techniques yield also the asymptotic stability of the solitary waves for subcritical KdV equations with a rather general nonlinearity (\cite{ MaMe, MaMe1}) which can be seen as a first step towards the {\it soliton resolution conjecture}, see {\it eg} \cite {Tao}.
 \end{remark}
 
 In order to see to what extent the long time dynamics of the KdV 
 equation is in some sense generic, we will consider as a toy model 
 the fractional KdV (fKdV) equation
 
 \begin{equation}\label{fKdV}
u_t+uu_x-D^\alpha u_x=0,\quad u(.,0)=u_0,
\end{equation}

where $\widehat{D^\alpha f}(\xi)=|\xi|^\alpha \hat f(\xi),\; \alpha>-1.$

Using the Fourier multiplier operator notation

$$\widehat{Lu}(\xi)=p(\xi)\hat{u}(\xi),\quad p(\xi)=|\xi|^\alpha,$$

it can be rewritten as

 \begin{equation}\label{fKdVbis}
u_t+uu_x-L u_x=0,\quad u(.,0)=u_0,
\end{equation}

When $\alpha=2$ (resp.$1$) \eqref {fKdV} reduces to the KdV (resp. Benjamin-Ono)  equation. If the symbol $|\xi|^\alpha$ is replaced by 

$$p(\xi)=\left( \frac{\tanh \xi}{\xi} \right)^{1/2} ,$$

one gets the so-called Whitham equation (\cite{Wh}) that models surface gravity waves in an appropriate regime. This symbol behaves like $|\xi|^{-1/2}$ for large frequencies.

 When surface tension is included in the Whitham equation, one gets 

$$p(\xi)=(1+\beta|\xi|^2)^{1/2}\left( \frac{\tanh \xi}{\xi} \right)^{1/2} ,\quad \beta\geq 0$$
which behaves like $|\xi|^{1/2}$ for large $|\xi|.$

The following quantities are formally conserved by the flow associated to \eqref{fKdV},
\begin{equation} \label{M}
M(u)=\int_{\mathbb R}u^2(x,t)dx,
\end{equation}
and the Hamiltonian

\begin{equation} \label{H}
H(u)=\int_{\mathbb R}\big( \frac{1}{2} |D^{\frac{\alpha}2}u(x,t)|^2-\frac{1}{6}u^3(x,t)\big) dx.
\end{equation}

One notices that the values $\alpha=1/3$ and $\alpha =1/2$ correspond respectively to the so-called energy critical and to the $L^2$ critical cases. Actually, equation \eqref{fKdV} is invariant under the scaling transformation 
\begin{equation} \label{scaling}
u_{\lambda} (x,t)=\lambda^{\alpha}u(\lambda x,\lambda^{\alpha+1}t),
\end{equation}
for any positive number $\lambda$. A straightforward computation shows that $\|u_{\lambda}\|_{\dot{H}^s}=\lambda^{s+\alpha-\frac{1}{2}}\|u_{\lambda}\|_{\dot{H}^s}$, and thus the critical index corresponding to \eqref{fKdV} is $s_{\alpha}=\frac{1}{2}-\alpha$. Thus, equation \eqref{fKdV} is $L^2$-critical for $\alpha=\frac{1}{2}$. On the other hand the Hamiltonian does not make sense in the {\it energy space} $H^{\alpha/2}(\R)$ when $\alpha <\frac{1}{3}$. The numerical simulations in \cite{KS} suggest that the Cauchy problem \eqref{fKdV} has global solutions (for arbitrary large suitably localized and smooth initial data) if and only if $\alpha>1/2.$ This has been rigorously proven when $\alpha\geq 1$ (see \cite{FLP1, FLP2}) but is an open problem when $1/2<\alpha<1.$ On the other hand, the local Cauchy problem is for $\alpha >0$ locally well-posed in $H^s(\R), s>\frac{3}{2}-\frac{3\alpha}{8},$  (\cite{LPS}).

More surprising is the fact that the resolution into solitary waves 
plus dispersion seems to be still valid when $\alpha>1/2$ as also 
suggested from the  numerical simulations in \cite{KS} from which we 
extract the following figures. In Fig.~\ref{fKdValpha06}, one can see 
the solution for the fKdV equation in the mass subcritical case 
$\alpha=0.6$ for the initial data $5\text{sech}^2 x$.
\begin{figure}[htb!]
    \includegraphics[width=\textwidth]{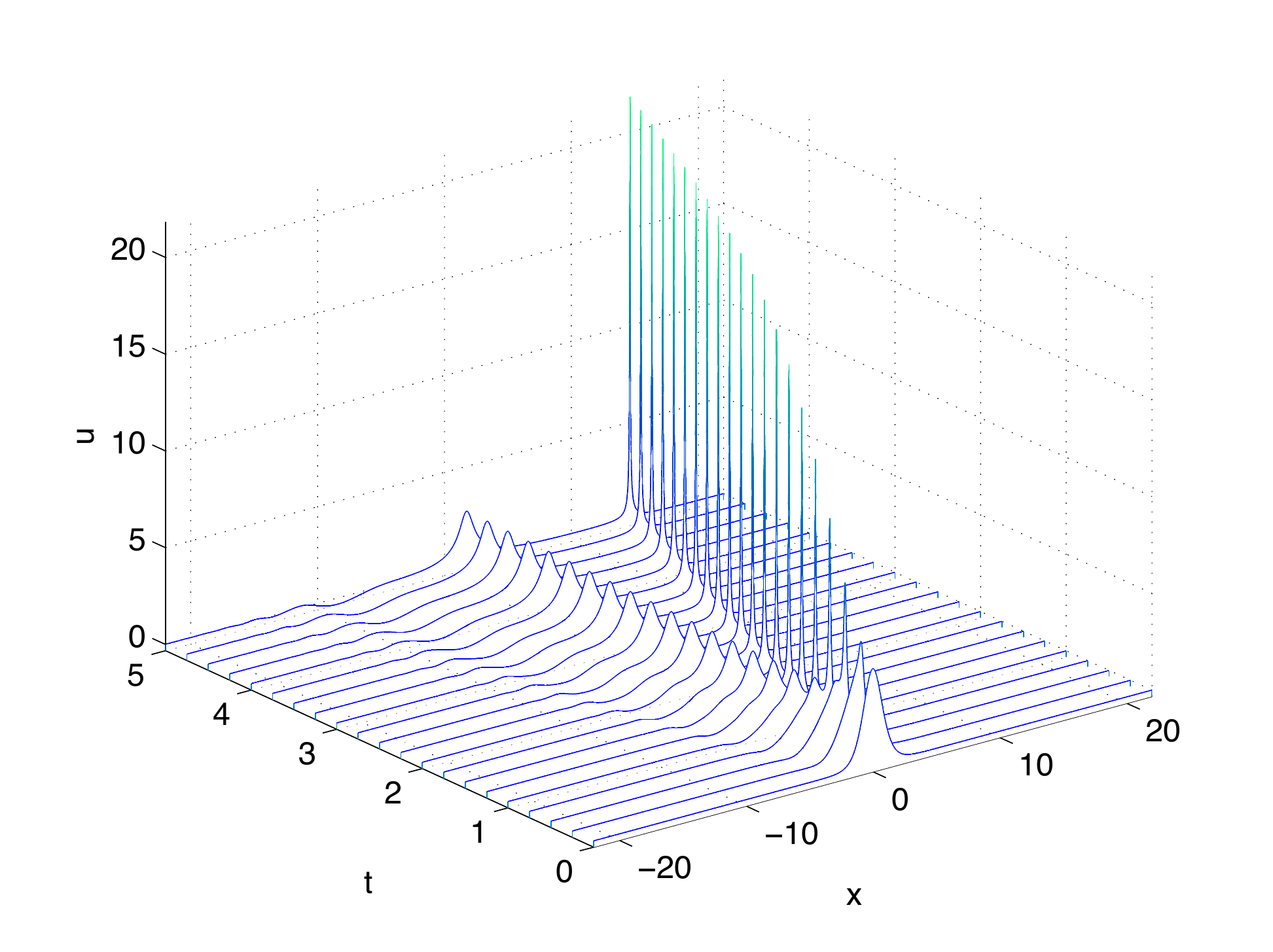}
 \caption{Solution to the fKdV equation  
 \eqref{fKdV} for $\alpha=0.6$, for the initial 
 data $\u_{0}=5\text{sech}^2 x.$}
 \label{fKdValpha06}
\end{figure}

\begin{figure}[htb!]
 \includegraphics[width=\textwidth]{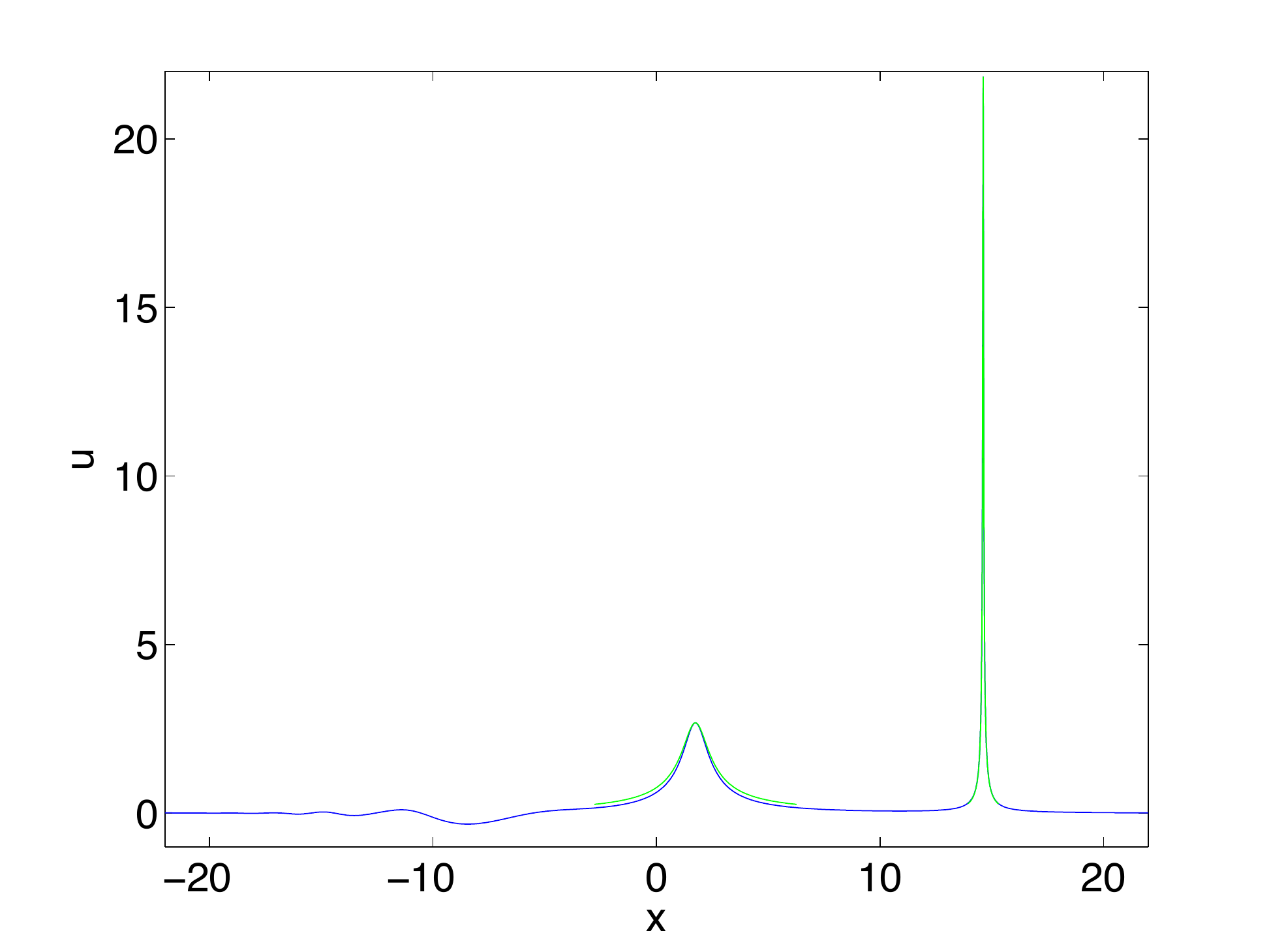}
 \caption{Solution to the fKdV equation  
 \eqref{fKdV} for $\alpha=0.6$, and the initial 
 data $\u_{0}=5\text{sech}^2 x$ for $t=5$ in blue, fitted solitons at the humps in green. }
 \label{fKdV06fit}
\end{figure}

In Fig.~\ref{fKdV06fit} we have fitted the humps with the computed solitary waves. This is an evidence for the above mentioned  {\it soliton resolution conjecture}.

A similar behavior seems to occur for the fractional BBM equation (fBBM)
\begin{equation}\label{fBBM}
u_t+u_x+uu_x+D^\alpha u_t=0
\end{equation}
in the subcritical case $\alpha>\frac{1}{3}.$ Again the simulations in  \cite{KS}) suggest that the soliton resolution also holds (see Figure \ref{fBBM05} below).

\begin{figure}[htb!]
 \includegraphics[width=\textwidth]{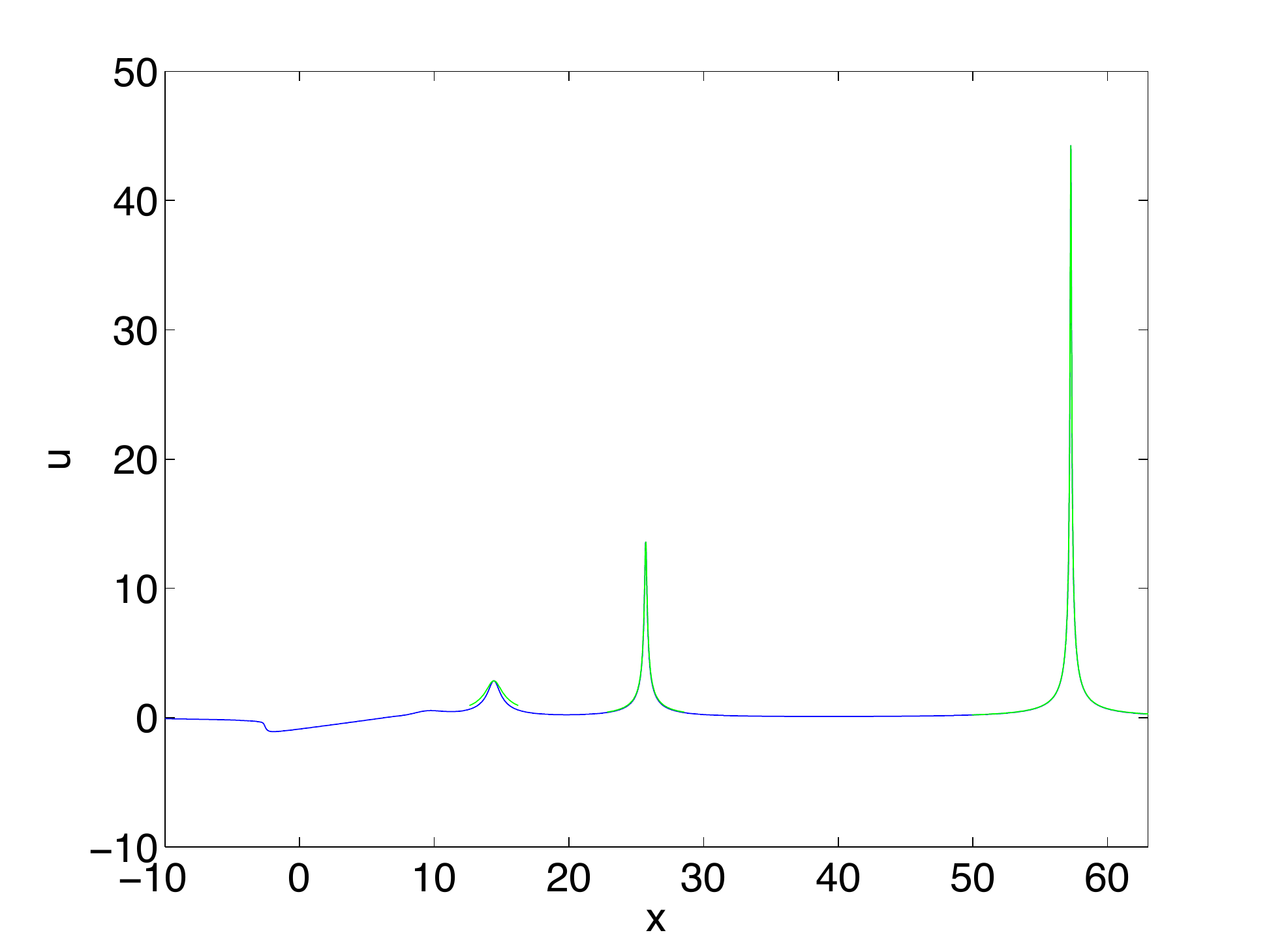}
 \caption{Solution to the fBBM equation  
 \eqref{fBBM} for $\alpha=0.5$, and the initial 
 data $\u_{0}=10\text{sech}^2 x$ for $t=10$ in blue, fitted solitons at the humps in green. }
 \label{fBBM05}
\end{figure}

\subsection {The Benjamin-Ono and Intermediate Long Wave Equations}

The  Intermediate Long Wave Equation (ILW) and the Benjamin-Ono equation (BO) are asymptotic models in an appropriate regime for a two-fluid system   when the depth $\delta$ of the  bottom layer is very large with respect to the upper one (ILW) or infinite (BO) (see \cite{BLS}).

The ILW corresponds (in the notations of \cite{ABFS}) to

$$p(\xi)=2\pi \xi \text{coth}\;(2\pi \delta \xi)-\frac{1}{\delta}$$

and the BO equation to

$$p(\xi)=2\pi |\xi|.$$

Alternatively they can be respectively written in convolution form

\begin{equation}\label{ILW}
u_t+uu_x+\frac{1}{\delta}u_x+\mathcal T_\delta(u_{xx})=0,
\end{equation}

$$\mathcal T_\delta(u)(x)=-\frac{1}{2\delta}PV \int_{-\infty}^{\infty} \text{coth}\left(\frac{\pi(x-y)}{2\delta}\right) u(y)dy$$

and

\begin{equation}\label{BO}
u_t+uu_x+\mathcal Hu_{xx}=0
\end{equation}

where $\mathcal H$ is the Hilbert transform

$$\mathcal Hu(x)=\frac{1}{\pi}PV \int_{-\infty}^\infty \frac{u(y)}{x-y}dy.$$

\subsubsection{The Benjamin-Ono equation} 

A striking difference between KdV and BO equations is that the latter 
is {\it quasilinear} rather that {\it semilinear}. This means that 
the Cauchy problem for BO cannot be solved by a Picard iterative 
scheme implemented on the integral Duhamel formulation, for  initial 
data in any Sobolev spaces $H^s(\R)$, $s\in \R$.  Alternatively, this implies that the flow map $u_0\mapsto u(t)$ cannot be 
{\it smooth} in the same spaces  \cite{MST4}, and actually not even locally Lipschitz \cite{KT2}. We will give a 
precise statements of those facts later on for the KP I equation which also is quasilinear in this sense.

The Cauchy problem has been proven to be globally well-posed in 
$H^s(\R), s>3/2$ by a compactness method using the various invariants\footnote{The existence of an infinite sequence of invariants (\cite{Ca, Mat3}) is of course a consequence of the complete integrability of the Benjamin-Ono equation.} of the equations (\cite {ABFS}) and actually in much bigger spaces (see \cite{Tao2, IK} and the references therein), in particular in the energy space $H^{1/2}(\R),$ by sophisticated methods based  on the dispersive properties of the equations.

Moreover it was proven in \cite{ABFS} that the solution $u_\delta$ of \eqref{ILW} with initial data $u_0$ converges as $\delta \to +\infty$ to the solution of the Benjamin-Ono equation \eqref{BO} with the same initial data.

Furthermore, if $u_\delta$ is a solution of \eqref{ILW} and setting 

$$v_\delta(x,t)=\frac{3}{\delta}u_\delta(x,\frac{3}{\delta}t),$$

$v_\delta$ tends as $\delta \to 0$ to the solution $u$
of the KdV equation

\begin{equation}\label{KdVnew}
u_t+uu_x+u_{xxx}=0
\end{equation}

with the same initial data.

\begin{remark}\label{HOBO}
The hierarchy of  conserved quantities of the Benjamin-Ono equation 
leads to a hierarchy of higher order BO equations $BO_n$ (by 
considering the associated Hamiltonian flow). Those equations have 
order $\frac{n}{2}$, $n=4,5, ...$ It was established in \cite{MP} that a family of order 3 equations containing $BO_6$ is globally well-posed in the energy space $H^1(\R).$ A similar result is expected for the whole hierarchy.
\end{remark}

Both the BO and the ILW equations  are classical examples of equations solvable par IST methods. The situation is however less satisfactory than for the KdV equation since the resolution of the Cauchy problem for the BO equation for instance needs a smallness condition on the initial data.

 The formal IST theory has been given by Ablowitz and Fokas \cite{AF}. They found the inverse spectral problem and Beals and Coifman \cite{BC1, BC2} observed that it is equivalent to a nonlocal $\bar{\partial}$ problem.

Unfortunately, the direct scattering problem can be only solved for 
small data and a complete theory for IST (as for the KdV equation) is a challenging open problem. In particular, one does not know (but expects) that any localized initial data decompose into a train of solitary waves and a dispersive tail. The rigorous theory of the Cauchy problem for small initial data is given in \cite{CW}.

As previously recalled, the BO equation has an infinite number of conserved quantities (\cite{Ca}, the first ones are displayed in \cite{ABFS}). The Hamiltonian flow of those invariants define the aforementioned Benjamin-Ono hierarchy. 

The BO equation possesses explicit soliton and multi-solitons 
(\cite{Ca2, Mat, Mat3}). The one soliton reads 
$$ Q(x)=\frac{4}{1+x^2}$$
and  is unique (up to translations) among all solitary wave solutions (\cite{AmTo}). Its slow (algebraic) decay is due (by Paley-Wiener type arguments) to the fact that the BO symbol $\xi|\xi|$ is not smooth at the origin. Actually similar arguments imply that the solution to the Cauchy problem cannot decay fast at infinity (see \cite{Io}). The BO solitary wave is orbitally stable (see \cite{BoSo, BSS} and the references therein).

On the other hand Kenig and Martel \cite{KeMa} have proven the 
asymptotic stability of the  BO solitary wave as well as that of the 
explicit multi-solitons in the energy space, a fact which reinforces 
the above conjecture on the long time dynamic of BO solutions. They  
do not use the integrability of the equation except for the exact 
expressions for the solitons (that help to study the spectral properties of the linearized operators).

\begin{remark}
The $H^1$ stability of the 2-soliton has been proven in \cite{NL} by variational methods.
\end{remark}

\begin{figure}[htb!]
 \includegraphics[width=\textwidth]{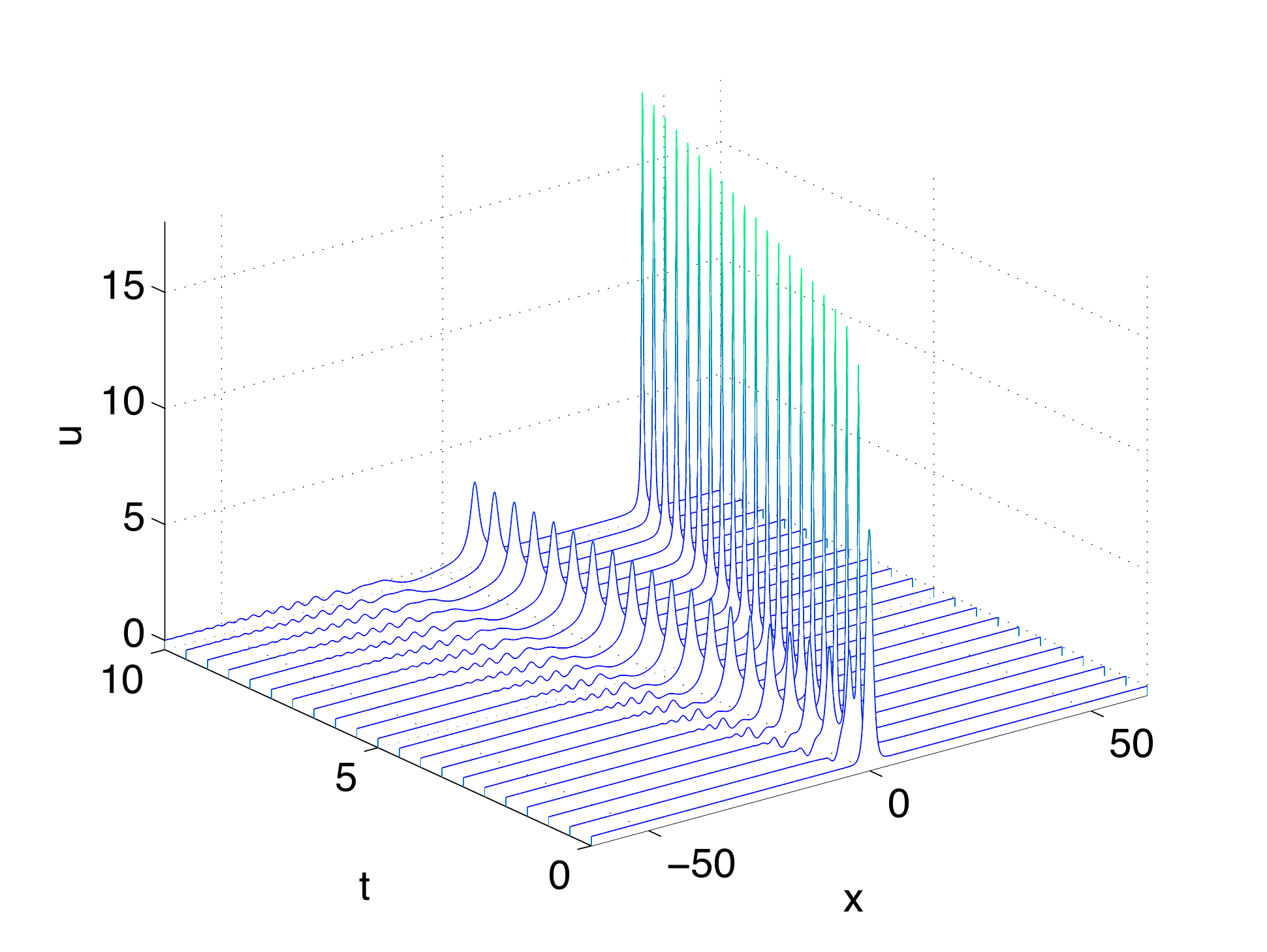}
 \caption{Solution to the BO equation  
for the initial 
 data $\u_{0}=10\text{sech}^2 x$ in dependence of $t$. }
 \label{BO10sech}
\end{figure}

We show the formation of solitons from localized initial data in 
Fig.~\ref{BO10sech}. Again there is a tail of dispersive oscillations 
propagating to the left. 

\subsubsection{The ILW equation}
 The formal IST theory has been given in \cite{KAS, KSA}. The direct 
 scattering problem is associated with a Riemann-Hilbert problem in a strip of the complex plane. As the BO equation, the ILW equation possesses an infinite sequence of conserved quantities (see {\it eg} \cite{LR}) which leads to a ILW hierarchy. They can be used to provide the global well-posedness of the Cauchy problem in $H^s(\R), s>3/2$, (\cite {ABFS}).
 
 Moreover,  the fact that for large $|\xi|$, 
 $|\xi|\text{coth}\;\xi=|\xi|(1+O(e^{-|\xi|}))$ implies that the well-posedness results are similar to those obtained for the BO equation, for instance in the energy space $H^{1/2}(\R).$

 On the other hand, we do not know of rigorous results for the Cauchy problem using the IST method. It is likely that they would require a smallness condition on the initial data.

Explicit N-soliton solutions are given in \cite{Mat, JE}. Contrary to 
those of the Benjamin-Ono equation, they decay exponentially at 
infinity (the ILW symbol is smooth). For instance, when the ILW equation is written in the form

\begin{equation}\label{ILW2}
u_t+2uu_x+\frac{1}{\delta}u_x-L_\delta u_x=0,
\end{equation}

where 

$$\widehat{L_\delta u}(\xi)=(\xi \text{coth}\;\xi \delta)\hat{u}(\xi),$$

 the 1-soliton  reads (see \cite{AlTo, J,JE}

$$u(x,t)=Q_{c,\delta}(x-ct),\quad Q_{c,\delta}(x)=\frac{a\sin a\delta}{\text{cosh}\; ax+\cos a\delta},\quad x\in \R$$

for arbitrary $c>0$ and $\delta >0,$ and $a$ is the unique positive solution of the transcendental equation

$$a\delta\; \text{coth}\; a\delta =(1-c\delta),\quad a\in 
(0,\pi/\delta).$$

Its uniqueness (up to translations) is proven in \cite{AlTo}. The orbital stability of this soliton is proven in \cite{AB1, ABH} (see also \cite{BoSo, BSS}). We do not know of asymptotic stability results for the 1 or N-soliton
similar to the corresponding ones for the BO equation. Those results should be in some sense easier than the corresponding ones for BO since the exponential decay of the solitons should make the spectral analysis of the linearized operators  easier.

\begin{figure}[htb!]
 \includegraphics[width=\textwidth]{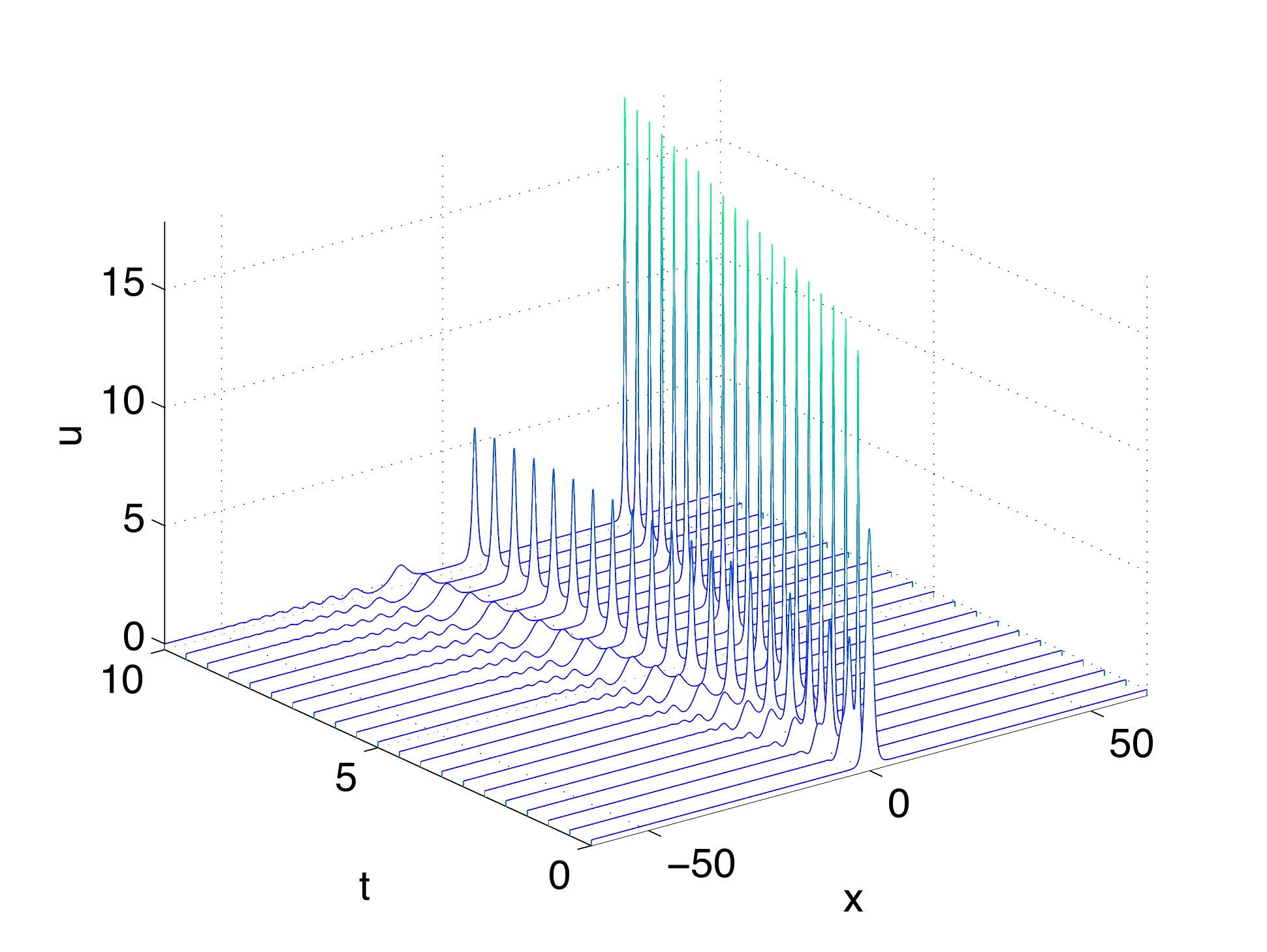}
 \caption{Solution to the ILW equation with $\delta=1$ 
for the initial 
 data $\u_{0}=10\text{sech}^2 x$ in dependence of $t$. }
 \label{IWL10sech}
\end{figure}

\begin{figure}[htb!]
 \includegraphics[width=\textwidth]{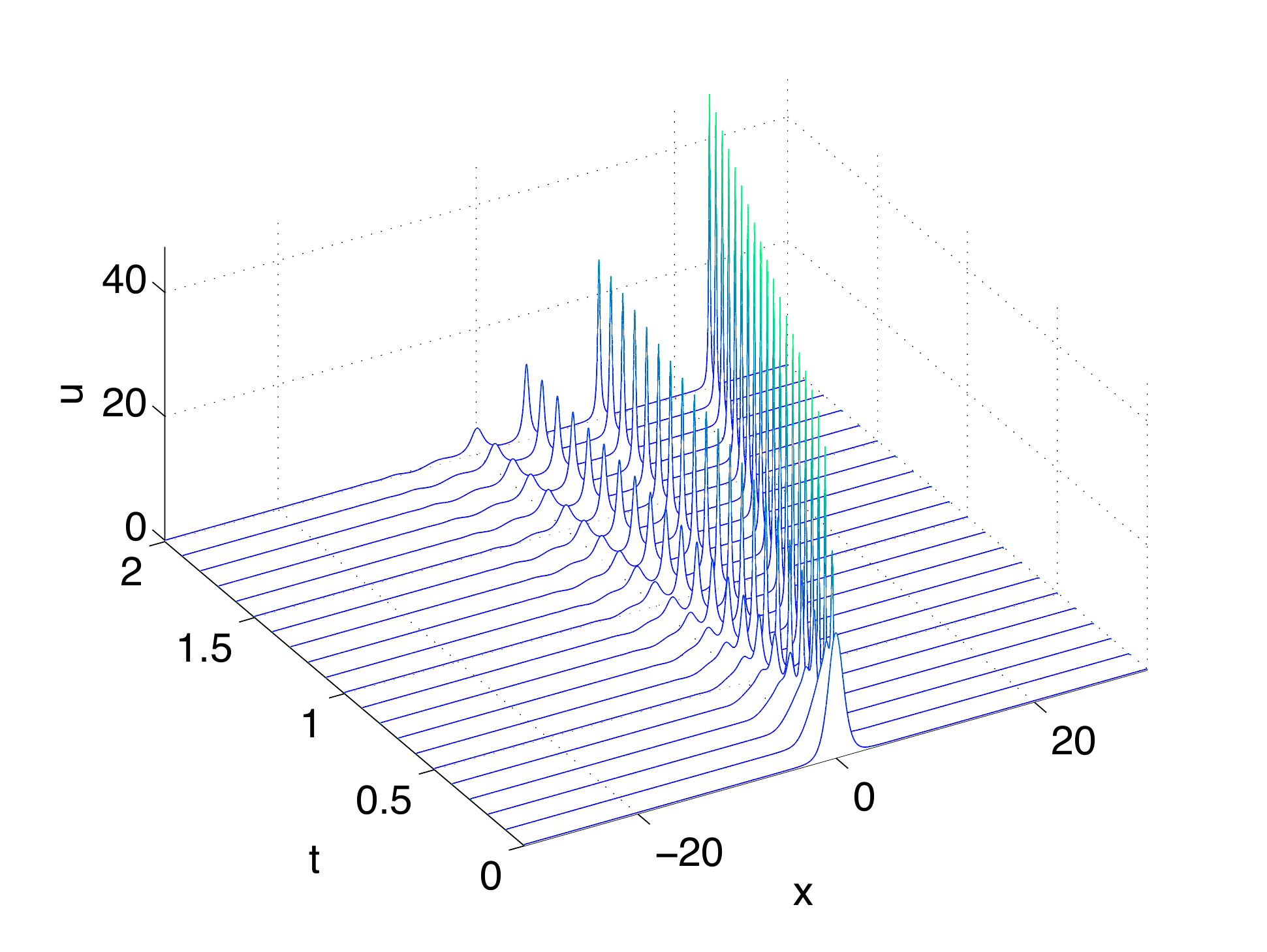}
 \caption{Solution to the ILW equation with $\delta=1$ 
for the initial 
 data $\u_{0}=20\text{sech}^2 x$ in dependence of $t$. }
 \label{IWL20sech}
\end{figure}

We show the decomposition of localized initial data into solitons and 
radiation in Fig.~\ref{IWL10sech} and  Fig.~\ref{IWL20sech}. Note that this case is numerically 
easier to treat with Fourier methods since the soliton solutions are 
more rapidly decreasing (exponentially instead of algebraically) than for the fKdV, fBBM and BO equations 
before. The different shape of the solitons is also noticeable in 
comparison to Fig.~\ref{BO10sech}. 

One also sees the change in the number of emerging solitons according to the size of the initial data. As in the case of BO, fKdV, fBBM equations, predicting the number of solitary waves which seem to form from a given (smooth and localized) initial data is a challenging open question.
 Except for 
the KdV case where the solitons are related to the discrete spectrum 
of the associated Schr\"odinger equation, there does not appear to be 
a clear characterization of the solitons emerging from given initial 
data for $t\to\infty$ even for integrable equations. 

\subsection{The Gross-Pitaevskii equation}

The Gross-Pitaevskii equation (GP) is a version of a nonlinear Schr\"{o}dinger equation (NLS), namely

\begin{equation}\label{GP}
i \partial_t \Psi = \Delta \Psi + \Psi (1 - |\Psi|^2) \quad\text{on} \quad \R^d \times \R,( d=1, 2, 3).
\end{equation}
It is a relevant model in nonlinear optics ("dark" and "black" solitons, optical vortices), fluid mechanics (superfluidity of Helium II), Bose-Einstein condensation of ultra cold atomic gases.

 At least on a formal level, the Gross-Pitaevskii equation is hamiltonian. The conserved Hamiltonian is a  {\bf Ginzburg-Landau energy}, namely
\begin{equation}\label{HGP}
E(\Psi) = \frac{1}{2} \int_{\R^d} |\nabla \Psi|^2 + \frac{1}{4} \int_{\R^d} (1 - |\Psi|^2)^2 \equiv \int_{\R^d} e(\Psi),
\end{equation}

associated to the natural energy space

$$\boE(\R^d) = \{ v \in H^1_{\text{loc}}(\R^d),  E(v) < + \infty \},$$

In order for $E(\Psi)$ to be finite,  $|\Psi|$ should in some sense tend to 1 at infinity. Actually this "non trivial" boundary condition provides (GP) with a richer dynamics  than in the case of null condition at infinity which, for a defocusing NLS type equation, is essentially governed by dispersion and scattering.

\vspace{1 cm}
For instance, in nonlinear optics, the "dark solitons" are localized nonlinear waves (or "holes")  which exist on a stable continuous wave background. The boundary condition  $|\Psi(x, \cdot)| \to 1$ is due to this non-zero background.

Similarly to the energy, the \textbf {momentum}
$\vec{P}(\Psi) = \frac{1}{2} \int_{\R^N} \langle i \nabla \Psi \ , \Psi  \rangle,$
is formally conserved. 

We will denote by $p=Im \int_{\R^N}i\Psi \bar{\Psi}_x$, the first component of $\vec{P}$, which is hence a scalar.

Justifying the momentum (and its conservation) is one of the difficulties one has to face when dealing with  the (GP) equation.  

We will restrict here to the one-dimensional case, $d=1$ since the GP equation is then completely integrable \cite{ZS}. We just consider

\begin{equation}\label{GP1D}
i\partial_t u +\partial^2_x u + (1-|u|^2)u = 0
\end{equation}
In fact Zakharov-Shabat (1973) consider the case when $|u (x,t)|^2 \rightarrow c >0$, $|x| \rightarrow \infty$ (propagation of waves through a condensate of constant density).

More precisely, the GP has a Lax pair $(B_u, L_u)$, where (for $c=1$)
\begin{equation}
L_u = i\left(\begin{array}{cc}
                       1+\sqrt{3}    & 0 \\
                       0                   & 1-\sqrt{3}\\
           \end{array}
   \right) \partial_x  +\left(\begin{array}{cc}
                      0                   & \bar{u}\\
                      u                   & 0 \\
                      \end{array}
      \right)
   \end{equation}

\begin{equation}
B_u = -\sqrt{3}\left(\begin{array}{cc}
                       1    & 0 \\
                       0                   & 1\\
           \end{array}
   \right) \partial^2_x  +\left(\begin{array}{cc}
                      \frac{|u|^2 -1}{\sqrt{3} + 1}                   &i \partial_x \bar{u}\\
                       - i \partial_x \bar{u}\                 & \frac{|u|^2 -1}{\sqrt{3} - 1}    \\
                      \end{array}
      \right)
   \end{equation}

So that u satisfies (GP) if and only if
\begin{equation}
\frac{d}{dt}L_u = [ L_u, B_u]
\end{equation}

 As a consequence, the 1 D Gross-Pitaevskii equation has an infinite number of (formally)  conserved energies $E_k, k\in \N$ and momentum $P_k, k\in \N.$

For instance,
$$I_5=\int_{\R}\lbrace 2|u|^6+6|u|^2|u_x|^2+(\frac{d}{dx}|u|^2)^2 +|u_{xx}|^2-2\rbrace.$$

It is of course necessary to prove rigorously that the $E_k$ and the $P_k$ are well defined and conserved by the GP flow, in a suitable functional setting.

We do not know of a {\it complete} resolution of the Cauchy problem by IST methods, including a possible decomposition into solitons. We will see however that one can get a proof of stability of the solitons by using IST techniques.

\vspace{1cm}

One can prove that the energy space in one dimension is identical to the {\it Zhidkov space} (see \cite{Z})

  $$X^1(\R) = \lbrace u\in H^{1}_{loc}(\R), u_x \in L^\infty(\R), (1-|u|^2) \in L^2(\R)\rbrace$$
  
  and actually Zhidkov \cite{Z} proved that the Cauchy problem for 
  \eqref{GP1D} is globally well-posed in $X^1(\R)$ \footnote {P. G\' erard has proven in fact  (\cite{Ge1}) that the Cauchy problem for the GP equation is globally well -posed in the energy space $\boE(\R^d)$ equipped with a suitable topology when $d\leq 4.$}.

The one-dimensional GP equation possesses two types of solitary waves of velocity $c$,  $0 \leq c < \sqrt{2}$ :
  \vspace{1 cm}
  
  \begin{itemize}
  \item The {\bf "dark"} solitons : 
  $ \v_c(x) \equiv \sqrt{\frac{2 - c^2}{2}} \th \Big( \frac{\sqrt{2 - c^2}}{2} x \Big) - i \frac{c}{\sqrt{2}}$.
\item The { \bf"black"} soliton : $$\v_0(x) = \th \Big( \frac{x}{\sqrt{2}} \Big).$$
\end{itemize}

Note that when $0< c < \sqrt{2},\:  \v_c(x)\not= 0, \forall x$ while the black soliton vanishes at  $x=0$.

  The orbital stability of the dark soliton has been proven in  \cite{BGS3} (see also \cite{Lin} for the cubic-quintic case).  

This case is easier since the dark soliton does not vanish and the momentum can be defined in a straightforward way.
 The orbital stability of the black soliton is more delicate since it 
 vanishes at $0.$ Both PDE and IST techniques provide the result, in a slightly different form though, and this is a good opportunity to compare them.

 The first method is used in \cite{BGSS3}. This is the  
 "Hamiltonian" method (see \cite{B, Jerry} for the stability of the 
 KdV solitary wave or \cite{CZ} for the stability of the focusing NLS ground states), that is one considers the  black soliton as a minimizer of the energy with fixed momentum.
 As previously noticed, serious difficulties arise here from the momentum (definition, conservation by the flow,...) because the black soliton vanishes at $0$.
  
   
    
    Given any $A > 0$, we consider on the energy space $X^1(\R)=X^1$ the distance $ d_{A, X^1}$ defined by
$$d_{A, X^1}(v_1, v_2) \equiv \| v_1 - v_2 \|_{L^\infty([- A, A])} + \| v_1' - v_2' \|_{L^2(\R)} + \| |v_1| - |v_2| \|_{L^2(\R)}.$$

 One can show that for $v\in X^1,$ $|v(x)|\rightarrow 1$ as $|x|\rightarrow \infty,$ (but not $v$ itself).

We have the orbital stability result (\cite{BGSS3}):
  
     \begin{theorem}\label{stab}
\label{stationnaire}
Assume that $v_0 \in X^1$ and consider the global in time solution
$v$ to (GP) with initial datum $v_0$. Given any numbers $\varepsilon > 0$ and
$A > 0$, there exists some positive number $\delta$, such that if
\begin{equation}
\label{pioneer}
d_{A, X^1}(v_0, \v_0) \leq \delta,
\end{equation}
then, for any $t \in \R$, there exist numbers $a(t)$ and $\theta(t)$ such that
\begin{equation}
\label{solaris}
d_{A, X^1} \big( v(\cdot + a(t), t), \exp (i \theta (t)) \ \v_0(\cdot) \big) < \varepsilon.
\end{equation}
\end{theorem}

  We also have a control of the shift $a(t)$:

\begin{theorem} \label{shift}
\label{vitesse}
Given any numbers $\varepsilon > 0$, sufficiently small, and $A > 0$, there exists some constant $K$, only depending on $A$, and some positive number $\delta > 0$ such that, if $v_0$ and $v$ are as in the previous Theorem, then
\begin{equation}
\label{eqshift}
|a(t)| \leq K \varepsilon (1 + |t|),
\end{equation}
for any $t \in \R$, and for any of the points $a(t)$ as above.
\end{theorem}

  \vspace{0.3cm}

We will not give a proof of the previous results (see \cite{BGSS3}) but just explain how to extend the definition of the momentum $P(v)=\frac{1}{2}\int_{\R}\langle iv,v'\rangle$
for general functions in $X^1.$

    \begin{itemize}
    \item Assume $v\in X^1$ and $v=\exp(i\phi).$
    \end{itemize}
    
    Then $\langle iv,v'\rangle=\phi ' \in L^2(\R)$ and $P(v)\equiv \frac{1}{2}\lbrack \phi(+\infty)-\phi(-\infty)\rbrack.$
    
    This is meaningful if $$v\in Z^1 = \Big\{ v \in X^1, \ {\rm s.t.} \ v_{\pm \infty} = \lim_{x \to \pm \infty} v(x) \ {\rm exist} \Big\},$$ 
but not for any arbitrary phase $\phi$ whose gradient is in $L^2$.

Let $$\tilde{X^1} = \{ v \in X^1, \ {\rm s.t.} \ |v(x)| >0, \ \forall x \in \R \}.$$
For $v\in \tilde{X^1}$, $v = \rho \exp (i \phi)$, so that $\langle i v, v' \rangle=\rho^2\phi'$. If $v \in \tilde{Z^1} \equiv \tilde{X^1} \cap Z^1$, one has 
\begin{equation}
\label{goret}
P(v) = \frac{1}{2} \int_\R \rho^2 \phi' = \frac{1}{2} \int_\R (\rho^2 - 1) \phi' + \frac{1}{2} \int_\R \phi' = \frac{1}{2} \int_\R (\rho^2 - 1) \phi'+ \frac{1}{2} \big[ \phi \big]_{- \infty}^{+ \infty}.
 \end{equation} 
 This is well controlled since by Cauchy-Schwarz
 $$|\int_{\R}(\rho^2-1)\phi^{\prime}|\leq \frac{2}{\delta^2}E(v)$$
 where
 $$\delta =inf\lbrace|v(x)|, x\in \R \rbrace.$$

        \begin{itemize}
        \item The momentum for maps with zeroes:
        \end{itemize}
        
        \begin{lemma}
\label{limite}
Let $v \in Z^1$. Then, the limit 
$$\boP(v) = \underset{R \to + \infty}{\lim} P_R(v) \equiv \lim_{R\to +\infty} \int_{- R}^{R} \langle i v, v' \rangle$$
exists. Moreover, if $v$ belongs to $\tilde{Z^1}$, then
$$\boP(v) = \frac{1}{2} \int_\R (\rho^2 - 1) \phi'+ \frac{1}{2} \big[ \phi \big]_{- \infty}^{+ \infty}.$$
\end{lemma}

\begin{itemize}
\item Example the black soliton  $\v_0$: since $\v_0$ is real-valued, $\langle i\v_0,\v_0' \rangle =0$ and $\boP(\v_0)=0.$
\end{itemize}

 \begin{itemize}
 \item Elementary observation. Let $V_0 \in Z^1$ and $w\in H^1(\R).$ Then, $V_0+w\in Z^1$ and
 
 $$\boP(V_0+w)=\boP(V_0)+\frac{1}{2}\int_{\R}\langle iw,w'\rangle +\int_{\R}\langle iw,V_0'\rangle.$$

 \item Note that besides $\boP(V_0)$ the integrals are definite ones.
  \end{itemize}
        
\vspace{0.3cm}
\begin{itemize}
\item {\bf The renormalized and the untwisted momentum.}
\end{itemize}

The renormalized momentum $p$ is defined for $v \in \tilde{X_1}$ by
\begin{equation}
\label{defnormmo}
p(v) = \frac{1}{2} \int_\R (\varrho^2 - 1) \varphi',
\end{equation}
 as seen before, if $v$ belongs to $\tilde{Z^1}$, then,
$$p(v) = \boP(v) - \frac{1}{2} \big[ \varphi \big]_{- \infty}^{+ \infty}.$$
If $v \in Z^1 \setminus \tilde{Z^1}$, the integral is a priori not well-defined since the phase $\varphi$ is not globally defined. Nevertheless, the argument $\arg{v}$ of $v$ is well-defined at infinity as an element of  $\R/ 2\pi\Z$. For $v \in Z^1$, one introduces \textbf{the untwisted momentum}
$$[p](v) = \Big( \boP(v) - \frac{1}{2} \big( \arg{v(+ \infty)} - \arg{v(- \infty)} \big) \Big) \ {\rm mod} \ \pi,$$
which is hence an element of $\R/ \pi\Z$. A remarkable fact concerning $[p]$ is that its definition extends to the whole space $X^1$, although for arbitrary maps in $X^1$, the quantity $\arg{v(+\infty)}-\arg{v(-\infty)}$ may not exist. 

\begin{lemma}
\label{limite2}
Assume that $v$ belongs to $X^1$. Then the limit
$$[p](v) = \lim_{R \to +\infty} \bigg( \int_{- R}^R \langle i v, v' \rangle - \frac{1}{2} \big( \arg{v(R)} - \arg{v(-R)} \big) \bigg) \ {\rm mod} \ \pi$$
exists. Moreover, if $v$ belongs to $\tilde{X^1}$, then
\begin{equation}
\label{bacon2}
[p](v) = p(v) \ {\rm mod} \ \pi.
\end{equation}
\end{lemma}
One has also:
\begin{itemize}
\item Let $V_0 \in X^1$ and $w \in H^1(\R)$. Then, $V_0 + w \in X^1$ and
\begin{equation}
\label{eq:andouille2}
[p](V_0 + w) = [p](V_0) + \frac{1}{2} \int_\R \langle i w, w' \rangle + \int_\R \langle i w, V_0' \rangle \ {\rm mod} \ \pi.
\end{equation}

\end{itemize}

\begin{itemize}
\item The evolution preserves the untwisted momentum:
\end{itemize}

\begin{lemma}
\label{grouic}
Assume $v_0 \in X^1$, and let $v$ be the solution to (GP) with initial datum $v_0$. Then,
$$[p](v(\cdot, t)) = [p](v_0), \ \forall t \in \R.$$
If moreover $v_0 \in Z^1$, then $v(t)$ belongs to $Z^1$ for any $t \in \R$, and
$$\boP(v(\cdot, t)) = \boP(v_0), \ \forall t \in \R.$$
\end{lemma}

The orbital stability of the dark soliton ($0<c<\sqrt{2})$ is based on the minimization problem:

$$E_{\min}(\p) = \inf\{ E(v), v \in \tilde{X^1} \; \text{s.t.} \; 
p(v) = \p \}.$$

For map having zeroes like $\v_0$, one should use the untwisted momentum and consider instead  
  \begin{equation}
\label{porcelet}
\mathfrak{E}_{\min} \Big( \frac{\pi}{2} \Big) \equiv \inf \big\{ E(v), v \in X^1 \ \text{s.t.}\;  [p](v) = \frac{\pi}{2} \; \text{ mod} \; \pi \big\}.
\end{equation}


\vspace{0.5cm}

We now turn to the second approach using IST in \cite{GZ} which avoids the renormalization by factors of modulus 1 in Theorem \ref{stationnaire}, at least
for sufficiently smooth and decaying perturbations. 

\begin{theorem}\label{GZ}

 Assume that the initial datum of (1.1) has the form:
$u_0(x) = U_0(x) +\epsilon u_1(x), U_0(x) = \text{tanh}\left(\frac{x}{\sqrt 2}\right),$ where $u_1(x)$ satisfies the following condition,
 
 $$\sup_{x\in \R}|(1+ x^2)^2\partial^ku_1(x)|\leq 1,\quad \text{for} \; k\leq 3.$$

 $$\forall t\in \R,\quad \exists \;y(t)\in \R,\quad |\tau_{y(t)}u(.,t)-U_0|_\infty \leq C \epsilon, \quad \text{for} \; 0\leq t<+\infty.$$
\end{theorem}
Using a classical functional analytic argument,  Theorem \ref{GZ} easily yields the orbital stability, at least for
sufficiently smooth and decaying perturbations \cite{GZ}:

\begin{corollary}\label{GZcor}
For every $\delta>0$ there exists $\epsilon >0$ such that, if 

 $$\sup_{x\in \R}|(1+ x^2)^2\partial^ku_1(x)|\leq 1,\quad \text{for} \; k\leq 3,$$
 
 then the solution $u$ of \eqref{GP1D} satisfies:
 
 $$\forall t\in \R,\; \exists y(t)\in \R,\quad d_E(\tau_{y(t)}u(t),U_0)\leq \delta.$$

\end{corollary}

\begin{remark}
As previously noticed Theorem \ref{GZ} improves on Theorem \ref{stab}, \ref{shift}  since it does not require a rotation factor. On the other hand it deals with a narrower class of perturbations and does not provide an evaluation of the time shift.
\end{remark}

One can combine PDE and IST techniques as for instance  in the long wave limit of the GP equation which we describe now.

The transonic (long wave, small amplitude) limit of KGP to KP I in 2D or KdV in 1D  is quite a generic phenomenon for nonlinear dispersive systems. A good example is that of the water wave systems (incompressible Euler with upper free surface).

While the {\it consistency} of the approximation is relatively easy to obtain, the {\it stability} (and thus the {\it convergence}) of the approximation is much more delicate, specially if one looks for the optimal error estimates on the correct time scales.


Kuznetsov and Zakharov \cite{ZK} have observed formally that the KdV equation provides  a good approximation of long-wave small amplitude solutions to the 1D Gross-Pitaevskii equation. We recall here briefly a rigorous proof of this fact \cite{BGSS1} (see also \cite{CR} for a different approach which dose not provide an error estimate and \cite{BGSS2} for a similar analysis for two-way propagation leading to a coupled system of KdV equations). We thus start from the 1D GP equation

\begin{equation}
\label{GP1}
i\partial_t\Psi+\partial_x^2\Psi=\Psi(|\Psi|^2-1)\quad \text{on}\quad \R\times \R, \quad \Psi(.,0)=\Psi_0,
\end{equation}

\begin{equation}\label{lim}
|\Psi(x, t)| \to 1, \text{as} \quad |x| \to + \infty,
\end{equation}
and recall the formal conserved quantities

$$E(\Psi) = \frac{1}{2} \int_{\R} |\partial_x \Psi|^2 + \frac{1}{4} \int_{\R} (1 - |\Psi|^2)^2 \equiv \int_{\R} e(\Psi).$$
$$P(\Psi) = \frac{1}{2} \int_{\R} \langle i \partial_x \Psi, \Psi \rangle.$$

One can prove  (see for instance  \cite{Ga2}) well-posedness in Zhidkov spaces  :

$$X^k(\R) = \{ u \in L^1_{\rm loc}(\R, \C), \ {\rm s.t.} \ E(u) < + \infty, \ {\rm and} \ \partial_x u \in H^{k - 1}(\R) \},$$

\begin{theorem}
\label{thm:existe}
Let $k \in \N^*$ and $\Psi_0 \in X^k(\R)$. Then, there exists a unique solution $\Psi(\cdot, t)$ in $\boC^0(\R, X^k(\R))$ to \eqref{GP1} with initial data $\Psi_0$. If $\Psi_0$ belongs to $X^{k+2}(\R)$, then the map $t \mapsto \Psi(\cdot, t)$ belongs to $\boC^1(\R, X^k(\R))$ and $\boC^0(\R, X^{k+2}(\R))$. Moreover, the flow map $\Psi_0 \mapsto \Psi(\cdot, T)$ is continuous on $X^k(\R)$ for any fixed $T \in \R$.
\end{theorem}
\begin{remark}
One can prove  conservation of energy and momentum (under suitable assumptions).
\end{remark}

We recall that if $\Psi$ does not vanish, one may write (Madelung transform)
$$\Psi = \sqrt{\rho} \exp i \varphi.$$
This leads to the hydrodynamic form of the equation, with $v = 2 \partial_x \varphi,$
\begin{equation}
\label{HDGP}
\left\{ \begin{array}{ll} \partial_t \rho + \partial_x (\rho v) = 0,\\
\rho (\partial_t v + v.\partial_x v) + \partial_x (\rho^2) = \rho 
\partial_x \Big( \frac{ \partial_x^2 \rho}{\rho} - \frac{|\partial_x 
\rho|^2}{2 \rho^2} \Big), \end{array} \right.
\end{equation}

which can be seen as a compressible Euler system with pressure law $ p(\rho)=\rho ^2$ and a quantum pressure term.

It is shown in \cite{BDS} that \eqref {GP1} or \eqref{HDGP}  can be  approximated by the {\it linear} wave equation. We justify here the long wave approximation on larger time scales $\mathcal{O}(\eps^{-3})$ and, following Kuznetsov and Zakharov,  introduce the slow variables :
\begin{equation}
\label{martinsepp}
X = \varepsilon(x + \sqrt{2} t), \ {\rm and} \ \tau = \frac{\varepsilon^3}{2 \sqrt{2}} t.
\end{equation}

 This corresponds to a reference frame traveling to the left with velocity $\sqrt{2}$ in the original variables $ (x, t)$. In this frame the left going wave is stationary while the right going wave has a speed $8\eps^{-2}$
and is appropriate to study waves traveling to the left (we need to impose additional assumptions which imply the smallness of the right going waves).

We  define the rescaled functions $N_\varepsilon$ and $\Theta_\varepsilon$ as follows
\begin{equation}
\label{slow-var}
\begin{split}
N_\varepsilon(X, \tau) & = \frac{6}{\varepsilon^2} \eta(x, t) = \frac{6}{\varepsilon^2} \eta \Big( \frac{X}{\varepsilon} - \frac{4 \tau}{\varepsilon^3}, \frac{2 \sqrt{2} \tau}{\varepsilon^3} \Big),\\
\Theta_\varepsilon(X, \tau) & = \frac{6 \sqrt{2}}{\varepsilon} \varphi(x, t) =
\frac{6 \sqrt{2}}{\varepsilon} \varphi \Big( \frac{X}{\varepsilon} - \frac{4
\tau}{\varepsilon^3}, \frac{2 \sqrt{2} \tau}{\varepsilon^3} \Big),
\end{split}
\end{equation}
where $\Psi = \varrho \exp (i \varphi)$ and $\eta = 1 - \varrho^2= 1 - |\Psi|^2$.

\begin{theorem}
\label{cochon}
\cite{BGSS1}
Let $\varepsilon > 0$ be given and assume that the initial data 
$\Psi_0(\cdot) = \Psi(\cdot, 0)$ belong to $X^4(\R)$ and satisfy the assumption
\begin{equation}
\label{grinzing1}
\| N_\varepsilon^0 \|_{H^3(\R)} + \eps \| \partial^4_x N_\varepsilon^0 \|_{L^2(\R)}+ \|\partial_x \Theta_\eps^0\|_{H^3(\R)} \leq K_0.
\end{equation}
Let $\boN_\varepsilon$ and $\boM_\eps$ denote the solutions to the Korteweg-de Vries equation
\begin{equation}
\label{KdVGP}
\partial_\tau N + \partial_x^3 N + N \partial_x N = 0
\end{equation}
with initial data $N_\varepsilon^0$, respectively $\partial_x \Theta_\eps^0$. There exist positive constants $\varepsilon_0$ and $K_1$, depending possibly only on $K_0$ such that, if $\varepsilon \leq \varepsilon_0$, we have for any $\tau \in \R$,
\renewcommand{\theequation}{\arabic{equation}}
\setcounter{equation}{18}
\begin{equation}
\begin{split}
\label{eq:fortis}
\| \boN_\varepsilon & (\cdot, \tau) - N_\varepsilon(\cdot, \tau ) \|_{L^2(\R)} + \| \boM_\eps(\cdot, \tau) - \partial_x \Theta_\varepsilon(\cdot, \tau ) \|_{L^2(\R)}\\
& \leq K_1 \big( \eps + \| N_\eps^0 - \partial_x \Theta_\eps^0 \|_{H^3(\R)} \big) \exp (K_1 |\tau|).
\end{split}
\end{equation}
\end{theorem}

\begin{itemize}
\item This is a convergence result to the KdV equation for appropriate {\it prepared} initial data.
\item Since the norms involved are translation invariant, the KdV approximation can be only valid if the right going waves are negligible. This is the role of the term $ \| N_\eps^0 - \partial_x \Theta_\varepsilon^0 \|_{H^3(\R)}.$
\item In particular, if $\| N_\eps^0 - \partial_x \Theta_\eps^0 \|_{H^3(\R)} \leq C \eps^\alpha$, with $\alpha > 0$, then the approximation is valid on a time interval $t \in [0, S_\eps]$ with $S_\eps = o(\eps^{- 3} |\log(\eps)|)$. Moreover, if $\| N_\eps^0 - \partial_x \Theta_\eps^0 \|_{H^3(\R)}$ is of order $\boO(\eps)$, then the approximation error remains of order $\boO(\eps)$ on a time interval $t \in [0, S_\eps]$ with $S_\eps = \boO(\eps^{- 3})$.
\end{itemize}

The functions $ N_{\eps}$ and $\partial_X \Theta_\varepsilon$ are rigidly constrained one to the other:

\begin{theorem}
\label{H3-controlbis}
Let $\Psi$ be a solution to GP in $\boC^0(\R, H^4(\R))$ with initial data $\Psi^0$. Assume that \eqref{grinzing1} holds. Then, there exists some positive constant $K$, which does not depend on $\varepsilon$ nor $\tau$, such that
\begin{equation}
\label{dobling1ter}
\| N_\varepsilon(\cdot, \tau) \pm \partial_X \Theta_\varepsilon(\cdot, \tau) \|_{L^2(\R)} \leq \| N_\varepsilon^0 \pm \partial_X \Theta_\varepsilon^0 \|_{L^2(\R)} + K \varepsilon ^2 \big( 1 + |\tau| \big),
*\end{equation}
for any $\tau \in \R$.
\end{theorem}

The approximation errors provided by the previous theorems diverge as time increases. Concerning the weaker notion of {\bf consistency}, we have the following result whose peculiarity is that the bounds are independent of time.

\begin{theorem}
\label{Bobby}
Let $\Psi$ be a solution to GP in $\boC^0(\R, H^4(\R))$ with initial data $\Psi^0$. Assume that  \eqref{grinzing1} holds. Then, there exists some positive constant $K$, which does not depend on $\varepsilon$ nor $\tau$, such that
\begin{equation}
\label{jerrard}
\| \partial_\tau U_\eps + \partial^3_X U_\eps + U_\eps \partial_X U_\eps \|_{L^2(\R)} \leq K(\eps + \|N_\eps^0 - \partial_X \Theta_\eps \|_{H^3(\R)}),
\end{equation}
for any $\tau \in \R$, where $U_\eps = \frac{N_\eps + \partial_x \Theta_\eps}{2}$.
\end{theorem}

\begin{itemize}
\item One gets explicit bounds for the traveling wave solutions
\end{itemize}

 $\Psi(x, t) = v_c(x + c t)$. Solutions  do exist for any value of the speed $c$ in the interval $[0, \sqrt{2})$. Next, we choose the wave-length parameter to be $\varepsilon = \sqrt{2 - c^2}$, and take as initial data $\Psi_\varepsilon$ the corresponding wave $v_c$. We consider the rescaled function
$$\nu_\varepsilon(X) = \frac{6}{\varepsilon^2} \eta_c \Big( \frac{x}{\varepsilon} \Big),$$
where $\eta_c \equiv 1 - |v_c|^2$. The explicit integration of the travelling wave equation for $v_c$ leads to the formula
$$\nu_\varepsilon(X) = \nu(x) \equiv \frac{3}{\ch^2 \big( \frac{X}{2} \big)}.$$
The function $\nu$ is the classical soliton to the Korteweg-de Vries equation, which is moved by the KdV flow with constant speed equal to $1$, so that
$$\boN_\varepsilon(X, \tau) = \nu(X - \tau).$$

On the other hand, we deduce from \eqref{slow-var} that $N^0_\varepsilon = \nu$, so that
$$N_\varepsilon(X, \tau) = \nu \Big( X -\frac{4}{\eps^2} \Big( 1 - \sqrt{1 - \frac{\eps^2}{2}} \Big) \tau \Big).$$
Therefore, we have for any $\tau\in \R$,
$$\| \boN_\varepsilon(\cdot, \tau) - N_\varepsilon(\cdot, \tau) \|_{L^2(\R)} = \boO(\eps^2 \tau).$$

This suggests that the $\eps$ error in the main theorem should be of 
order $\eps^2$. This is proved in \cite{BGSS2} at a cost of higher regularity on the initial data (and also for a two way propagation, described by a system of two KdV equations).

\vspace{0.3cm}
We now give some elements of the proofs which rely on energy methods.

 We first write the equations for $N_\varepsilon$ and
$\Theta_\varepsilon$:
\begin{equation}
\label{slow1-0}
\partial_X N_\varepsilon - \partial_x^2 \Theta_\varepsilon + \frac{\varepsilon^2}{2} \Big( \frac{1}{2} \partial_\tau N_\varepsilon + \frac{1}{3} N_\varepsilon \partial_X^2 \Theta_\varepsilon + \frac{1}{3} \partial_X N_\varepsilon \partial_X \Theta_\varepsilon \Big) = 0,
\end{equation}
and
\begin{equation}
\label{slow2-0}
\partial_X \Theta_\varepsilon - N_\varepsilon + \frac{\varepsilon^2}{2} \Big( \frac{1}{2} \partial_\tau \Theta_\varepsilon + \frac{\partial_X^2 N_\varepsilon}{1 - \frac{\varepsilon^2}{6} N_\varepsilon} + \frac{1}{6} (\partial_X \Theta_\varepsilon)^2 \Big) + \frac{\varepsilon^4}{24} \frac{(\partial_X N_\varepsilon)^2}{(1 - \frac{\varepsilon^2}{6} N_\varepsilon)^2} = 0.
\end{equation}

The leading order in this expansion is provided by $N_\varepsilon - \partial_X \Theta_\varepsilon$ and its spatial derivative, so that an important step is to keep control on this term. In view of \eqref{slow1-0} and \eqref{slow2-0} and d'Alembert decomposition, we are led to introduce the new variables $U_\varepsilon$ and $V_\varepsilon$ defined by
$$U_\varepsilon = \frac{N_\varepsilon + \partial_X \Theta_\varepsilon}{2}, \ {\rm and} \ V_\varepsilon = \frac{N_\varepsilon - \partial_X \Theta_\varepsilon}{2},$$
and compute the relevant equations for $U_\varepsilon$ and $V_\varepsilon$,
\begin{equation}
\label{slow1}
\partial_\tau U_\varepsilon + \partial_X^3 U_\varepsilon + U_\varepsilon \partial_X U_\varepsilon = - \partial_X^3 V_\varepsilon + \frac{1}{3}\partial_X \big( U_\varepsilon V_\varepsilon \big) + \frac{1}{6} \partial_X \big( V_\varepsilon^2 \big) - \varepsilon^2 R_{\varepsilon},
\end{equation}
and
\begin{equation}
\label{slow2}
\partial_\tau V_\varepsilon + \frac{8}{\varepsilon^2} \partial_x V_\varepsilon = \partial_X^3 U_\varepsilon + \partial_X^3 V_\eps + \frac{1}{2} \partial_x (V_\varepsilon^2) - \frac{1}{6} \partial_x (U_\varepsilon)^2 -\frac{1}{3} \partial_X (U_\varepsilon V_\varepsilon) + \varepsilon^2 R_{\varepsilon},
\end{equation}

where the remainder term $R_{\varepsilon}$ is given by the formulae
\begin{equation}
\label{grouin}
R_{\varepsilon} = \frac{N_\varepsilon\partial_X^3 N_\varepsilon}{6 (1 - \frac{\varepsilon^2}{6} N_\varepsilon)} + \frac{(\partial_X N_\varepsilon) (\partial_X^2 N_\varepsilon)}{3 (1 - \frac{\varepsilon^2}{6} N_\varepsilon)^2} + \frac{\varepsilon^2}{36} \frac{(\partial_X N_\varepsilon)^3}{(1 - \frac{\varepsilon^2}{6} N_\varepsilon)^3}.
\end{equation}

The main step is to show  that the RHS of the equation for $U_\eps$ is small in suitable norms. In particular one must show that $V_\eps$ which is assumed to be small at time $\tau =0$ remains small, and that $U_\eps$, assumed to be bounded at time $\tau=0$, remains bounded in appropriate norms.
We use in particular various conservation laws provided by the integrability of the one-dimensional Gross-Pitaevskii equation.

\begin{itemize}
\item For instance we use the conservation of momentum and energy to get the $L^2$ estimates.
\item It turns out that the other conservation laws behave as higher order energies and higher order momenta. We use them to get :
\end{itemize}
\begin{theorem}
\label{H3-control}
Let $\Psi$ be a solution to (GP) in $\boC^0(\R, H^4(\R))$ with initial data $\Psi^0$. Assume that \eqref{grinzing1} holds. Then, there exists some positive constant $K$, which does not depend on $\varepsilon$ nor $\tau$, such that
\begin{equation}
\label{grinzing1bis}
\| N_\varepsilon(\cdot, \tau) \|_{H^3(\R)} + \eps \| \partial_X^4 N_\varepsilon(\cdot, \tau) \|_{L^2(\R)} + \|\partial_X \Theta_\eps(\cdot, \tau) \|_{H^3(\R)} \leq K,
\end{equation}
and
\begin{equation}
\label{dobling1bis}
\| N_\varepsilon(\cdot, \tau) \pm \partial_X \Theta_\varepsilon(\cdot, \tau) \|_{H^3(\R)} \leq K \big( \| N_\varepsilon^0 \pm \partial_X \Theta_\varepsilon^0 \|_{H^3(\R)} + \eps \big),
\end{equation}
for any $\tau \in \R$.
\end{theorem}


 The proof of the previous result led to a number of facts linked to integrability which have independent interest:
 
 \begin{itemize}
\item It provides expressions for the invariant quantities of GP 
\footnote{They appear in pairs: generalized energies $\boE_k$ and 
generalized momenta $\boP_k$.} and prove that they are well-defined 
in the spaces $X^k(\R)$. Their expressions are not a straightforward consequence of the induction formula of Zakharov and Shabat since many renormalizations have to be performed to give them a sound mathematical meaning.
\item It stablishes rigorously that they are conserved by the GP flow in the appropriate functional spaces.
\item It displays a striking relationship between the conserved quantities of the Gross-Pitaevskii equation and the KdV invariants:
$$\boE_k(N, \partial_X \Theta) - \sqrt{2} \boP_k(N, \partial_x \Theta) = E_k^{KdV} \Big( \frac{N - \partial_X \Theta}{2} \Big) + \boO(\eps^2).$$
\end{itemize}

It would be interesting to investigate further connections between the IST theories of the KdV and GP equations.

\begin{remark}\label{ISTGP}
We do not know of a rigorous result by IST methods describing the qualitative behavior of a solution (solitons+radiation,...) of the GP equation corresponding to a smooth and localized initial data.
\end{remark}

\section{The Kadomtsev-Petviashvili equation}

The Kadomtsev-Petviashvili equations are universal asymptotic models 
for dispersive systems in the weakly nonlinear, long wave regime, 
when the wavelengths in the transverse direction are much larger than in the direction of propagation.

The  (classical) Kadomtsev-Petviashvili (KP) equations read

\begin{equation}\label{KP}
(u_t+u_{xxx} +u u_x )_x \pm u_{yy} =0.
\end{equation}

Actually the (formal) analysis) in \cite{KaPe} consists in looking for a {\it  weakly transverse} perturbation 
of the one-dimensional transport equation

\begin{equation}\label{transp}
u_t+u_x=0.
\end{equation}

This perturbation amounts to adding a nonlocal term, leading to 

\begin{equation}\label {perttransp}
u_t+u_x+\frac{1}{2}\partial_x^{-1}u_{yy}=0.
\end{equation}

Here the operator $\partial_x^{-1}$ is defined via Fourier transform,

$$\widehat{\partial_x^{-1}f}(\xi)=\frac{i}{\xi_1}\widehat{f}(\xi),\,\text{where}\;\xi=(\xi_1,\xi_2).$$

When this same formal procedure is applied to the KdV equation written in the context of water waves (where $T\geq 0$ is the Bond number measuring the surface tension effects)

\begin{equation}\label {KdVWW}
u_t+u_x+uu_x+(1-\frac{1}{3}T)u_{xxx}=0, \;x\in \R,\;t\geq0,
\end{equation}

 this yields the KP equation in the form

\begin{equation}\label{KPbrut}
u_t+u_x+uu_x+(1-\frac{1}{3}T)u_{xxx}+\frac{1}{2}\partial_x^{-1}u_{yy}=0.
\end{equation}

By change of frame and scaling, \eqref{KPbrut} reduces to \eqref{KP} 
with the $+$ sign (KP II) when $T<\frac{1}{3}$ and the $-$ sign (KP I) when $T>\frac{1}{3}$.

\vspace{0.3cm}
Of course the same formal procedure could be applied to {\it any} one-dimensional weakly nonlinear dispersive  equation of the form

\begin{equation}\label {gen}
u_t+u_x+f(u)_x-Lu_x=0, \;x\in \R,\;t\geq0,
\end{equation}

where $f(u)$ is a smooth real-valued function (most of the time polynomial) and $L$ a linear operator taking into account the dispersion and defined in Fourier variable by

\begin{equation}\label {L}
\mathfrak F( Lu)(\xi)=p(\xi)\mathfrak F u(\xi),
\end{equation}

where the symbol $p(\xi)$ is real-valued. The KdV equation corresponds for instance to $f(u)=\frac{1}{2}u^2$ and $p(\xi)= -\xi^2.$  

This leads to a class of generalized KP equations

\begin{equation}\label {gKP}
u_t+u_x+f(u)_x-Lu_x+\frac{1}{2}\partial_x^{-1}u_{yy}=0, \;x\in \R,\;t\geq0.
\end{equation}

Thus one could have KP versions of the Benjamin-Ono, Intermediate Long Wave, Kawahara, etc... equations, but only the KP I and KP II equations are completely integrable (in some sense).

Let us notice, at this point, that alternative models  to KdV-type equations
\eqref{gen} are the equations of Benjamin--Bona--Mahony (BBM) type
\begin{equation}\label{gBBM}
u_{t}+u_{x}+uu_{x}+Lu_{t}=0
\end{equation}
with corresponding two-dimensional ``KP--BBM-type models'' (in the case $p(\xi)\geq 0$)
\begin{equation}\label{1.5}
u_{t}+u_{x}+uu_{x}+Lu_{t}+\partial_{x}^{-1}\partial_{y}^{2}u=0
\end{equation}
or, in the {\it derivative form}
\begin{equation}\label{1.10}
(u_{t}+u_x+uu_{x}+Lu_{t})_{x}+ \partial_{y}^{2}u=0
\end{equation}
and free group
\[
S(t)=e^{-t(I+L)^{-1}[\partial_{x}+ \partial_{x}^{-1}\partial_{y}^{2}]}\,\, .
\]

It was only after the seminal paper \cite{KaPe} that Kadomtsev-Petviashvili type equations have been  derived as asymptotic weakly nonlinear models (under an appropriate scaling) in various physical situations (see \cite{AS} for a formal derivation in the context of water waves  \cite{La, La2}, \cite {LS} for a rigorous one in the same context) and \cite{Ka} in a different context.
\begin{remark}
In some physical contexts (not in water waves!) one could consider higher dimensional transverse perturbations, which amounts to replacing $\partial_x^{-1}u_{yy}$ in \eqref{genKP} by $\partial_x^{-1}\Delta^{\perp}u$, where $\Delta^{\perp}$ is the Laplace operator in the transverse variables.

For instance, as in the one-dimensional case the KP I equation (in both two and three dimensions) also describes after a suitable scaling the long wave {\it transonic} limit of the Gross-Pitaevskii equation (see \cite {BGS} for the solitary waves and \cite{CR} for the Cauchy problem).
\end{remark}

Note again that in the classical KP equations, the distinction  between KP I and KP II arises from the {\it sign} of the dispersive term in $x$.

\subsection{KP by Inverse Scattering}

It is usual in the Inverse Scattering community to write the Kadomstsev-Petviashvili equations as 

\begin{equation}\label{KPIST}
\partial_x(\partial_t u+6u\partial_x u+\partial^3_{x} u)=-3\sigma^2\partial_t^2u,
\end{equation}

where $\sigma^2=1$ corresponding to KP II and $\sigma^2=-1$ to KP I.

\subsubsection{The KP II equation}

The direct scattering problem is associated to the heat equation with the initial potential $u_0(x,y):$

\begin{equation}\label{KPIIdirect}
-\partial_y\phi+\partial_x^2\phi+2ik\partial_x\phi+u_0\phi=0,\quad \phi_{|k|\to \infty}=1,
\end{equation}

and the scattering data are calculated by 

$$F(k)=(2\pi)^{-1}\text{sign}(\mathcal Re\; 
k)\int_{\R^2}u_0(x,y)\phi(x,y,k)\text{exp}\lbrace -i(k+\bar 
{k})x-(k^2-\bar{k}^2)y\rbrace dxdy.$$

The time evolution of the scattering data is trivial :

$$\mathcal F(k,t)=F(k)\text{exp} (4it(k^3+\bar{k}^3)).$$

The inverse scattering problem, that is the reconstruction of the potential $u(x,y,t)$ reduces to a $\bar{\partial}$ problem :

 \begin{equation}
    \label{InvKPII}
    \left\lbrace
    \begin{array}{l}
    \partial_k\phi=\psi F(-k)\text{exp}(itS), \\
    \partial_k\psi=-\phi F(k)\text{exp}(-itS)

    \end{array}\right.
    \end{equation}
    
    \begin{math}
    \begin{pmatrix} \phi\\
    \psi
    \end{pmatrix}
    \end{math}
    $\to$ \begin{math}
    \begin{pmatrix} 1\\
    1
    \end{pmatrix}
    \end{math}
as $|k| \to \infty,$

and where 

$$S=4(k^3+\bar{k}^3)+(k+\bar{k})\xi-i(k^2-\bar{k}^2)\eta, \;\xi=x/t,\;\text{and}\; \eta=y/t.$$

It turns out that the direct scattering problem can be solved only 
for small data in spaces of type $L^1\cap L^2,$ yielding global existence of uniformly bounded (in the space of $L^2$ functions with bounded Fourier transform) solutions of KP II provided $u_0$ has small derivatives up to order $8$ in $L^1\cap L^2(\R^2)$  (\cite {W}). We will see that PDE methods provide a much better result.

\begin{remark}
In \cite{Gr}, Grinevich has proven that the direct spectral problem is nonsingular for real nonsingular exponentially decaying at infinity, arbitrary large potentials. Unfortunately, this does not mean that the solution of the direct scattering problem belongs to an appropriate functional class for existence of an inverse scattering problem when $t>0.$ In fact, the direct problem and inverse problem are different
and the solvability of the first one does not give the automatic solvability of the second one. 
\end{remark}

\begin{remark}
Since KP II type equations do not have localized solitary waves (\cite{deBS}), one expects the large time behavior of solutions to  be just governed by scattering.  In particular, one can conjecture safely than the global solutions of KP II (that exist by the result of Bourgain, see \cite{Bourgain3} and the discussion below) should decay in the sup-norm as $1/t. $ This is also suggested by our numerical simulations \cite{KS2}.

\end{remark}

A very precise asymptotics as $t\to \infty$ is given in \cite{Ki2} 
(see also \cite{Ki}) for a specific class of scattering data. It 
differs according to different domains in the $(x,y,t)$ space, 
expressed in terms of the variables $\xi=x/t$ and $\eta =y/t$. The 
main term of the asymptotics has order $O(1/t)$ (which is exactly the 
decay rate of the free linear evolution, see \cite{Sa}) and rapidly 
oscillates. In one of the domains, the decay is $o(1/t).$ It is not 
clear however how the hypothesis on the scattering data translate to the space of initial data.

On the other hand, the Inverse Scattering method has been used 
formally in \cite{AV} and rigorously in \cite {AV2} to study the 
Cauchy problem for the KP II equation with nondecaying data along a line, that is $u(0,x,y)= u_{\infty}(x-vy)+ \phi(x,y)$ with $\phi(x,y)\to 0$ as $x^2+y^2\to \infty$ and $u_{\infty}(x)\to 0 $ as $|x|\to \infty$.  Typically, $u_{\infty}$ is the profile of a traveling wave solution $U({\bf k}.{\bf x}-\omega t)$ with its peak localized on the moving line ${\bf k}.{\bf x}=\omega t.$ It is a particular case of the $N$- soliton of the KP II equation discovered by Satsuma \cite{Sat} (see the derivation and the explicit form when $N=1,2$ in the Appendix of \cite{NMPZ}). As in all results obtained for KP equations by using the Inverse Scattering method, the initial perturbation of the non-decaying solution is supposed to be small enough in a weighted $L^1$ space (see \cite{AV2} Theorem 13).

\subsubsection{The KP I equation}
The direct scattering problem for KP I is associated to the Schr\"{o}dinger operator with potential

$$i\psi_t+\psi_{xx}=-u\psi.$$

As for the KP II case, there is a restriction on the size of the initial data to solve the direct scattering problem (see \cite{Ma}). In \cite{Zh} the
nonlocal Riemann-Hilbert problem for inverse scattering is shown to 
have a solution leading to the global solvability of the Cauchy problem (with a smallness condition on the initial data).
A formal asymptotic of small solutions is given in \cite{MaSaTa}. It would be interesting to provide a rigorous proof of this result.

It is proven in \cite{Su} that the solution constructed by  the IST belongs to the Sobolev spaces $H^s(\R^2), s\geq 0,$  provided the initial data is a small function in the Schwartz space $\mathcal S (\R^2),$ thus not assuming the zero mass constraint  (see subsection 4.2.1 below) contrary to the  result in \cite{Zh} (see also \cite{FoSu1} where the IST solution is shown to be $ C^\infty$ for a small Schwartz initial data).

Finally the Cauchy problem of the background of a one-line soliton is solved  formally (for small initial perturbations) in \cite{FP}.

\subsubsection{Conservation laws}

The KP equations being integrable have an infinite set of (formally) 
conserved quantities. For instance,  (see\cite{ZS3})  the KP I equation   has a Lax pair
representation. This in turn provides an algebraic procedure generating an
infinite sequence of conservation laws. More precisely, if $u$ is a formal
solution of the KP I equation then
$$
\frac{d}{dt}\Big[\int \chi_n\Big]=0,
$$
where $\chi_1=u$, $\chi_2=u+i\partial_{x}^{-1}\partial_y u$ and for $n\geq 3$,
$$
\chi_n=\Big(\sum_{k=1}^{n-2}\chi_k\, \chi_{n+1-k}\Big)+\partial_x\chi_{n-1}+
i\partial_{x}^{-1}\partial_{y}\chi_{n-1}\,\, .
$$
For $n=3$, we find the conservation of the $L^2$ norm, $n=5$ corresponds to
the energy functional giving the Hamiltonian structure of the KP I
equation, that is the following quantities are well defined and conserved by the flow 
(in an appropriate functional setting, see \cite{MST2})

\begin{eqnarray*}
M(\phi)
& = &
\int_{\R^2}
|u|^{2},
\\[0.5cm]
E(u) & = &
\frac{1}{2} \int_{\R^2} u_x^2
+
\frac{1}{2} \int_{\R^2} (\partial_x^{-1} u_y)^2 
-
\frac{1}{6}
\int_{\R^2} u^3 ,
\\[0.5cm]
F(u) & = & \frac{3}{2} \int_{\R^2} u_{xx}^2 + 5 \int_{\R^2} u_y^2
+ \frac{5}{6}\int_{\R^2} (\partial^{-2}_x u_{yy} )^2
 -\frac{5}{6} \int_{\R^2} u^2 \partial^{-2}_x u_{yy}
\\[0.4cm]
& & -\frac{5}{6} \int_{\R^2} u (\partial_x^{-1} u_y)^2
 +\frac{5}{4} \int_{\R^2} u^2\, u_{xx} + \frac{5}{24} \int_{\R^2} u^4 .
\end{eqnarray*}

As was noticed in \cite{MST2}, there is a serious analytical
obstruction to give sense to $\chi_9$ as far as $\R^2$ is considered as a
spatial domain. More precisely the conservation
law which controls $\|u_{3x}(t,\cdot)\|_{L^2}$
involves the $L^2$ norm of the term
$\partial_{x}^{-1}\partial_{y}(\phi^2)$which has no sense for a nonzero function $\phi$ in $H^3(\R^2)$ say. Actually one easily checks that if $\partial_{x}^{-1}\partial_{y}(\phi^2)\in L^2(\R^2),$ then   $\int_{\R^2} \partial _y(u^2) dx = \partial_y \int_{\R^2}u^2 dx \equiv 0,\; \forall y\in \R,$  which, with $u\in L^2(\R^2),$ implies that $u\equiv 0.$ Similar obstructions occur for the higher order ``invariants".

The fact that the invariants $\xi_n, \;n\geq 9$ do not make sense for a  nonzero function yields serious difficulties to define the so-called KP hierarchy.

\subsection{KP by PDE methods}
The basic difference between KP I and KP II as far as PDE techniques are concerned, is that KP I is {\it quasilinear} while  KP II is {\it semilinear}. We recall that this  means that the Cauchy problem for KP I cannot be solved by a Picard iterative scheme implemented on the integral Duhamel formulation, for any initial data in very  general spaces (that is the Sobolev spaces $H^s(\R^2), \forall s\in \R,$ or the anisotropic ones $H^{s_1,s_2}(\R^2), \forall s_1, s_2 \in \R)$ .  Alternatively, this implies that the flow map $u_0\mapsto u(t)$ cannot be 
{\it smooth} in the same spaces. Here are precise statements of those results from \cite{MST}.

\begin{theorem}\label{th1}
Let $s\in\R$  and $T$ be a positive real number.
Then there does not exist a space $X_T$ continuously
embedded in $C([-T,T],H^{s}(\R))$
such that there exists $C>0$ with
\begin{eqnarray}\label{contr1}
\|S(t)\phi\|_{X_T}\leq C\|\phi\|_{H^{s}(\R)},\quad \phi\in
H^{s}(\R),
\end{eqnarray}
and
\begin{eqnarray}\label{contr2}
\left\|\int_{0}^{t}S(t-t')\left[u(t')u_{x}(t')\right]dt'\right\|_{X_T}
\leq
C\|u\|_{X_T}^{2},\quad u\in X_T.
\end{eqnarray}
\end{theorem}

Note that (\ref{contr1}) and (\ref{contr2}) would be needed to
implement a Picard iterative scheme on the integral (Duhamel) formulation of the equation in the space $
X_T $. As a consequence of Theorem \ref{th1} we can obtain the
following result.

\begin{theorem}\label{ill}
Let $(s_1,s_2)\in \R^2$ (resp. $ s\in \R $). 
Then there exists no $T>0$ such that KPI  admits a unique local
solution defined on the interval $[-T,T]$ and such that the flow-map 
$$
S_t : \phi\longmapsto u(t), \quad t\in [-T,T]
$$ 
for (\ref{KdV}) is $C^{2}$ 
differentiable at zero from $H^{s_1,s_2}(\R^2) $ to $H^{s_1,s_2}(\R^2)$,
 (resp. from $ H^s(\R^2) $ to $ H^s(\R^2) $). 
\end{theorem}

\begin{remark}
It has been proved in \cite{KT} that the flow map cannot be uniformly continuous in the energy space.
\end{remark}

\begin{proof}
We merely sketch it (see \cite{MST} for details). Let

$$\sigma(\tau,\xi,\eta)=\tau-\xi^3-\frac{\eta^2}{\xi},$$

$$\sigma_1(\tau_1,\xi_1,\eta_1)=\sigma(\tau_1,\xi_1,\eta_1),$$

$$\sigma_2(\tau_{1},\xi,\eta_1,\tau_1,\xi_1,\eta_1)=\sigma(\tau-\tau_1,\xi-\xi_1,\eta-\eta_1).$$

We then define

$$
\chi(\xi,\xi_{1},\eta,\eta_{1}):=
3\xi\xi_{1}(\xi-\xi_1)-
\frac{(\eta\xi_1-\eta_{1}\xi)^{2}}{\xi\xi_{1}(\xi-\xi_1)}.
$$
Note that $\chi(\xi,\xi_{1},\eta,\eta_{1})=\sigma_1+\sigma_2-\sigma$.
The ``resonant'' function $\chi(\xi,\xi_{1},\eta,\eta_{1})$
plays an important role in our analysis. The ``large'' set of
zeros of $\chi(\xi,\xi_{1},\eta,\eta_{1})$ is responsible for the ill-posedness issues. In contrast, the corresponding resonant function for the KP II equation is

$$
\chi(\xi,\xi_{1},\eta,\eta_{1}):=
3\xi\xi_{1}(\xi-\xi_1)+
\frac{(\eta\xi_1-\eta_{1}\xi)^{2}}{\xi\xi_{1}(\xi-\xi_1)}.
$$

Since it is essentially the sum of two squares, its  zero set is small and this is the key point to establish the crucial bilinear estimate in Bourgain's method (\cite{Bourgain3}).
\end{proof}

\begin{remark}
 It is worth noticing that  the  property of the resonant set of the KP II equation was used by Zakharov \cite{Za2} to construct a Birkhoff formal form for the {\it periodic} KP II equation with small initial data. On the other hand, the fact that for the KP I equation the corresponding resonant set is non trivial is crucial in the construction of the counter-examples of  \cite{MST} and is apparently an obstruction to the Zakharov construction for the periodic KP I equation.
\end{remark}

Since the next property  of KP equations is only based on the presence of the operator $\partial_x^{-1}\partial_y^2$ we will consider it in the context of {\it generalized} KP type equations :

\begin{equation}\label {gentypeKP}
u_t+u_x+f(u)_x-Lu_x+\frac{1}{2}\partial_x^{-1}u_{yy}=0, \;(x,y)\in \R^2,\;t\geq0,
\end{equation}

where $\widehat{Lu}(\xi)=p(\xi)\hat{u}(\xi),$ $p$ real and and $f(u)$ is a nonlinear function, for instance $f(u)=\frac {1}{q+1}u^{q+1}.$

An important particular case is the {\it generalized KP equation}

\begin{equation}\label {genKP}
u_t+u_x+u^ru_x+u_{xxx}\pm\partial_x^{-1}u_{yy}=0, \;(x,y)\in \R^2,\;t\geq0,
\end{equation}

where $r\in \N$ or $r=\frac{p}{q},$ $p,q$ relatively prime integers, $q$ odd.

\subsubsection{The zero mass constraint}

In \eqref{gentypeKP} or \eqref{genKP}, it is implicitly assumed that the operator $\partial_{x}^{-1}\partial_{y}^{2}$
is well defined, which a~priori imposes a constraint on the solution $u$, which, in some sense, has to be an $x$-derivative.
This is achieved, for instance, if  $u\in {\mathcal S}'(\R^2)$ is such that
\begin{equation}\label{1.6}
\xi_{1}^{-1}\,\xi_{2}^{2}\,\widehat{u}(t,\xi_1,\xi_2)\in {\mathcal S}'(\R^2)\, ,
\end{equation}
thus in particular if $\xi_{1}^{-1}\,\widehat{u}(t,\xi_1,\xi_2)\in {\mathcal S}'(\R^2)$.
Another possibility to fulfill the constraint is to write $u$ as
\begin{equation}\label{1.7}
u(t,x,y)=\frac{\partial}{\partial x}v(t,x,y),
\end{equation}
where $v$ is a continuous function having a classical derivative with respect
to $x$, which, for any fixed $y$ and $t\neq 0$, vanishes when $x\rightarrow \pm \infty$.
Thus one has
\begin{equation}\label{1.8}
\int_{-\infty}^{\infty}u(t,x,y)dx=0,\qquad y\in\R,\,\,\, t\neq 0,
\end{equation}
in the sense of generalized Riemann integrals.
Of course the differentiated version of \eqref{gentypeKP}, \eqref{genKP}, for instance
\begin{equation}\label{1.9}
(u_{t}+u_{x}+uu_{x}-Lu_{x})_{x}+ \partial_{y}^{2}u=0,
\end{equation}
can make sense without any constraint of type (\ref{1.6}) or  (\ref{1.8}) on $u$, 
and so does the Duhamel integral representation of \eqref{genKP}, \eqref{gentypeKP}, for instance
\begin{equation}\label{1.11}
u(t)=S(t)u_{0}-\int_{0}^{t}S(t-s)(u(s)u_{x}(s))ds,
\end{equation}
where $S(t)$ denotes the (unitary in all Sobolev spaces $H^{s}(\R^2)$) group associated with \eqref{gentypeKP},
\begin{equation}\label{1.12}
S(t)=e^{-t(\partial_{x}-L\partial_{x}+\partial_{x}^{-1}\partial_{y}^{2})}\,\, .
\end{equation}

In particular, the results established on the Cauchy problem for KP type equations  which use the Duhamel (integral ) formulation (see for instance \cite{Bourgain} for the KP II equation and  \cite{ST4} for the KP II BBM equation)  are valid {\it without} any constraint on the initial data.

 One has however to be careful {\it in which sense} the {\it differentiated equation} is taken. For instance let us consider the integral equation
 
 \begin{equation}\label {dudu}
 u(x,y,t)=S(t) u_0(x,y)-\int_0^t S(t-t')\lbrack u(x,y,t')u_x(x,y,t)\rbrack dt',
 \end{equation}
 
where $S(t)$ is here the KP II group, for initial data $u_0$ in $H^s(\R^2),$ $s>2$, (thus $u_0$ does not satisfy any zero mass constraint), yielding a local solution $u\in C(\lbrack 0,T\rbrack; H^s(\R^2))$.

By differentiating \eqref{dudu} first with respect to $x$ and then with respect to $t,$ one obtains the equation

$$\partial _t\partial_xu+\partial_x(uu_x)+\partial_x^4u+\partial_y^2u =0 \quad\text{in} \quad C(\lbrack 0,T\rbrack; H^{s-4}(\R^2)).$$

However, the identity $ \partial_t\partial_xu =\partial_x\partial_t u$ holds only in a very weak sense, for example in $\mathfrak D'((0,T)\times \R^2).$

On the other hand, a constraint has to be imposed when using the Hamiltonian
formulation of the equation. In fact, the Hamiltonian for \eqref{1.9} is
\begin{equation}\label{1.13}
\frac{1}{2}\int\left[-u\,Lu+(\partial_{x}^{-1}u_y)^2+u^2+\frac{u^3}{3}\right]
\end{equation} 
and the  Hamiltonian associated with (\ref{1.10}) is
\begin{equation}\label{1.14}
\frac{1}{2}\int\left[(\partial_{x}^{-1}u_y)^2+u^2+\frac{u^3}{3}\right]\, .
\end{equation}

It has been established in \cite{MST1} that, for a rather  general class of KP or KP--BBM
equations,  the solution of the Cauchy problem obtained for 
\eqref{1.9}, \eqref{1.10} (in an appropriate functional setting) satisfies the zero-mass constraint in $x$  for any $t \neq 0$ (in a sense to be made precise below), even if the initial data does not. This is a manifestation of the infinite speed of propagation inherent to KP equations. Moreover, KP type equations display a striking {\it smoothing effect}  : if the initial data belongs to the space $L^1(\R^2)\cap H^{2,0}(\R^2)$ and if $u\in C(\lbrack 0,T\rbrack; H^{2,0}(\R^2))$  \footnote {We will see below  that KP type equations (in particular the classical KP I and KP II equations) do possess solutions in this class.} is a solution in the sense of distributions, then, for any $t>0,$  $u(.,t)$ becomes a {\it continuous} function of $x$ and $y$ (with zero mean in $x$).  Note that the space  $L^1(\R^2)\cap H^{2,0}(\R^2)$ is not included in the space of continuous functions.

\begin{figure}
[!htbp]
\begin{center}
\includegraphics[width=0.7\textwidth]{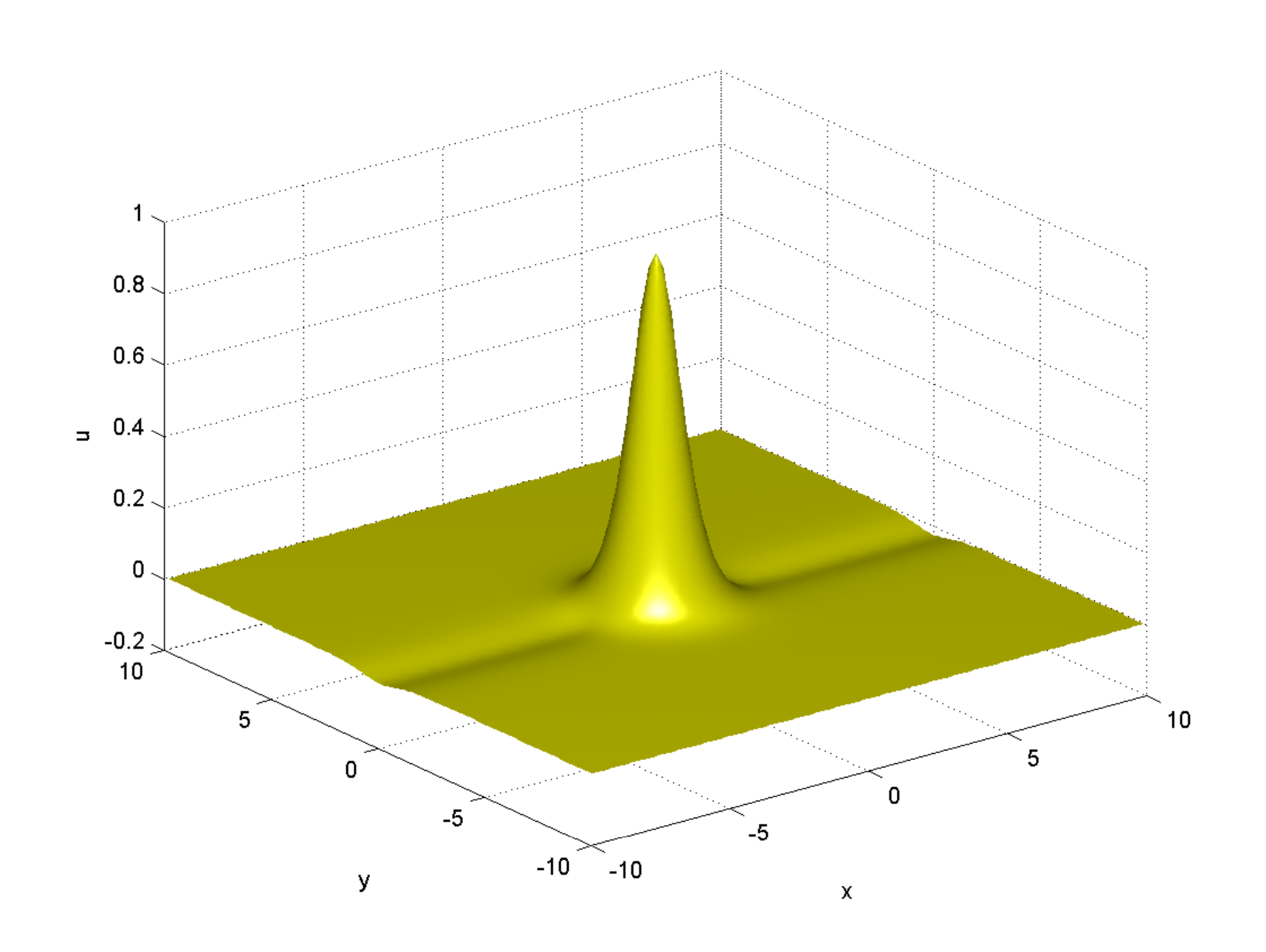} 
\caption{Solution at time
$t =4.8\times 10^{-4}$ to the KP-I
solution with initial data $\text{sech }^{2}\sqrt{x^{2}+y^{2}}$, 
which are not subject to the zero mass constraint.}
\label{kpnocon}
\end{center}
\end{figure}

The key point when proving those results is a careful analysis of  the fundamental solution 
of KP-type equations \footnote{ In the case of KP II, one can use the 
explicit form of the fundamental solution found in \cite{Re}.} which turns out to be  an $x$ derivative of a continuous function of $x$ and $y$, $C^1$ with respect to  $x,$ which, for fixed $t\neq 0$ and $y$, tends to zero as
$x\rightarrow \pm \infty$. Thus its generalized Riemann integral in $x$
vanishes for all values of the transverse variable $y$ and of $t\neq 0$. A
similar property can be established for the solution of the nonlinear
problem (see \cite{MST1}). Those results have been checked in the 
numerical simulations of \cite{KSM2} as can be seen in 
Fig.~\ref{kpnocon} taken from this reference. It can be seen that for 
initial data not satisfying the constraint, 
after an arbitrary short time some sort of infinite trench forms the 
integrall over which just ensures that the constraint holds at all 
$t$. 

 We have already referred  to \cite{FoSu1}, \cite{Su}  for a   
 rigorous approach to the Cauchy problem with (small) initial data which do not satisfy the zero-mass condition via the Inverse Spectral Method  in the integrable case.

Nevertheless, the singularity at $\xi_1=0$ of the dispersion relation of KP type  equations make them rather {\it bad} asymptotic models. First the singularity at $\xi_1=0$ yields a very bad approximation of the dispersion relation of the original system (for instance the water wave system) by that of the KP equation. 

Another drawback is the poor error estimate between the KP solution and the solution of the original problem. This has been established clearly in the context of water waves (see \cite{La, La2,LS}).

\subsubsection{The Cauchy problem by PDE techniques}

All the KP type equations can be viewed as a linear skew-adjoint perturbation of the Burgers equation. Using this structure, it is not difficult (for instance by a compactness method)  to prove that the Cauchy problem is locally well-posed for data in the Sobolev spaces $H^s(\R^2),$ $s>2$  (see \cite{U}, \cite{Sa}, \cite{IN} for results in this direction).

Unfortunately, this kind of result does not suffice to obtain the {\it global} well-posedness of the Cauchy problem. This would need to use the {\it conservation laws} of the equations. For {\it general} KP type equations, there are only two of them, the conservation of the $L^2$ norm and the conservation of the energy (Hamiltonian). For the general equation \eqref{genKP} where $f(u) =\frac{1}{p+1}u^{p+1},$ and without the transport term $u_x$ (which can be eliminated by a change of variable), we recall that the Hamiltonian reads

\begin{equation}\label{Ham}
E(u)=\frac{1}{2}\int\left[-u\,Lu+(\partial_{x}^{-1}u_y)^2+\frac{u^{p+2}}{p+2}\right],
\end{equation} 

and for the classical KP I/II equations

\begin{equation}\label{HamKP}
E(u)=\frac{1}{2}\int\left[u_x^2\pm(\partial_{x}^{-1}u_y)^2-\frac{u^{3}}{3}\right],
\end{equation} 

where the $+$ sign corresponds to KP I and the $-$ sign to KP II.

We recall  that the ``integrable" KP I and KP II equations possess 
more conservation laws, but only a finite number of them make sense rigorously (see above).

In any case, it is clear that for KP II type equation, the Hamiltonian is useless to control any Sobolev norm, and to obtain the {\it global}  well-posedness of the Cauchy problem one should consider $L^2$ solutions, a very difficult task. On the other hand, for KP I type equations, one may hope (for a not too strong nonlinearity) to have a global control in the {\it energy} space, that is

\begin{equation}\label{energysp}
E= \lbrace u\in L^2,\; uLu\in L^2,\;\partial_x^{-1}u_y\in L^2 \rbrace.
\end{equation}

For the usual KP I equation, $E$ reduces to 

$$Y=\lbrace u\in L^2,\; u_x\in L^2,\;\partial_x^{-1}u_y\in L^2 \rbrace.$$

The problem is thus reduced to proving the local well-posedness of the Cauchy problem in spaces of very low regularity,  a difficult matter  which has attracted a lot of efforts in the recent years.

\begin{remark}
By a standard compactness method, one obtains easily the existence of global weak finite energy solutions (without uniqueness) to the KP I equation (see {\it eg} \cite{Tom}).
\end{remark}

A fundamental step for KP II is made in \cite{Bourgain} who proved that the Cauchy problem for the KP II equation is locally (thus globally in virtue of the conservation of the $L^2$ norm) for data in $L^2(\R^2),$ and even in $L^2(\T^2).$ This result is  based on an  iterative method implemented on the Duhamel formulation, in the functional framework of  the {\it Fourier restriction} $X^{s,b}$ spaces of Bourgain (see a nice description of this framework in \cite{Gi}). The basic {\it bilinear estimate} which aims to regain the loss of one $x$-derivative  uses in a crucial way the fact  (both in the periodic and full space case) that the dispersion relation of the KP II equation induces the triviality of a  {\it resonant set} (the zero set of the aforementioned resonant function). With in particular the injection of various linear dispersive estimates (see for instance \cite{Sa,BAS}), Bourgain's result was later improved (see  \cite{TaTz}, \cite{IM}, \cite{Ha}, \cite{HHK} and the references therein) to allow the case of initial data in negative order Sobolev spaces. 

We also would like to mention the paper by Kenig and Martel \cite{KM} where the Miura transform is used to prove the global well-posedness of a modified KP II equation.

Moreover it was proven in \cite{MST5} that the Cauchy problem for KP II is globally well-posed with initial data
\begin{equation}\label{1.2}
u(0,x,y)=\phi(x,y)+\psi_c(x,y),
\end{equation}
where $\psi_c$ is the profile\footnote{This means that
$\psi(x-ct,y)$ solves KP II.} of a non-localized (i.e. not
decaying in all spatial directions) traveling wave of the KP II
equation moving with speed $c\neq 0$.

We recall (\cite {deBS}) that, contrary to the KP I equation, the KP II equation does not possess any {\it localized in both directions} traveling wave solution.

The background solution $\psi_c$ could be for instance the line soliton (1-soliton) of the Korteweg-
de Vries (KdV) equation
\begin{equation}\label{KdV}
s_c(x,y)=\frac{3c}{2}\,{\rm cosh}^{-2}\Big(\frac{\sqrt{c}\, x}{2}\Big),
\end{equation}

or the N-soliton solution of the KdV equation, $N\geq 2.$

The KdV N-soliton is of course considered as a two dimensional (constant
in $y$) object. 

There are two suitable settings for an initial perturbation $\phi$, either localized in $x$ and periodic in $y$ ({\it eg} $\phi \in L^2(\R\times \T)$) or localized ({\it eg} $\phi \in L^2(\R^2).)$

 Solving the Cauchy  problem in both those functional settings can be 
 viewed as a preliminary step towards the study of the dynamics of 
 the KP II equation on the background of a non fully localized 
 solution, in particular towards  a PDE proof of the nonlinear   
 stability of the KdV soliton or N-soliton  with respect to 
 transversal perturbations governed by the KP II flow.  This has been 
 established in \cite{AV2} Proposition 17 by Inverse Scattering 
 methods. The advantage of the PDE approach  is that it can be straightforwardly applied to non integrable equations such as  the higher order KP II equations (see \cite{ST2}, \cite{ST3}).


%

We now state the main result of \cite{MST5} in the two aforementioned functional settings.

\begin{theorem}\label{theo1}
The Cauchy problem associated with the KP II equation is globally well-posed in
 $ H^s(\R\times\T) $ for any $ s\ge 0$.
 \end{theorem}
\begin{theorem}\label{theo2}
Let $\psi_c(x-ct,y)$ be a solution of the KP II equation such that for some
 $s \geq 0 $,
$$
J^s \psi_c\,:\, \R^2\longrightarrow \R
$$
is bounded and belongs to $ L^2_x L^\infty_y $   \footnote{The bounds can of course depend on the propagation speed $c$.}.
Then for every $\phi\in H^s(\R^2) $ there exists a unique  global solution
$u$ of KP II with initial data \eqref{1.2} satisfying for all $T>0$,
$$
[u(t,x,y)-\psi_c(x-ct,y)]\in X^{1/2+,s}_T \cap X^{3/4+,s}_T\cap C([0,T];H^s(\R^2)) .
$$
Furthermore, for all $ T> 0$,  the map $ \phi \mapsto u $ is continuous
from $ H^s(\R^2)$ to $ C([0,T];H^s(\R^2)))$.
\end{theorem}

\begin{remark}\label{example}
As was previously noticed, the hypothesis on $\psi_c$ in Theorem 
\ref{theo2} is satisfied by the N-soliton solutions of the KdV 
equation, but not by a  function $\psi$ which is non-decaying along a 
line $\lbrace (x,y)\vert x-vy=x_0\rbrace,$ as for instance the 
line-soliton of the KP II equation which has the form $\psi (x-vy-c t).$

However, the change of variables $(X=x, Y=x-vy)$ transforms the KP II equation into 

$$u_t-2vu_Y+v^2u_X+u_{XXX}+uu_X+\partial_X^{-1}u_{YY}=0$$

and the analysis applies to this equation with aninitial datum  which is  a localized (in $(X,Y)$) perturbation of  $\psi(X).$
\end{remark}


As was previously mentioned the Cauchy problem for the KP I  cannot be solved by a Picard iteration on the integral formulation and one has to implement  instead sophisticated {\it compactness methods} to obtain the local well-posedness and  to use the conservation laws to get global solutions.

For the classical KP I equation, the first global well-posedness 
result for arbitrary large initial data in a suitable function space  was obtained in \cite{MST2} in the space 

$$
Z=\{
\phi\in L^{2}(\R^2)~:~
\|\phi\|_{Z}<\infty
\},
$$
where
$$
\|\phi\|_{Z}
=
|\phi|_{2}+|\phi_{xxx}|_{2} + |\phi_y |_2
+ |\phi_{xy}|_{2}
+
|\partial_{x}^{-1}\phi_{y}|_{2}
+
|\partial_{x}^{-2}\phi_{yy}|_{2}.
$$
By an anisotropic Sobolev embedding theorem (cf. \cite{bes})
if $\phi\in Z$ then $\phi, \phi_{x}\in L^{\infty}(\R^2),$  so the global solution is uniformly bounded in space and time.
Moreover, if $\phi\in Z$
then the first three formal invariants  $M(\phi),  E(\phi), F(\phi)$ are well defined and conserved.
Furthermore it is easily checked that any finite energy  solitary 
waves (in particular the lumps, see below) of the KP I equation belong to $ Z $.

The proof is based on  a rather sophisticated compactness method and 
uses the first invariants of the KP I equation to get global in time 
bounds. As already mentioned, only a small number of the formal invariants make sense and in order to control $|\phi_{xxx}|_{2}$
one is thus led to introduce a {\it quasi-invariant} (by skipping the non well defined terms) which eventually will provide the desired bound. There are also serious technical difficulties to {\it justify} rigorously the conservation of the meaningful invariants along the flow and to control the remainder terms

The result of \cite{MST2} was extended by Kenig \cite{K} (who considered initial data in a larger space), and by Ionescu, Kenig and Tataru \cite{IKT} who proved that the KP I equation is globally well-posed in the energy space $Y$.

Moreover it is  proven in \cite{MST3}  that the Cauchy problem for 
the KP I equation is globally well-posed
for initial data which are localized perturbations (of arbitrary size) of a
{\bf non-localized} i.e. not decaying in all directions) traveling wave
solution (e.g. the KdV line solitary wave or the Zaitsev solitary waves which
 are localized in $x$ and $y$ periodic or conversely (see Section 4.2 below).

\subsubsection{Long time behavior}
The results above do not give information on the behavior of the  
global solution  for large time. Actually no result in this direction 
is known by PDE techniques. However, one can make precise the large time behavior  of small solutions to the generalized KP equations \eqref{genKP} when $r\geq 2$.


Actually, it is shown in \cite{HNS, Ni} that for initial data small in an appropriate weighted Sobolev space, \eqref{genKP} for $r\geq 2$ has a unique global solution satisfying

$$|u(.,t)|_\infty\leq C (1+|t|)^{-1} (\text{Log} (1+|t|))^\kappa,$$

$$|\partial_x u(.,t)|_\infty \leq C(1+|t|)^{-1},$$

where $\kappa=1$ when $r=2$ and $\kappa=0$ when $r>2.$

This result does not distinguish between the KP I and KP II case since it relies on the (same) large time asymptotic of the KP I and KP II groups $S_{\pm}(t)=e^{it(-\partial_x^3\pm \partial{-1}_x\partial_y^2)},$ namely (see \cite{Sa})

$$|S_{\pm}(t)\phi|_{\infty}\leq \frac{C}{|t|} |\phi|_1.$$

\begin{remark}
Since this phenomena does not happen for the classical KP I and KP II 
equation, we do not comment on the possible  {\it blow-up} in finite 
time of solutions to the generalized KP I equation \eqref{genKP} when $r\geq\frac{4}{3}.$ We refer to \cite{Liu} for a theoretical study and to \cite{KP2,KS2} for numerical simulations. 
\end{remark}

\subsection{Solitary waves}

We are interested here in localized solitary wave solutions to the KP equations, that is solutions of KP equations of the form

$$u(x,y, t)= \psi_c(x-ct, y),$$

where $y$  is the transverse variable and $c>0$ is  the solitary wave velocity.

The solitary wave is said to be {\it localized} if $\psi_c$ tends to zero at infinity in all directions. For such solitary waves, the  {\it energy space} $Y$ is natural. Recall that

$$Y=Y(\R^2)=\lbrace u\in L^2(\R^2),\; u_x \in L^2(\R^2),\;\partial _x^{-1}u_y \in L^2(\R^2)\rbrace,$$

and throughout this section we will deal only with {\it finite energy solitary waves}.

Due to its integrability properties, the KP I equation possesses a localized, finite energy, explicit solitary wave, the {\it lump}:

\begin{equation}\label{Lump}
\phi_c(x-ct,y)=\frac{8c(1-\frac{c}{3}(x-ct)^2+\frac{c^2}{3}y^2)}{(1+\frac{c}{3}(x-ct)^2+\frac{c^2}{3}y^2)^2}.
\end{equation}

The formula was found in  \cite{MZBM} where one can also find a study  of the interaction of lumps. The interactions do not result in a phase shift as in the
case of line solitons (KdV solitons). More general rational solutions of the KPI equation were subsequently found (\cite{Kr, SatAb, PeSt,Pel, Sat}). These solutions were incorporated within the framework of the IST in \cite{VA} where it was observed that, in general, the spectral characterization of the potential must include, in addition to
the usual information about discrete and continuous spectrum, an integer-valued topological quantity (the {\it index} or winding number), defined by an appropriate two-dimensional integral involving both the solution of the KP equation and the corresponding scattering eigenfunction.

Another interesting explicit solitary wave of the KP I equation which 
is {\it localized in $x$ and periodic in $y$} has been found by 
Zaitsev \cite{Zai}. It reads

\begin{equation}\label{Za}
Z_c(x, y)=12\alpha^2\frac{1-\beta \cosh (\alpha x)\cos (\delta y)}{\lbrack \cosh (\alpha x)-\beta\cos (\delta y)\rbrack^2},
\end{equation}

where 

$$(\alpha, \beta)\in ( 0, \infty)\times (-1, +1),$$

and the propagation speed is given by

$$c=\alpha ^2\frac{4-\beta ^2}{1-\beta^2}.$$

Let us observe that the transform $\alpha\rightarrow i\alpha$,
$\delta\rightarrow i\delta$, $c\rightarrow ic$ produces solutions of
the KP I equation which are periodic in $x$ and localized in $y$.

No real non-singular rational solutions are known for KP II. Moreover, it was established in \cite{deBS} that no localized solitary waves exist for the KP II equation (and generalized KP II equations).

For obvious (stability) issues it is important to characterize the solitary waves by variational principles. We will consider in fact the  slightly more  general class of generalized KP I equations

\begin{equation}\label {genKPI}
u_t+u_x+u^ru_x+u_{xxx}-\partial_x^{-1}u_{yy}=0, \;(x,y)\in \R^2,\;t\geq0,
\end{equation}

where again $r\in \N$ or $r=\frac{p}{q},$ $p,q$ relatively prime integers, $q$ odd.

 Solitary waves are looked for in the energy space $Y(\R^2)$ which can also be defined (see \cite{deBS})  as the closure of the space of $x$ derivatives of smooth and compactly supported functions in $\R^2$ for the norm
 
$$\| \partial_x f \|_{Y(\R^2)} \equiv \Big( \| \nabla f \|_{L^2(\R^2)}^2 + \| \partial_x^2 f \|_{L^2(\R^2)}^2 \Big)^\frac{1}{2}.$$

The equation of a solitary wave $\psi$ of speed $c$ is given by

\begin{equation}
\label{SW}
c\partial_x \psi - \psi \partial_x \psi - \partial_x^3 \psi + \partial_x^{- 1} (\partial_y^2 \psi) = 0,
\end{equation}

which implies 
\begin{equation}
\label{SW2}
c\partial_{xx} \psi - (\psi \partial_x \psi)_x - \partial_x^4 \psi +  \partial_y^2 \psi = 0,
\end{equation}

When $\psi \in Y(\R^2)$, the function $\partial_x^{- 1} \partial_y ^2\psi$ is well-defined (see \cite{deBS}), so that \eqref{SW} makes sense.

Given any $c > 0$, a solitary wave $\psi_c$ of speed $c$ is deduced from a solitary wave  $\psi_1$ of velocity $1$ by the scaling
\begin{equation}
\label{scalingsw}
\psi_c(x, y) = c \psi_1(\sqrt{c} x, c y).
\end{equation}

\vspace{0.3cm}
We now introduce the important notion of {\it ground state} solitary waves.

We set

$$E_{KP}(\psi) = \frac{1}{2} \int_{\R^2} (\partial_x \psi)^2 + \frac{1}{2} \int_{\R^2} (\partial_x^{-1}\partial_y \psi)^2 - \frac{1}{2(p+2)} \int_{\R^2} \psi^{p+2},$$
and we define the action

$$S(N) = E_{KP}(N) + \frac{c}{2} \int_{\R^2} N^2.$$

We call {\it ground state}, a solitary wave $N$ which minimizes the action $S$ among all finite energy non-constant solitary waves of speed $c$ (see \cite{deBS} for more details). 

It was proven in \cite{deBS} that ground states exist if and only if $c>0$ and $1\leq p<4$. Moreover, when $1\leq p<\frac{4}{3},$ the ground states are minimizers of the Hamiltonian  $E_{KP}$ with prescribed mass ($L^2$ norm).

\begin{remark}
When $p=1$ (the classical KP I equation), we should emphasize that it 
is unknown (but conjectured) whether the lump solution is a ground 
state. This important issue is of course related to the {\it uniqueness} of the ground state or of localized solitary waves, up to symmetries. A similar question stands for the focusing nonlinear Schr\"{o}dinger equation but it can be solved there because the ground state is shown to be radial and and uniqueness follows from (non trivial!) ODE arguments (see \cite{Kw}). Of course such arguments cannot work in the KP I case since the lump (or the ground states) are not radial.
\end{remark}

It turns out that qualitative properties of solitary waves which are that of the lump solution can be established for a large class of KP type equations. Ground state solutions are shown in \cite{deBS2} to be {\it cylindrically symmetric}, that is radial with respect to the transverse variable up to a translation of the origin.

On the other hand, {\it any} finite energy solitary wave is 
infinitely smooth (at least when the exponent $p$ is an integer) and decays with an {\it algebraic} rate $r^{-2}$  \cite{deBS2}. Actually the decay rate is sharp in the sense that a solitary wave {\it cannot} decay faster that $r^{-2}.$

Moreover a precise {\it asymptotic expansion} of the solitary waves  has been obtained by Gravejat (\cite{Gra}).

\begin{remark}

The fact whether the ground states are or are not minimizers of the 
Hamiltonian is strongly linked to the {\it orbital stability} of the 
set $\mathcal G_c$ of ground states of velocity $c$. We recall  that the {\it uniqueness}, up to the obvious symmetries, of the ground state of velocity $c$ is a difficult open problem, even for the classical KP I equation.

Saying that $\mathcal G_c$ is orbitally stable in $Y$ means that for $\phi_c\in \mathcal G_c,$ then for all $\epsilon >0, $ there exists $\delta>0$ such that if $u_0\in Y$\footnote {Some extra regularity on $u_0$ is actually needed.} is such that $||u_0-\phi_c||_Y\leq \delta,$ then the solution $u(t)$ of the Cauchy problem initiating from $u_0$ satisfies

$$\sup_{t\geq 0} \inf_{\psi \in \mathcal G}||u(t)-\psi||_Y\leq \epsilon.$$

Of course, the previous inequality makes full sense only if one knows that the Cauchy problem is globally well-posed. As we have previously seen, this is the case for the classical KP I equation \cite{MST3}, even in the energy space \cite{IKT}, but is still an open problem for the generalized KP I equation  when $1<p<4/3.$

Actually it is proved in \cite{deBS3} that the ground state solitary 
waves of the generalized KP I equations \eqref{genKPI} are orbitally stable in dimension two if and only if $1\leq p<\frac{4}{3}$. 
\end{remark}

\subsection{Transverse stability of the line soliton}

The KP I and KP II equations behave quite differently with respect to the {\it transverse} stability of the KdV $1$-soliton.

Zakharov \cite{Za} has proven, by exhibiting an explicit perturbation
using the integrability, that the  KdV $1$-soliton is {\it nonlinearly} unstable for the KP I flow. Rousset and Tzvetkov \cite{RT1, RT2, RT3} have given an alternative proof of this result, which does not use the integrability, and which can  thus be implemented on other problems ({\it eg} for nonlinear Schr\"{o}dinger  equations).

The {\it nature} of this instability is not known (rigorously) and one has to rely on   numerical simulations (\cite{KS2, PeSt2}).

On the other hand, Mizomachi and Tzvetkov \cite{MiTz} have recently proved the $L^2(\R\times \T)$ orbital stability of the KdV $1$-soliton for the KP II flow. The perturbation is thus localized in $x$ and periodic in $y$.  The proof involves  in particular the Miura transform used in \cite{KM} to established the global well-posedness for a modified KP II equation. 

Such a result is not known (but expected) for a perturbation which is localized in $x$ and $y$.

The transverse stability of the KdV soliton has been numerically 
studied in \cite{KS2} from which the figures in this subsection are 
taken. 
We considered perturbations of the form
\begin{equation}
    u_{p}=6(x-x_{1})\exp(-(x-x_{1})^{2})\left(\exp(-(y+L_{y}\pi/2)^{2})
    +\exp(-(y-L_{y}\pi/2)^{2})\right),
    \label{pert}
\end{equation}
which are in the 
Schwartz class for both variables and satisfy the zero mass constraint. They are of the same order of magnitude as the KdV 
soliton, i.e., of order $0(1)$, and thus
test the nonlinear stability of the KdV soliton. As discussed in 
\cite{KS2}, the computations are carried out in a doubly periodic 
setting, i.e., $\mathbb{T}^{2}$ and not on $\mathbb{R}^{2}$. 

\begin{figure}
[!htbp]
\begin{center}
\includegraphics[width=\textwidth]{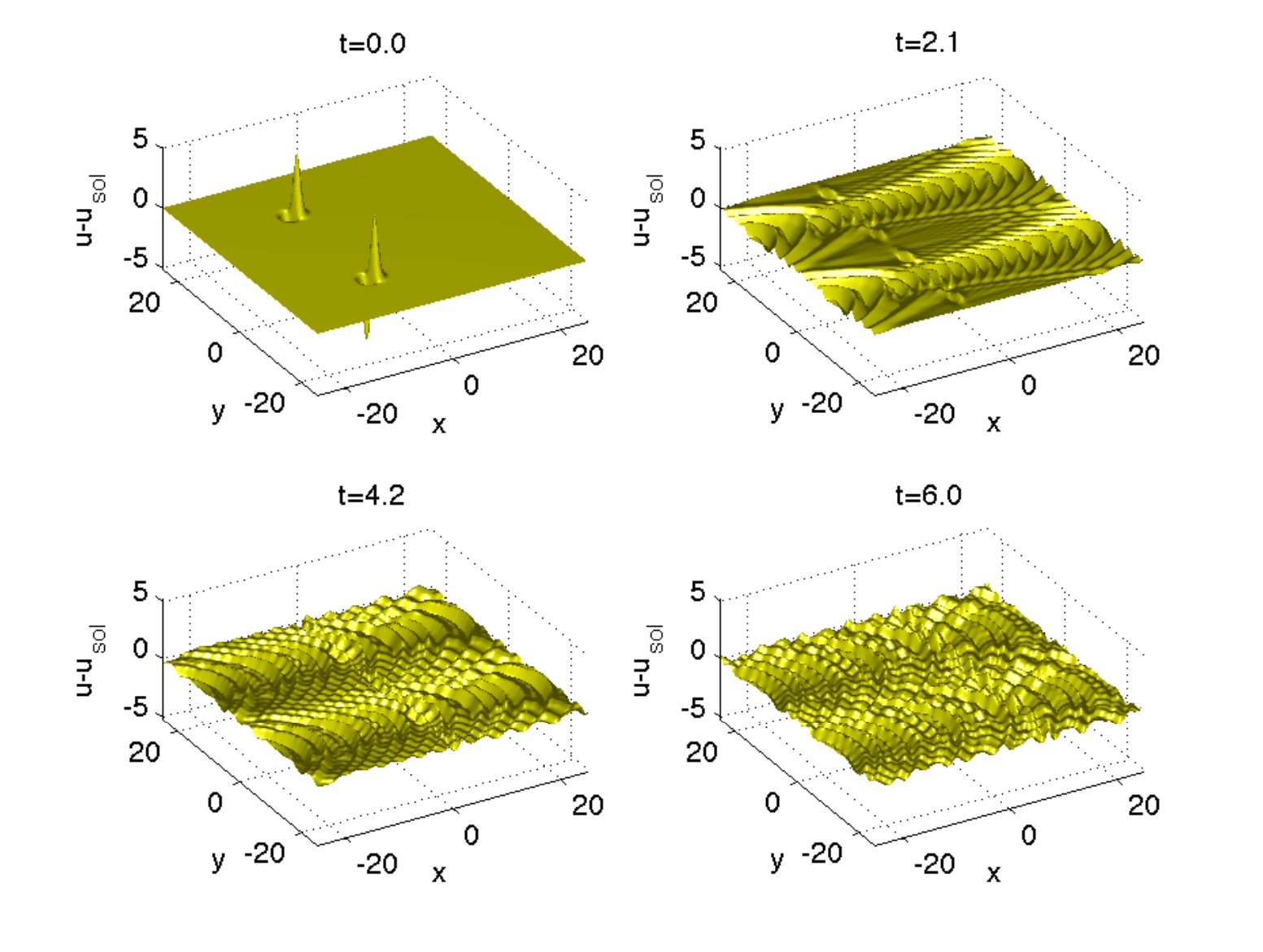} 
\caption{Difference of the solution to the  KP II equation  for initial 
data given by the KdV soliton plus perturbation, 
$u(0,x,y)=u_{sol}(x+2L_{x},0)$ and perturbation $u_{p}=6(x-x_{1})\exp(-(x-x_{1})^{2})\left(\exp(-(y+L_{y}\pi/2)^{2})
    +\exp(-(y-L_{y}\pi/2)^{2})\right)$, $x_{1}=-L_{x}$ and the 
KdV soliton for various values of $t$.}
\label{kpIIsolperd}
\end{center}
\end{figure}

For KP II we consider the initial data 
$u_{0}(x,y)=12\mbox{sech}^{2}(x+2L_{x})-u_{p}$ ($x_{1}=x_{0}/2$), 
i.e., a 
superposition of the KdV soliton and the not aligned perturbation 
which leads to the situation shown in Fig.~\ref{kpIIsolperd}. The 
perturbation is dispersed in the form of tails to infinity which 
reenter the computational domain because of the imposed periodicity. 
The soliton appears to be unaffected by the perturbation which 
eventually seems to be smeared out in the background of the soliton.

\begin{figure}
[!htbp]
\begin{center}
\includegraphics[width=\textwidth]{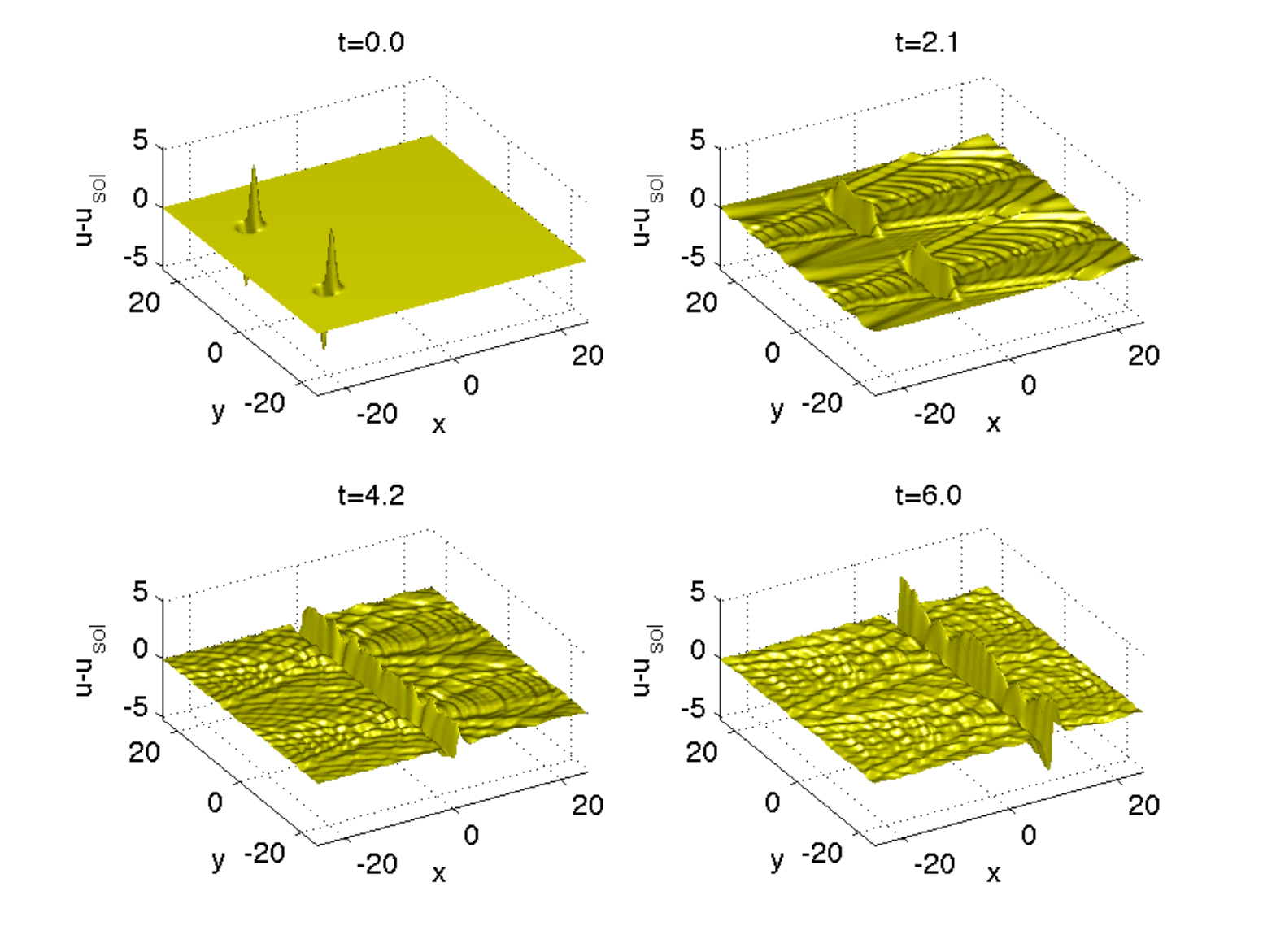} 
\caption{Difference of the solution to the  KP II equation  for initial 
data given by the KdV soliton $u_{sol}(x+2L_{x},0)$ 
plus perturbation $u_{p}=6(x-x_{1})\exp(-(x-x_{1})^{2})\left(\exp(-(y+L_{y}\pi/2)^{2})
    +\exp(-(y-L_{y}\pi/2)^{2})\right)$, $x_{1}=-2L_{x}$, and the 
KdV soliton for various values of $t$.}
\label{kpIIsolper}
\end{center}
\end{figure}

The situation is somewhat different if the perturbation and the initial 
soliton are centered around the same $x$-value initially, i.e., the 
same situation as above with $x_{1}=x_{0}=-2L_{x}$. 
In Fig.~\ref{kpIIsolper} we show the 
difference between the numerical solution and the KdV soliton 
$u_{sol}$ for several times for this case.
It can be seen that the initially localized perturbations spread in 
$y$-direction, i.e., orthogonally to the direction of propagation and 
take finally themselves the shape of a line soliton. 
It appears that the perturbations lead eventually to a KdV soliton of slightly 
higher mass. As discussed in \cite{KS2}, different types of 
perturbation all indicate the stability of the KdV soliton for KP II.

It was shown in \cite{RT1, RT2, RT3, Za} that the KdV soliton is nonlinearly unstable 
against transversal perturbations in the KP I setting if its mass is 
above a critical value. The proof in \cite{Za} relies on the integrability of the KP I equation, but the methods in  \cite{RT1,RT2,RT3} apply to general dispersive equations.

\begin{figure}
[!htbp]
\begin{center}
\includegraphics[width=\textwidth]{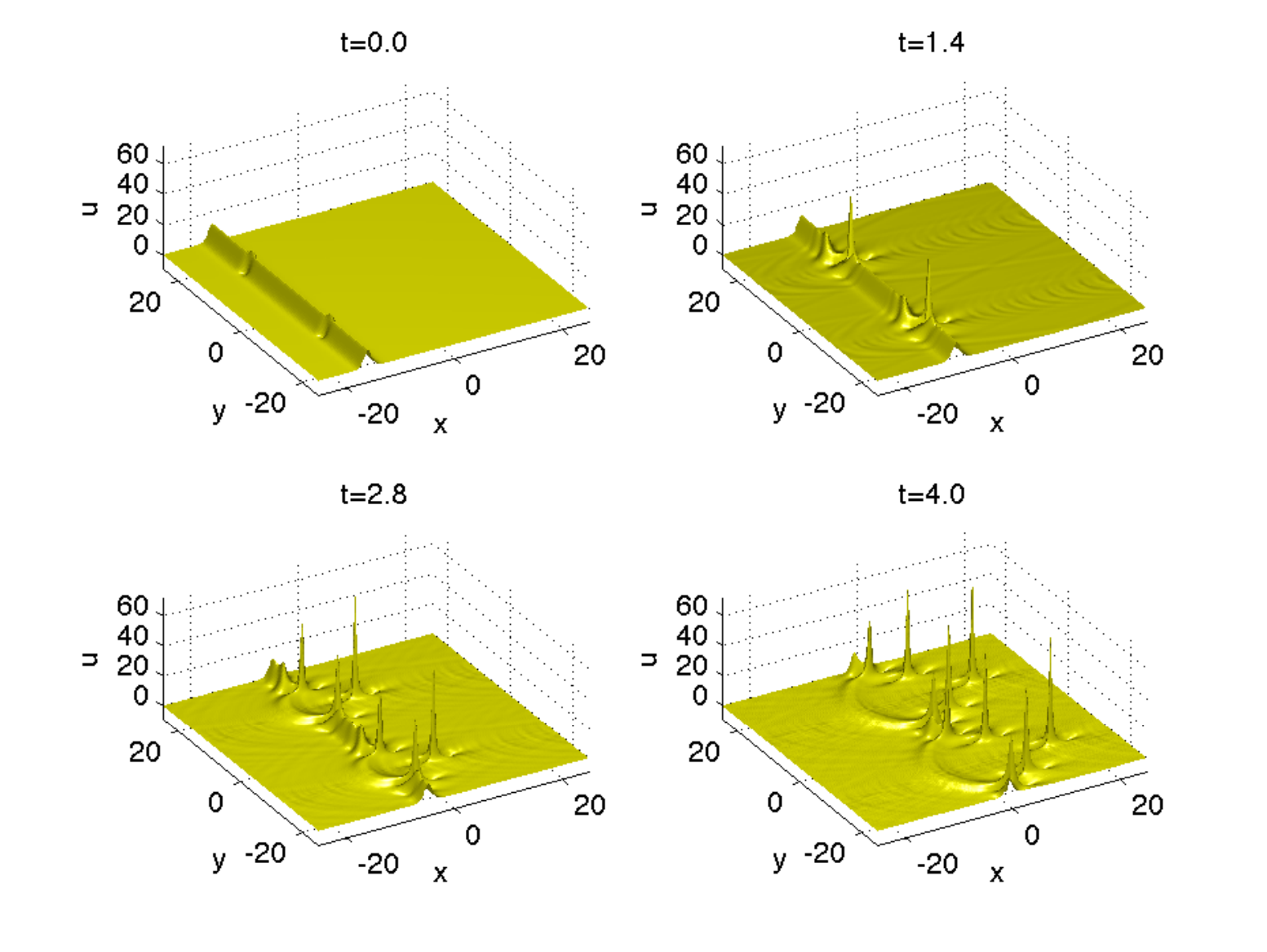} 
\caption{Solution to the  KP I equation  for initial 
data given by the KdV soliton $u_{sol}(x+2L_{x},0)$ 
plus perturbation $u_{p}=6(x-x_{1})\exp(-(x-x_{1})^{2})\left(\exp(-(y+L_{y}\pi/2)^{2})
    +\exp(-(y-L_{y}\pi/2)^{2})\right)$, $x_{1}=-2L_{x}$,  for various values of $t$.}
\label{kpIsolper}
\end{center}
\end{figure}

However, the type of the instability is 
unknown. Therefore in \cite{KS2}, this question was addressed 
numerically. 
In Fig.~\ref{kpIsolper} we show the KP I solution for the perturbed initial data of a 
line soliton with the perturbation (\ref{pert}) and 
$x_{1}=x_{0}=-2L_{x}$, 
i.e., the same setting as studied in Fig.~\ref{kpIIsolper} for KP II. Here the 
initial perturbations develop into 2 lumps which are traveling with 
higher speed than the line soliton. The formation of these lumps 
essentially destroys the line soliton which leads to the formation of 
further lumps. It appears plausible that for sufficiently long times 
one would only be able to observe lumps and small perturbations which 
will be radiated to infinity if studied on $\mathbb{R}^{2}$. 

\begin{figure}
[!htbp]
\begin{center}
\includegraphics[width=0.8\textwidth]{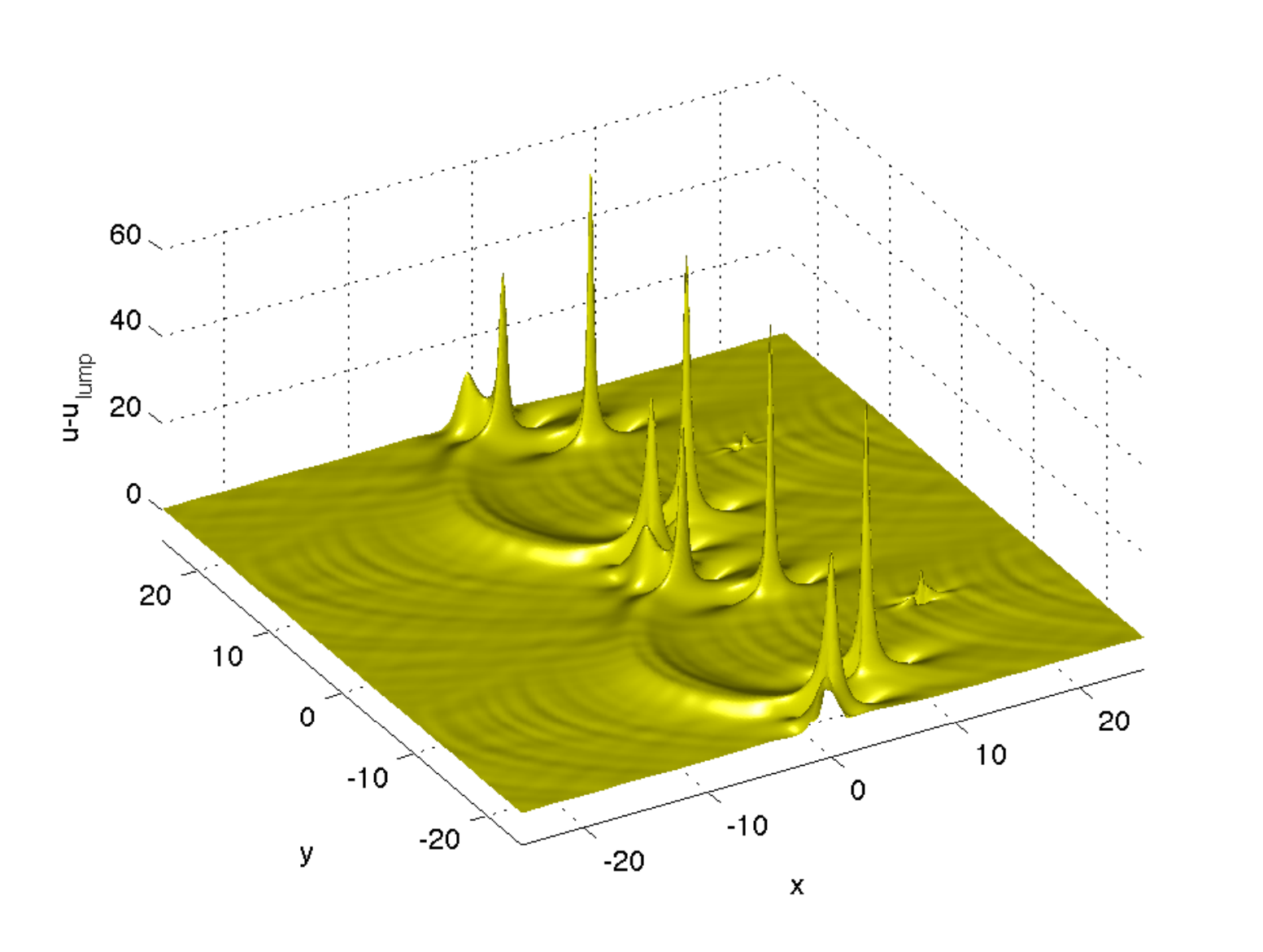} 
\caption{Difference of the solution to the  KP I equation  in 
Fig.~\ref{kpIsolper} for $t=4$ and two lump solutions fitted at the 
peaks farthest to the right. Only very small peaks remain of these 
`lumps' indicating that they develop asymptotically into true lumps.}
\label{kpIsolperl}
\end{center}
\end{figure}

We can give some numerical evidence for the validity of the 
interpretation of the peaks in Fig.~\ref{kpIsolper} as lumps in an 
asymptotic sense. We can identify numerically a certain peak, i.e., 
obtain the value and the location of its maximum. With these 
parameters one can study the difference between the KP solution and 
a lump with these parameters to see how well the lump fits the peak. 
This is illustrated for the two peaks, which formed first and which 
have therefore traveled the largest distance in Fig.~\ref{kpIsolperl}. 
A multi-lump 
solution might be a better fit, but here we mainly want to illustrate 
the concept which obviously cannot fully apply at  the studied small 
times. Nonetheless Fig.~\ref{kpIsolperl} illustrates convincingly 
that the observed peaks will asymptotically develop into lumps.

\begin{figure}
[!htbp]
\begin{center}
\includegraphics[width=\textwidth]{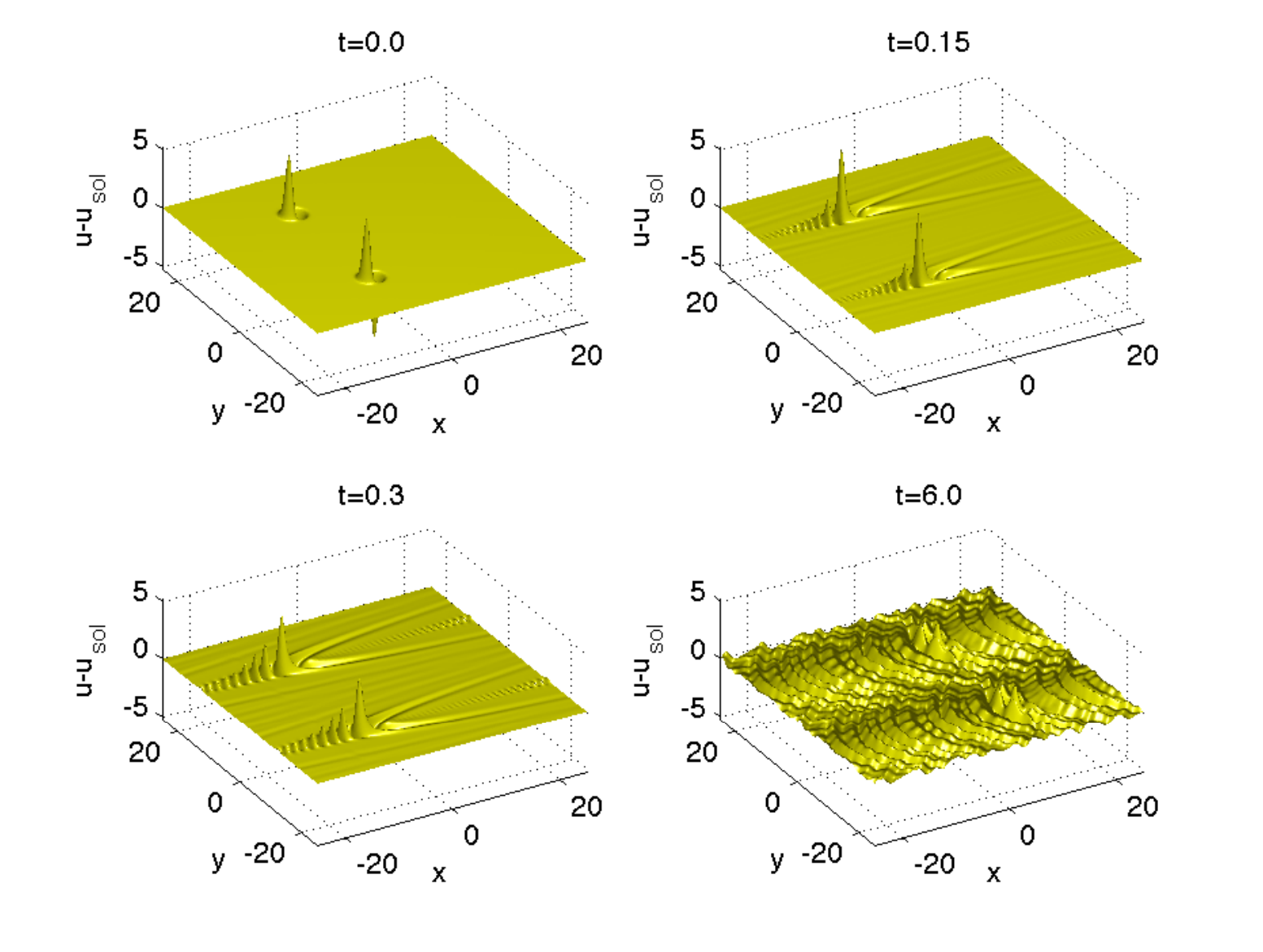} 
\caption{Difference of the solution to the  KP I equation  for initial 
data given by the KdV soliton $u_{sol}(x+2L_{x},0)$ 
plus perturbation $u_{p}=6(x-x_{1})\exp(-(x-x_{1})^{2})\left(\exp(-(y+L_{y}\pi/2)^{2})
    +\exp(-(y-L_{y}\pi/2)^{2})\right)$, $x_{1}=-L_{x}$ and the 
KdV soliton for various values of $t$.}
\label{kpIsolperd}
\end{center}
\end{figure}

Further examples for the nonlinear instability of the KdV soliton in 
this setting are given in \cite{KS2}. However, the nonlinear stability 
discussed in \cite{RT1, RT2} can be seen in 
Fig.~\ref{kpIsolperd} where the KP I solution is given for 
a perturbation of the line soliton as before 
$u_{0}(x,y)=12\mbox{sech}^{2}(x+2L_{x})+u_{p}$, but this 
time with $x_{1}=x_{0}/2$, i.e., perturbation and soliton are well 
separated.
The figure shows the difference between KP I solution and line 
soliton. It can be seen that the soliton is essentially stable on the shown time 
scales. The perturbation leads to algebraic tails towards 
positive $x$-values and to dispersive oscillations as studied in 
\cite{KSM2}. Due to the imposed periodicity both of these cannot 
escape the computational domain and appear on the respective other 
side. The important point is, however, that though the oscillations 
of comparatively large amplitude hit the line soliton quickly after 
the initial time, its shape is more or less unaffected till $t=6$.  
The KdV soliton eventually decomposes into lumps once it comes close 
to the boundaries of the computational domain, but this appears to be 
a numerical effect related to the imposed periodicity in $x$.

\section{The Davey-Stewartson systems}

The Davey-Stewartson (DS) systems are derived as asymptotic models in the so-called {\it modulation regime}  from various physical situations (water waves, plasma physics, ferromagnetism, see \cite {DS, AS, Co, Le}). They provide also
a good approximate solution to general quadratic hyperbolic systems using diffractive geometric optics \cite{Co,CoLa}.
They have the general form, where $a, b, c, \nu_1, \nu_2$ are real parameters depending on the physical context

\begin{equation}\label{DSgen}
\begin{array}{ccc}
i\partial_{t}\psi +a\partial_x^2\psi+b\partial_y^2 \psi=(\nu_1|\psi|^2+\nu_2\partial_x\phi)\psi,\\
\partial_x^2\phi+c\partial_y^2\phi=-\delta\partial_x|\psi|^2,
\end{array}
\end{equation}
where one can assume (up to a change of unknown) $a>0$ and $\delta>0.$

Using the terminology of \cite{GS}, one says that \eqref{DSgen} is

\begin{center}

elliptic-elliptic if \quad ($\text{sgn}\; b, \text{sgn}\; c)=(+1,+1),$

hyperbolic-elliptic  if \quad ($\text{sgn}\; b, \text{sgn}\; c)=(-1,+1),$

elliptic-hyperbolic  if \quad ($\text{sgn}\; b, \text{sgn}\; c)=(+1,-1),$

hyperbolic-hyperbolic  if \quad ($\text{sgn}\; b, \text{sgn}\; c)=(-1,-1).$

\end{center}

It is worth noticing that the Davey-Stewartson systems are "degenerate" versions of a more general class of systems describing the interaction of short and long waves, the Benney-Roskes, Zakharov-Rubenchik systems (\cite{BR, ZR}). None of those systems is known to be integrable.

It turns out that a very special case of the hyperbolic-elliptic and the  elliptic-hyperbolic DS systems are completely integrable. They are then classically known respectively as the DS II and DS I systems. Since we want to compare IST and PDE methods, we will focus on the hyperbolic-elliptic and elliptic-hyperbolic cases (referred to as DS II type and DS I type). We refer to \cite{GS, Ci, Ci2} for results on the elliptic-elliptic DS systems.

We will from now on write the DS II system in the form

\begin{equation}
    \label{DSII}
\begin{array}{ccc}
i 
\partial_{t}\psi+\partial_{xx}\psi-\partial_{yy}\psi+2\rho\left(\beta\Phi+\left|\psi\right|^{2}\right)\psi & = & 0,\quad \psi :  \R^2\times \R \to \C,
\\
\partial_{xx}\Phi+\partial_{yy}\Phi+2\partial_{xx}\left|\psi\right|^{2} & = & 0, \quad \Phi :  \R^2\times \R \to \R,\\
\psi(.,0)=\psi_0,
\end{array}
\end{equation}
where  the integrable DS II system corresponds to $\beta=1$ and $\rho$ takes the values $-1$ (focusing) and $1$ (defocusing).

The DS II system can be viewed as a nonlocal cubic nonlinear Schr\"{o}dinger equation. Actually one can solve $\Phi$ as

$$\Phi=2\lbrack(-\Delta)^{-1}\partial_{xx}\rbrack|\psi|^2,$$

where $(-\Delta)^{-1}\partial_{xx}=R_1^2$ is a zero order operator with Fourier symbol $-\frac{\xi_1^2}{|\xi|^2}$ and is thus bounded in all $L^p(\R^2)$ spaces, $1<p<\infty$ and all  Sobolev spaces $H^s(\R^2),$ allowing to write \eqref{DSII} as

\begin{equation}\label{DDSIIeq}
i 
\partial_{t}\psi+\partial_{xx}\psi-\partial_{yy}\psi+2\rho\left(2\beta R_1^2 (|\psi|^2)+\left|\psi\right|^{2}\right)\psi =  0.
\end{equation}

One easily finds that  (\ref{DSII}) has two formal conservation laws, the $L^2$ norm

$$ \int_{\R^2}|\psi(x,y,t)|^2dx dy= \int_{\R^2}|\psi(x,y,0)|^2dx dy$$

and the energy (Hamiltonian)
\begin{eqnarray}
    E(\psi(t))&=&\int_{\R^2}\left[ 
|\partial_{x}\psi|^2-|\partial_{y}\psi|^2-\rho (|\psi|^2+ 
\beta\Phi)|\psi|^2)\right] dxdy\nonumber\\
&=&E(\psi(0)).
    \label{energy}
\end{eqnarray}

Note that the integrable case $\beta=1$ is 
distinguished by the fact that the same hyperbolic operator appears 
in the linear and in the nonlinear part. In this case 
the equation is invariant under the transformation $x\to y$ and $\psi 
\to \bar{\psi}$ and \eqref{DDSIIeq} can be written in a "symmetric" form as

\begin{equation}\label{box}
i \partial_{t}\psi+\Box\psi-2\rho\lbrack(\Delta^{-1}\Box) |\psi|^2\rbrack\psi=0,
\end {equation}

where $\Box =\partial_{xx}-\partial_{yy}.$

This extra symmetry in the integrable case could be responsible for properties
 (existence of localized lump solutions and to  blow-up phenomena in case of DS II and existence of coherent "dromion" structures in case of DS I) that exist in the integrable case and might not persist in the non integrable cases as the following discussion will suggest.
 
 We summarize now some issues discussed in \cite{KS3} where one can also find many numerical simulations.
 
 \subsection{DS II type systems}
 
 Systems \eqref{DDSIIeq} (whatever the value of $\rho$ or $\beta$) can be seen as nonlocal variants of the {\it hyperbolic nonlinear Schr\"{o}dinger equation}
 \begin{equation}\label{HNLS}
 i\psi_t+\psi_{xx}-\psi_{yy}\pm |\psi|^2\psi=0, \quad \psi(\cdot, 0)=\psi_0
 \end{equation}

 and actually one can obtain (using Strichartz estimates in the Duhamel formulation) exactly the same results concerning the Cauchy problem (see \cite{GS}). Namely the Cauchy problem for \eqref{DDSIIeq} is locally well-posed for initial data $\psi_0$ in $L^2(\R^2)$ or $H^1(\R^2),$ and {\it globally} if $|\psi_0|_2$ is small enough. Nevertheless, since the existence time does not depend only on $|\psi_0|_2$ but on $\psi_0$ in a more complicated way, one cannot infer from the conservation of the $L^2$ norm that the $L^2$ solution is  a global one. Actually, proving (or disproving) the global well-posedness of the Cauchy problem for \eqref{HNLS} is an outstanding open problem.
 
As for the KP II equation, the inverse scattering problem for the (integrable) DS II is a $\bar{\partial}$ problem.

 It turns out that in the integrable case ($\beta=1$) inverse scattering techniques provide far reaching results which seem out of reach of purely PDE methods.
 
 In particular, Sung (\cite{Su1, Su2, Su3, Su4}) has proven the following

\begin{theorem}\label{DSSung}
Assume that $\beta =1.$ Let $\psi_0\in \mathcal S(\R^2).$ Then \eqref {DSII} possesses a unique global solution $\psi$ such that the mapping $t\mapsto \psi(\cdot,t)$ belongs to $C^\infty(\R, \mathcal S(\R^2))$ in the two cases:

(i) Defocusing.

(ii) Focusing and $|\widehat{\psi_0}|_1|\widehat{\psi_0}|_\infty<C,$ where $C$ is an explicit constant.

Moreover, there exists $c_{\psi_0}>0$ such that

$$|\psi(x,t)|\leq \frac{c_{\psi_0}}{|t|}, \quad (x,t)\in \R^2\times \R^*.$$
\end{theorem}


\begin{remark}
1. Sung obtains in fact the global well-posedness (without the decay rate) in the defocusing case under the assumption that $\hat \psi_0\in L^1(\R^2)\cap L^\infty(\R^2)$ and $\psi_0\in L^p(\R^2)$ for some $p\in [1,2),$ see \cite{Su4}.

2. Recently, Perry \cite{Pe} has given a more precise asymptotic behavior in the defocusing case for initial data in $H^{1,1}(\R^2)=\lbrace f\in L^2(\R^2) \;\text{such that}\; \nabla f, (1+|\cdot|)f \in L^2(\R^2)\rbrace,$ proving that the solution obeys the asymptotic behavior in the $L^\infty(\R^2)$ norm: 

$$\psi(x,t)=u(\cdot,t)+o(t^{-1}),$$

where $u$ is the solution of the linearized problem. 

\end{remark}

On the other hand, using an {\it explicit} localized lump like 
solution (see below) and a pseudo-conformal transformation, Ozawa \cite{Oz} has proven that the integrable focusing   DS II system possesses an $L^2$ solution that blows up in finite time $T^*$. In fact the mass density $|\psi(.,t)|^2$ of the solution converges as $t\to T^*$ to a Dirac measure with total  mass $|\psi(.,t)|^2_2=|\psi_0|^2_2$  (a weak form of the conservation of the $L^2$ norm). Every regularity breaks down at  the blow-up point but the solution persists after the blow-up time and disperses in the sup norm when  $t\to \infty $ as $t^{-2}.$ The numerical simulations in \cite{KS} suggest that a blow-up in finite time may also happen for other initial data, {\it eg} a sufficiently large Gaussian. Other numerical simulations suggest that the finite time blow-up does not persist in the {\it non integrable}, $\beta \neq 1$ case, both for the defocusing and focusing cases.

The family of lump solutions (solitons) to the integrable focusing  
DS II system reads (\cite{APP, MZBM,AC})

\begin{equation}\label{lump}
 \psi(x,y,t) = 2c \frac{\exp \left( -2i(\xi x - \eta y + 
2(\xi^{2}-\eta^{2} )t)\right)}{|x + 4\xi t + i(y + 4\eta t) +
 z_{0}|^2+|c|^2}
\end{equation}
 where $(c,z_{0})\in \mathbb{C}^2$ and $(\xi,\eta)\in\mathbb{R}^2$ 
 are constants. The lump moves with constant 
 velocity $(-4\xi, -4\eta)$ and decays as $(x^2+y^2)^{-1}$ for 
 $x,y\to\infty$.
 
 As explained in \cite{APP}, there is a one-to-one correspondence 
 between the lumps and the pole of the matrix solution of the direct 
 scattering problem. It is shown formally in \cite{GK} and rigorously 
 in \cite{Ki} that the lump is {\it unstable} in the following sense. 
 The soliton structure of the scattering data is unstable with 
 respect to a small compactly supported perturbation of the 
 soliton-like potential. It was also proven in \cite{PeSu, PeSu2} 
 that the lump is spectrally unstable.

\begin{figure}[htb!]
\centering
\includegraphics[width=0.45\textwidth]{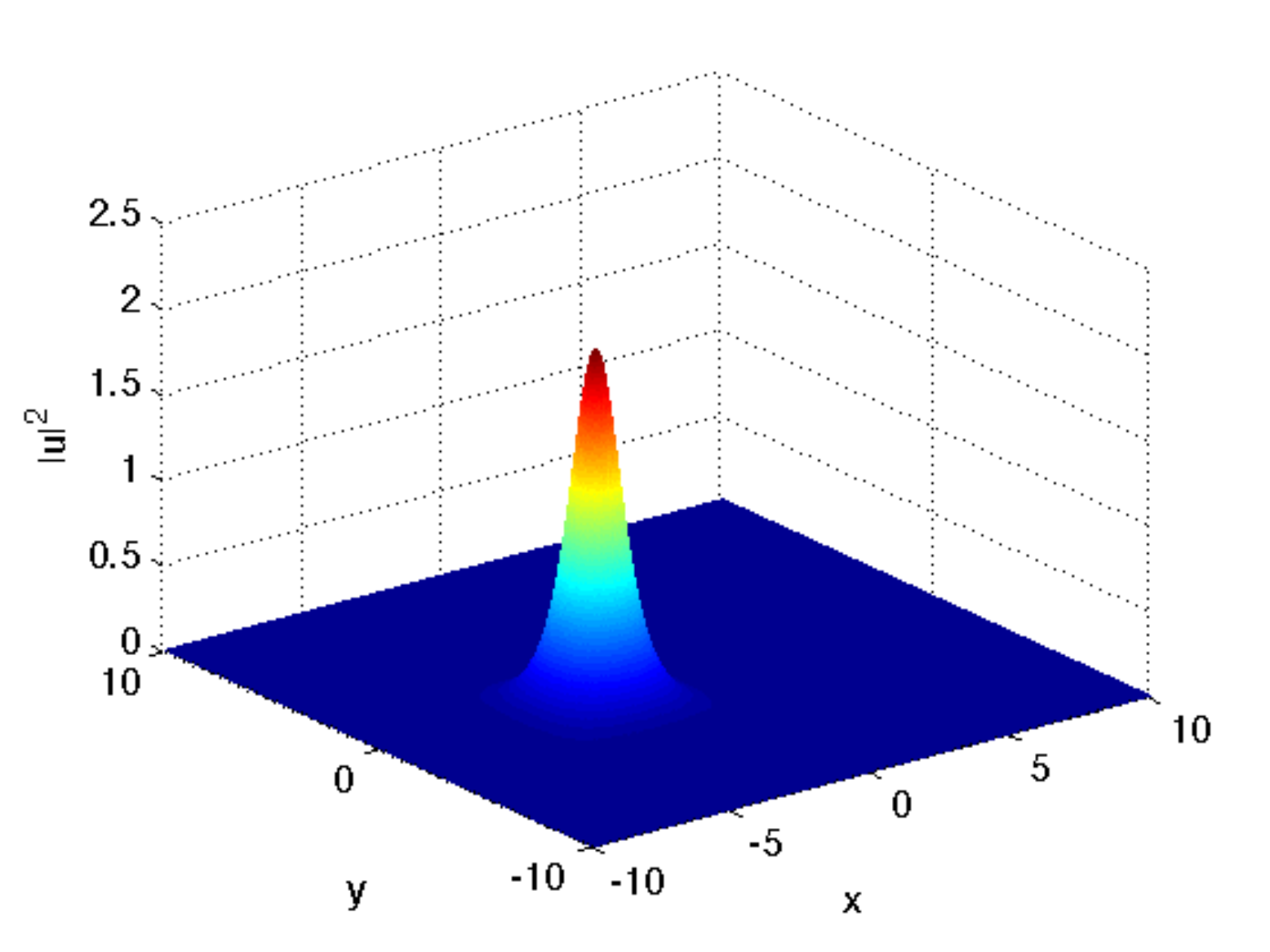}
\includegraphics[width=0.45\textwidth]{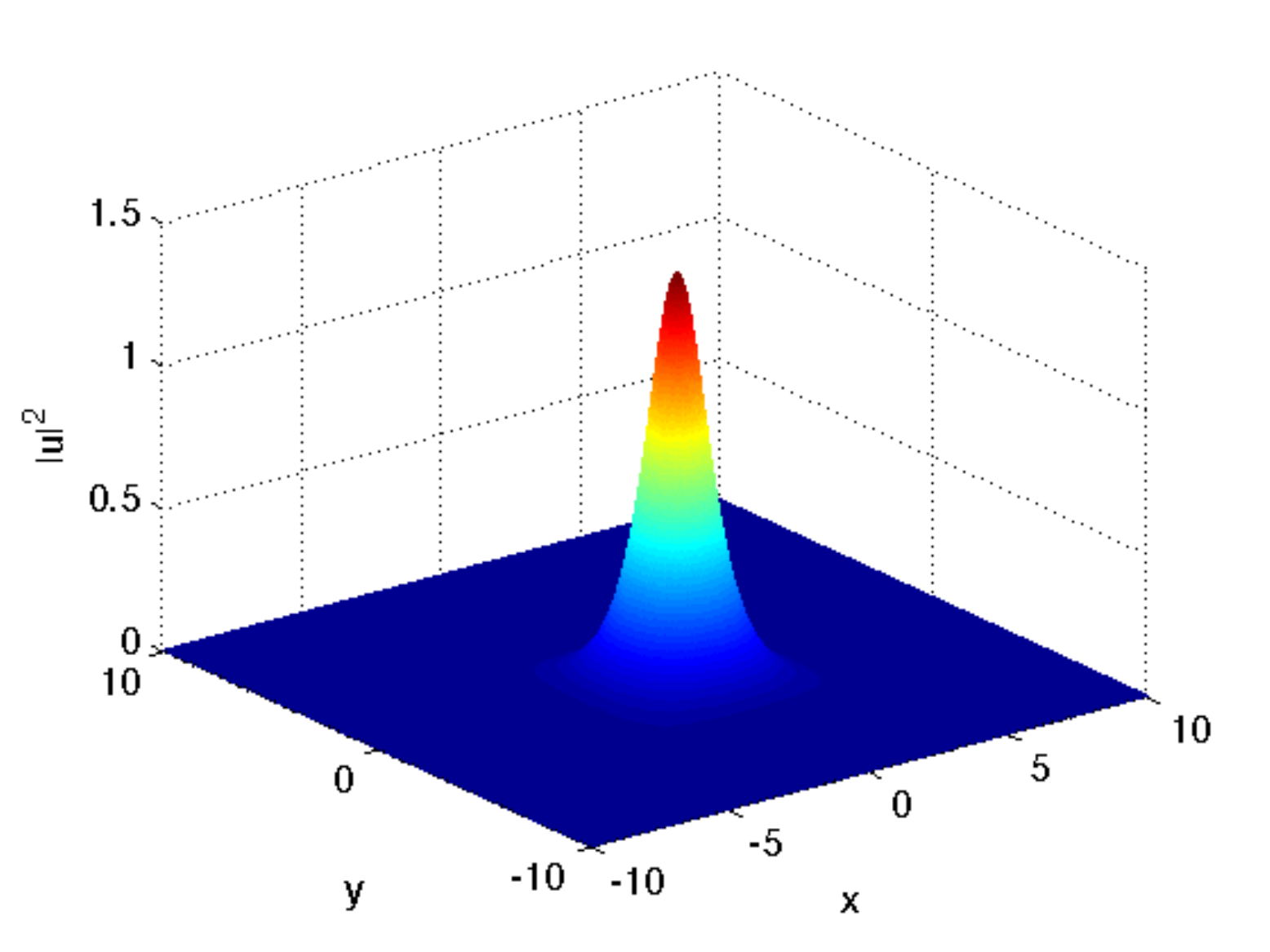}
\includegraphics[width=0.45\textwidth]{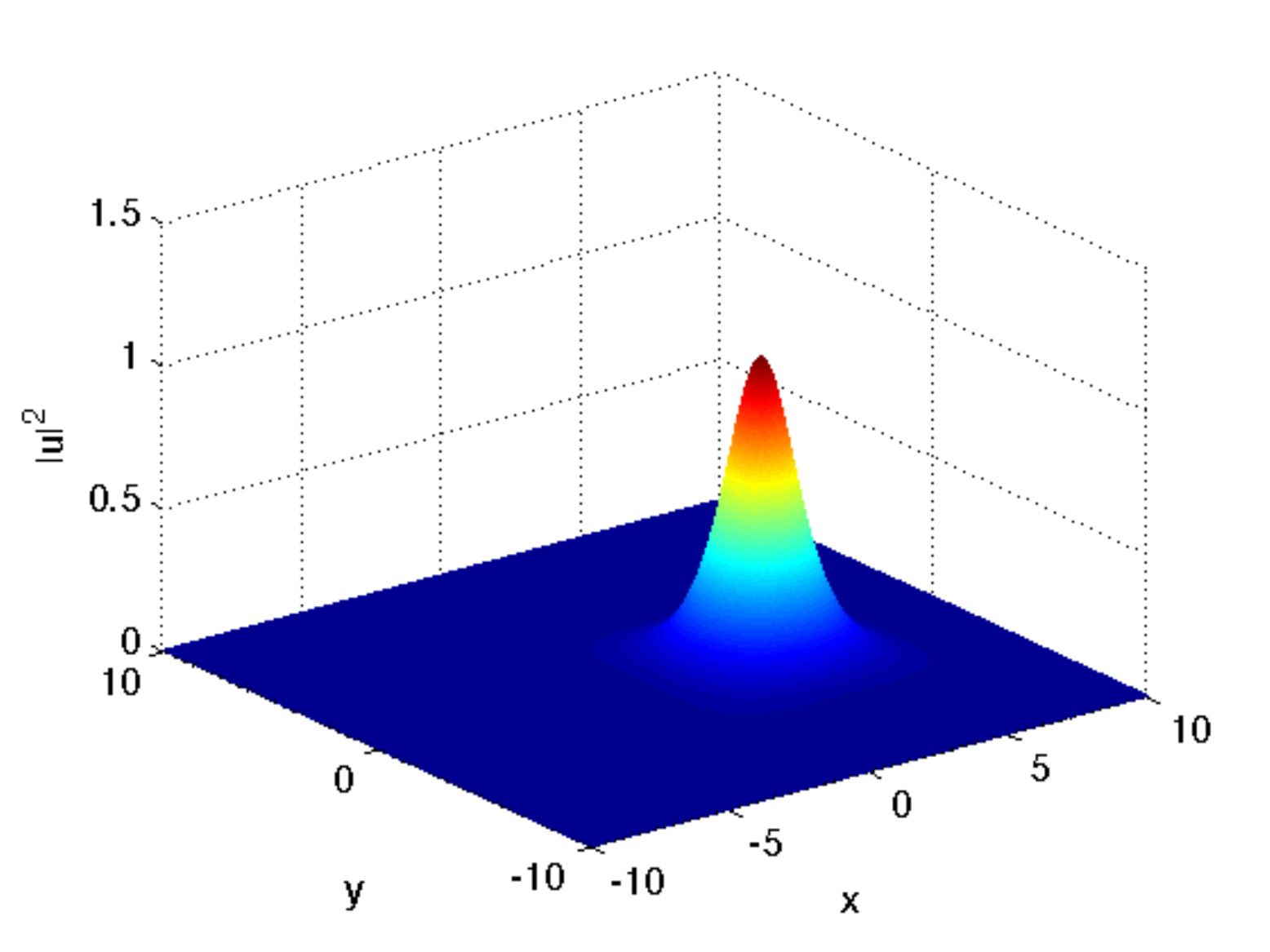}
\includegraphics[width=0.45\textwidth]{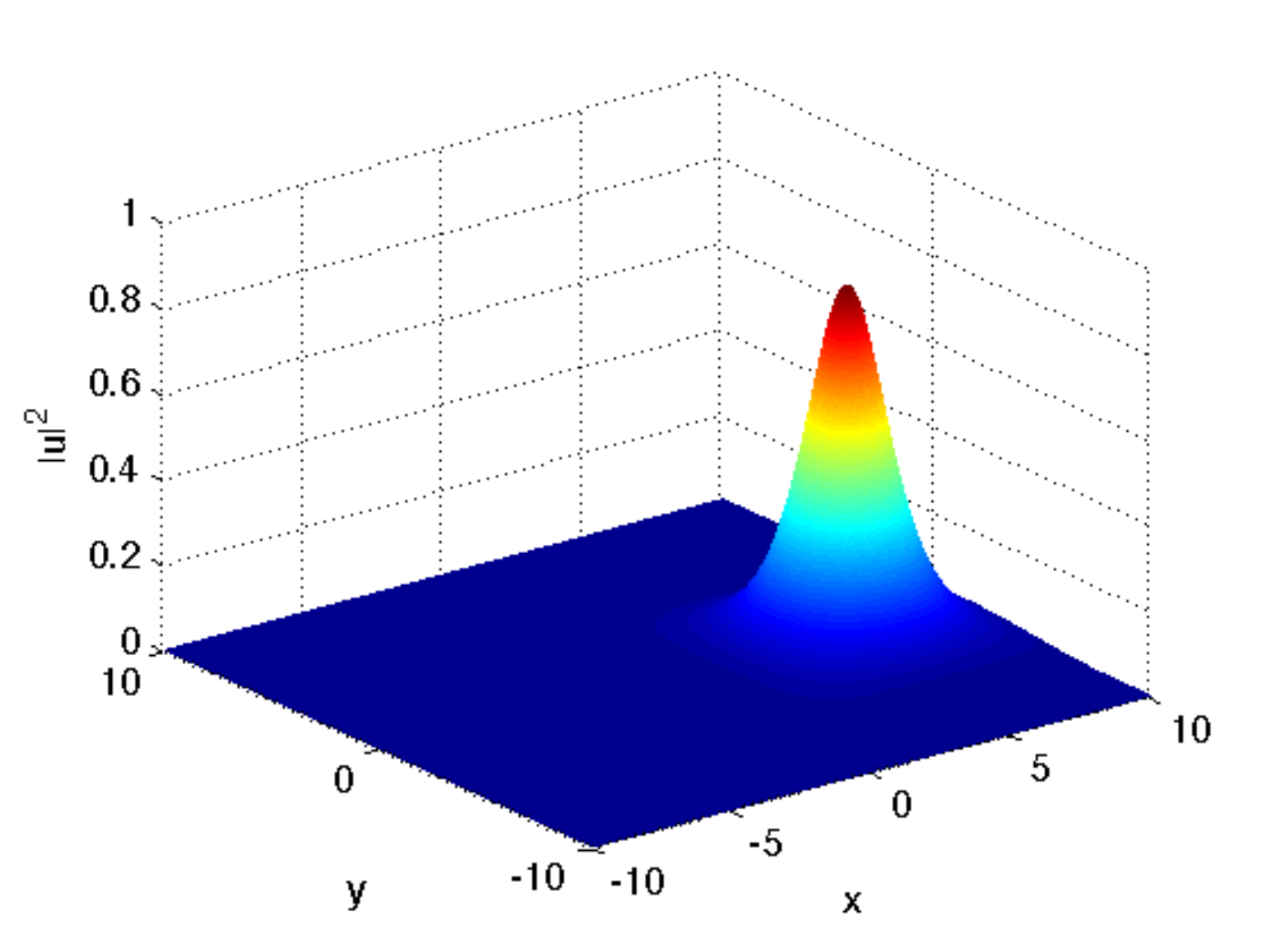}
\caption{Solution to the focusing DS II equation (\ref{DSII}) 
for an initial condition of the form $\psi(x,y,-3) = 0.9\psi_l$.} 
\label{ulupg1}
\end{figure}

The stability of the lump was 
 numerically studied in \cite{MFP,KMR}. It was shown that the lump is 
 both unstable against an $L^{\infty}$ blow-up in finite time and 
 against being dispersed away. In Fig.~\ref{ulupg1} taken from 
 \cite{KMR}, we consider an initial condition of the form 
$\psi(x,y,-3) = 0.9\psi_l,$ where $\psi_{l}$ is 
the lump solution (\ref{lump}) with $c=1$.
The solution travels at the same speed as before, but its amplitude varies, 
growing and decreasing successively.

\begin{figure}[htb!]
\centering
\includegraphics[width=0.45\textwidth]{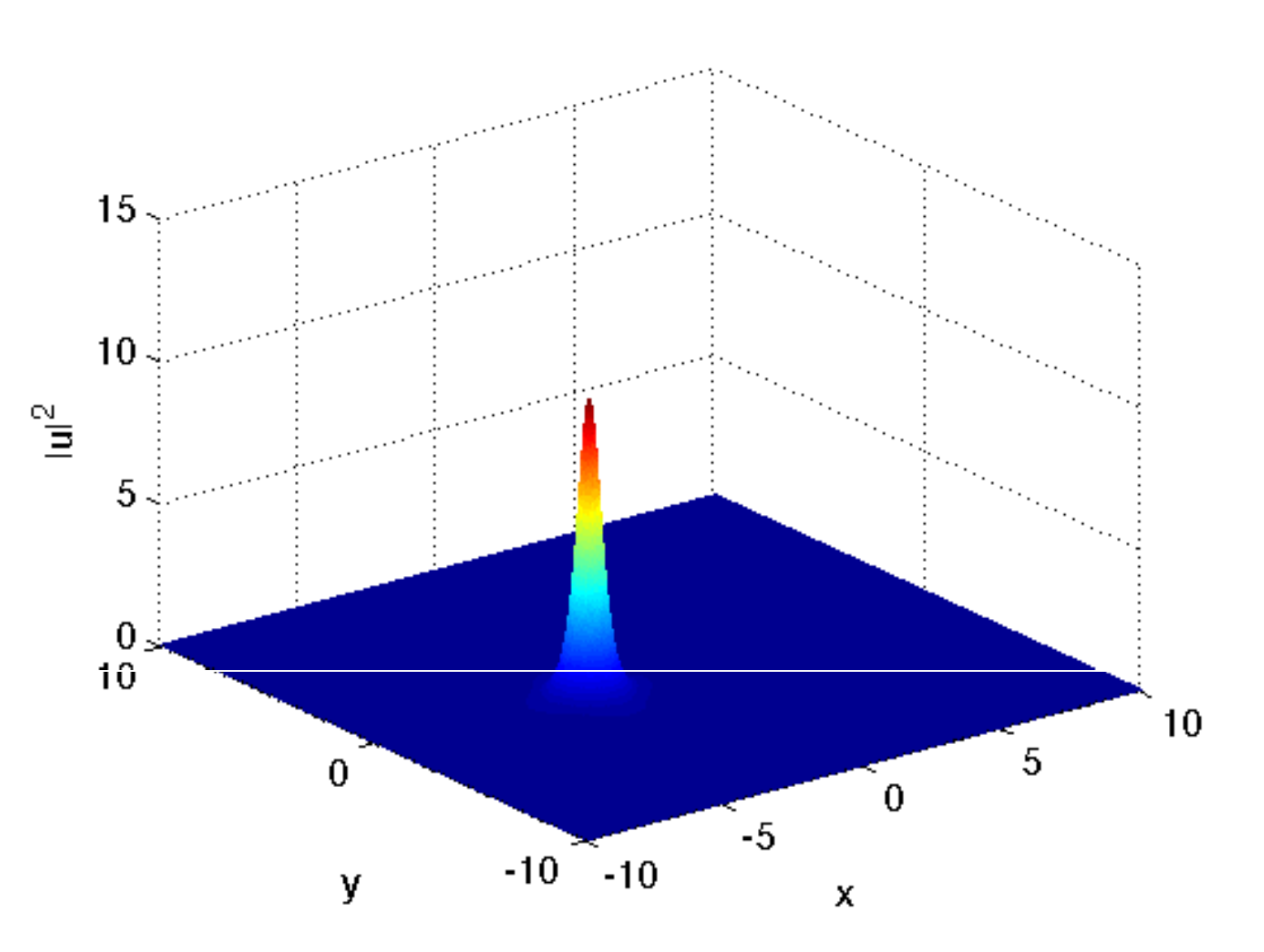}
\includegraphics[width=0.45\textwidth]{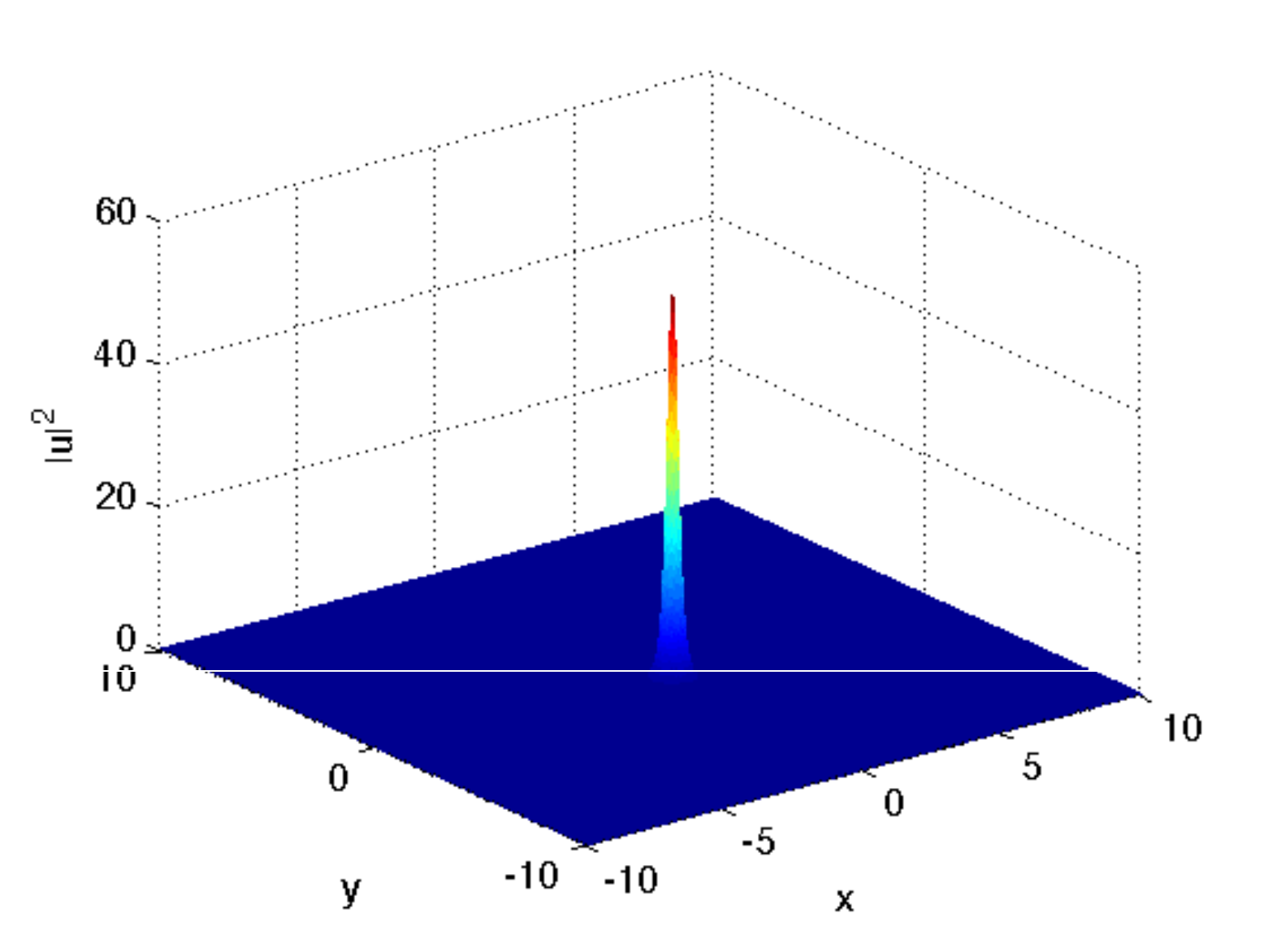}
\caption{Solution to the focusing DS II equation (\ref{DSII}) 
for an initial condition of the form $\psi(x,y,-3) = 1.1\psi_l$.} 
\label{ulupg11}
\end{figure}

If instead the initial data $\psi(x,y,-3) = 0.9\psi_l,$ is 
considered, the solution appears to blow up in finite time as can be 
seen in Fig.~\ref{ulupg11}.  Note that also the Ozawa 
solution \cite{Oz} was in \cite{KMR} numerically shown to be unstable 
against both an earlier blow-up and being dispersed away. In 
\cite{KS3} it was shown that the blow-up of the Ozawa solution is not 
generic. 

\subsubsection{Localized solitary waves of DS II type systems}
 It is well known that hyperbolic NLS equations such as \eqref{HNLS} do not 
possess solitary waves of the form $e^{i\omega t}\psi (x)$ where $\psi$ is localized (see \cite {GS1}).

It was furthermore proven in \cite{GS1} that  non trivial solitary waves may exist for DS II type systems only when $\rho =-1$ (focusing case) and $\beta\in (0,2).$ Note that the (focusing) integrable case corresponds to $\beta =1.$ Moreover  solitary waves with  {\it radial} (up to translation) profiles can exist only when $\rho=-1$ and $\beta=1$, that is in the focusing integrable case.

Those results  (and  the numerical simulations in \cite{KMR, KS}) suggest that localized solitary waves for the focusing DS II systems exist only in the integrable case and this might be due to the new symmetry of the system we were alluding to above in this case. 

To summarize, one is led to conjecture that neither the existence of 
the lump nor the associated Ozawa blow-up persist in the focusing DS II non integrable case. One can also conjecure that the solution of DS II type solution is global and decays in the sup norm as $1/t$ as was shown by Sung and Perry in the  (integrable) DS II case.

\begin{remark}\label{nonloc} The previous results and conjectures do not of course exclude the existence of nontrivial, non localized,
traveling waves. In the context of the hyperbolic cubic NLS, see  for instance Remark 2.1 in \cite{GS1} or \cite{KNZ}. We do not know of existence of  similar solutions for DS II type systems.
\end{remark}


 \subsection{DS I type systems}

 The DS I type systems are  quite different from the other DS systems. Actually, solving the hyperbolic equation for $\phi$
 (with suitable conditions at infinity) yields a {\it loss of one derivative} in the nonlinear term and the resulting NLS type equation is no more semilinear.
 Even proving the rigorous conservation of the Hamiltonian leads to serious problems. We describe now how to solve the equation for $\phi$ in  a $H^1$ framework (see \cite{GS}).

The  elliptic-hyperbolic DS system can be written after scaling as 

\begin{equation*}\label{DSEH}\left\{\begin{aligned}
 &i\partial_t \psi +\Delta \psi=\chi |\psi|^2\psi+b\psi\phi_x\\
 &\phi_{xx}-c^2\phi_{yy}=\frac{\partial}{\partial x}|\psi|^2.
 \end{aligned}\right.
 \end{equation*}

 The integrable DS I system corresponds to $\chi+\frac{b}{2}=0.$

\vspace{0.3cm}
 We now solve the equation for $\phi.$

 Let $c>0.$  Consider the equation 
\begin{equation}\label{ondes}
\frac{\partial^2\phi}{\partial x^2}-c^2\frac{\partial^2\phi}{\partial y^2}=f\quad \text{in}\; \R^2
\end{equation}
with the boundary condition
\begin{equation}\label{bounhyp}
\lim_{\xi\to+\infty} \phi(x,y)=\lim_{\eta\to +\infty}\phi(x,y)=0
\end{equation}
where $\xi=cx-y$ and $\eta= cx+y.$

Let $K=K_c$ the kernel
$$K(x,y;x_1,y_1)=\frac{1}{2}H(c(x_1-x)+y-y_1)H(c(x_1-x)+y_1-y)$$
where $H$ is the usual Heaviside function.

\begin{lemma}\label {hypphi}(\cite{GS}). 
For every $f\in L^1(\R^2),$ the function $\phi=\mathcal K(f)$ defined by
\begin{equation}\label{solvphi}
\phi(x,y)=\int_{\R^2}K(x,y;x_1,y_1)f(x_1,y_1)dx_1dy_1
\end{equation}
is continuous on $\R^2$ and satisfies \eqref{ondes} in the sense of distributions.
Moreover,
$\phi\in L^\infty(\R^2), (\partial \phi/\partial x)^2-c^2(\partial \phi/\partial y)^2\in L^1(\R^2)$ and we have the following estimates
\begin{equation}\label{est1}
\sup_{(x,y)\in \R^2} |\phi(x,y)|\leq \int_{\R^2}|f| dxdy
\end{equation}
\begin{equation}\label{est2}
\int_{\R^2} \left|\left(\frac{\partial \phi}{\partial x}\right)^2-c^2\left(\frac{\partial \phi}{\partial y}\right)^2\right|dx dy\leq \frac{1}{2c}\left(\int_{\R^2}|f| dx dy\right)^2.
\end{equation}
\end{lemma}

\begin{remark}

1. No condition is required as $\xi$ or $\eta$ tends to $-\infty.$
2.  In general, $\nabla \phi\notin L^2(\R^2)$ even if $f\in C_0^\infty(\R^2),$ but Lemma \ref{hypphi} allows to solve the $\phi$ equation as soon as $\psi\in H^1(\R^2)$ for instance.

\end{remark}

The DS I type system  possesses the formal Hamiltonian 
\begin{equation*}\label{Hamil}
E(t)=\int_{\R^2}\left[|\nabla 
\psi|^2+\frac{\chi}{2}|\psi|^4+\frac{b}{2}(\phi_x^2-c^2\phi_y^2)\right]dxdy.
\end{equation*}
Lemma \ref{hypphi} allows to prove that this Hamiltonian makes sense in an $H^1$ setting for $\psi$ (see (GS). Proving its conservation  on the time interval of the solution is an open problem as far as we know (this would lead to global existence of a weak $H^1$ solution).

\subsubsection{DS I type by PDE methods}

The first local well-posedness result is due to Linares-Ponce \cite{LiPo2} and we summarize below the best known results, due to Hayashi-Hirata \cite{NH1, NH2} and Hayashi \cite{H}.

After rotation, one can write the DS-I type systems as

\begin{equation*}\label{DSI}\left\{\begin{aligned}
 &i\partial_t \psi +\Delta \psi=i(c_1+\frac{c_2}{2}) |\psi|^2\psi-\frac{c_2}{4}\left(\int_x^\infty\partial_y|\psi|^2dx'+\int_y^\infty \partial_x|\psi|^2dy'\right) \psi\\
 &+\frac{c_2}{\sqrt 2}\left((\partial_x \phi_1)+\partial_y \phi_2)\right)\psi,
 \end{aligned}\right.
 \end{equation*}
 
 where $c_1,c_2\in \R$ and $\phi$ satisfies the radiation conditions
 
 $$\lim_{y\to \infty}\phi(x,y,t)=\phi_1(x,t),\quad \lim_{x\to \infty}\phi(x,y,t)=\phi_2(y,t).$$
 
 \vspace{0.5cm}

\begin{theorem}
(\cite{H})

Assume
$\psi_0\in H^2(\R^2),\; \phi_1\in C(\R;H^2_x)\;\text{and}\; \phi_2\in C(\R;H^2_y).$ Then there exist  $T>0$ and a unique solution $\psi\in C(\lbrack 0,T\rbrack;H^1)\cap L^\infty(0,T;H^2)$ with initial data $\psi_0.$
\end{theorem}

\begin{itemize}
\item The proof uses in a crucial way the  {\it smoothing properties} of the Schr\"{o}dinger group.
\end{itemize}

The next result concerns global existence and scattering of small solutions  in the weighted Sobolev space

\ $$H^{m,l}=\lbrace f\in L^2(\R^2); |(1-\partial_x^2-\partial_y^2)^{m/2}(1+x^2+y^2)^{l/2} f|_{L^2}<\infty\rbrace.$$

\begin{theorem} (\cite{NH2}) Let $\psi_0\in H^{3,0}\cap H^{0,3},\; \partial_x^{j+1}\phi_1\in C(\R;L^\infty_x),\; \partial_y^{j+1}\phi_2\in C(\R;L^\infty_y),\; (0\leq j\leq 3),$
"small enough". Then

\begin{itemize}
\item There exists a solution $\psi\in L^\infty_{\text {loc}}(\R;  H^{3,0}\cap H^{0,3})\cap C(\R; H^{2,0}\cap H^{0,2}).$
\item Moreover
$$||\psi(\cdot,t)||_{L^\infty}\leq C(1+|t|)^{-1}(||\psi||_{H^{3,0}}+||\psi||_{H^{0,3}}).$$

There exist $ u^{\pm}$ such that

$$||\psi(t)-U(t)u^{\pm}||_{H^{2,0}}\to 0,\quad \text{as}\; \to \pm \infty.$$
\end{itemize}

where $U(t)=e^{it(\partial_x^2+\partial_y^2}).$
\end{theorem}

\subsubsection{DS I by IST. Comparison with elliptic-hyperbolic DS}

Contrary to the results of the previous subsection which were valid for arbitrary values of $\beta$ we focus here on the integrable case, $\beta =1.$

The first set of results concerns coherent structures (dromions) for DS I with nontrivial conditions (on $\phi$) at infinity. The existence of dromions is established in \cite{BLMP, FoSa} and the perturbations of the dromion are investigated in \cite {Ki3} (see also Section 7 of \cite{Ki}).

 We do not know of any study of dromions by PDE techniques or of existence of similar structures in the non-integrable case. Actually they might have no physical relevance.

Concerning the Cauchy problem, the global existence and uniqueness of a solution $\psi\in C(\R;\mathcal S(\R^2))$ of DS I for data $\psi_0\in \mathcal S(\R^2),\; \phi_1, \phi_2\in C(\R;\mathcal S(\R))$ is proven in \cite{FoSu2}. Under a smallness condition, the solutions with trivial boundary conditions $\phi_1=\phi_2=0$ disperse as $1/t$ (Kiselev \cite{Ki3}, see also \cite{Ki}). A  precise asymptotics is also given.
 The numerical simulations in \cite{BB} for general DS I type systems confirm the dispersion of solutions of DS I with trivial boundary conditions and 
suggest that the dromion is not stable with respect to the coefficients, that is it does not persist in the non-integrable case.

\section{Final comments}
We briefly comment here on two other integrable equations.

\subsection{The Ishimori system}

The Ishimori systems were introduced in \cite{Ishi} as 
two-dimensional generalizations of the Heisenberg equation in 
ferromagnetism. They read 

\begin{equation}
    \label{Ishi I}
\begin{array}{ccc}
S_t=S\wedge(S_{xx}-S_{yy})+b(\phi_xS_y+\phi_yS_x)\\
\phi_{xx}+\phi_{yy}=2S\cdot(S_x\wedge S_y)\\
S(\cdot,0)=S_0, 
\end{array}
\end{equation}

\begin{equation}
    \label{Ishi II}
\begin{array}{ccc}
S_t=S\wedge(S_{xx}+S_{yy})+b(\phi_xS_y+\phi_yS_x)\\
\phi_{xx}-\phi_{yy}=-2S\cdot(S_x\wedge S_y)\\
S(\cdot,0)=S_0, 
\end{array}
\end{equation}

where $S$ is the spin,  $S(\cdot,t): \R^2\to \R^3,$ $|S|^2=1,$ $S\to (0,0,1)$ as $|(x,y)|\to \infty $ and $\wedge$ is the wedge product in $\R^3.$

The coupling potential $\phi$ is a scalar unknown related to the topological charge density $q(S)=2S\cdot(S_x\wedge S_y).$ $b$ is a real coupling constant. When $b=1$, \eqref{Ishi I}, \eqref{Ishi II} are completely integrable (\cite{KoMa, KoMa2}).

Note the formal analogy \footnote{which is clearer after the stereographic projection below.} of \eqref{Ishi I}, \eqref{Ishi II} with DS II and DS I type systems.

It is proven in \cite{So} that the Cauchy problem for \eqref{Ishi I} 
(for arbitrary values of $b$) is locally  well-posed in $H^m(\R^2), 
m\geq 3$ provided the initial spins are almost parallel. Under a 
stronger regularity assumption on the initial data it is furthermore 
proven that   the solution  is global and converges to a solution of the linear  "hyperbolic" Schr\"{o}dinger equation as t goes to infinity.

The idea is to reduce (\ref{Ishi I}), (\ref{Ishi II}) to a nonlinear (hyperbolic) Schr\"{o}dinger type equation by the stereographic projection

$$u: \R^2\to \C,\;u=\frac{S^1+iS^2}{1+S^3},$$

reducing \eqref{Ishi I} to

\begin{equation}
    \label{Ishi Ist}
    \left\lbrace
    \begin{array}{l}
iu_t+u_{xx}-u_{yy}=\frac{2 \bar{u}}{1+|u|^2}(u_x^2-u_y^2)+ib(\phi_xu_y+\phi_xu_x)\\
\Delta \phi=4i\frac{u_x\bar{u_y}-\bar{u_x}u_y}{(1+|u|^2)^2}\\
u(\cdot, 0)=u_0,
 \end{array}\right.
\end{equation}

with the condition $|u(x,y)|\to 0$ as $|(x,y)|\to \infty. $

The initial condition on $S$ becomes thus a smallness condition on $u_{\vert t=0}.$ 

The local well-posedness of the Cauchy problem for \eqref{Ishi I} for {\it arbitrary} initial data in $H^m(\R^2), m\geq 4$ is proven in \cite{KN}. This result is improved in \cite{BIK} where local well-posedness is proven for arbitrary large initial data in $H^s(\R^2), s>\frac{3}{2}$ having a range that avoids a neighborhood of the north pole.

Concerning the IST method for \eqref{Ishi I}, the rigorous justification of the procedure in \cite{KoMa, KoMa2} is not trivial and Sung \cite{Su5} used instead a gauge transform which relates  \eqref{Ishi I} (when $b=1$) to the focusing integrable DS II system. This allows to prove the global well-posedness of the Cauchy problem for small initial data (in a different functional setting than \cite{So}). 

The connection between the (integrable) Ishimori system \eqref{Ishi 
I} is nicely used in \cite{BIK} to prove the {\it global} well-posedness of the Cauchy problem for \eqref{Ishi I} in the {\it defocusing} case, that is when the target of $S$ is no more the sphere $\mathbb{S}^2$, but the hyperbolic space $\mathbb{H}^2=\lbrace (x,y,z)\in \R^3; x^2-y^2-z^2=1, x>0\rbrace.$ The gauge transform relates in this case the {\it defocusing} Ishimori system \eqref{Ishi I} to the {\it defocusing} DS II system. Such a result is not known in the non integrable case, $b\neq 1.$

We do not know of any finite time blow-up for \eqref{Ishi I} or of any rigorous result (by PDE or IST methods) for \eqref{Ishi II}.

\subsection{The Novikov-Veselov equation}

The Novikov-Veselov system

\begin{equation}
   \label{NVor}
    \left\lbrace
    \begin{array}{l}
    v_t=4\text{Re}\;(4\partial_z^3 v+\partial_z(vw)-E\partial_zw),\\
\partial_{\bar {z}}w=-3\partial_zv,\quad v=\bar{v}
\end{array}\right.
    \end{equation}
    
    where $\partial_z=\frac{1}{2}\left(\partial_x-i\partial_y\right),\quad \partial_{\bar {z}}=\frac{1}{2}\left(\partial_x+i\partial_y \right),$

\vspace{0.5cm}
was introduced in \cite{NV1, NV2} as a two dimensional  analog of the  Korteweg-de Vries equation, integrable via the inverse scattering transform for the following 2-dimensional stationary Schr\"{o}dinger equation
at a fixed energy $E$:

\begin{equation}\label{Sch}
L\psi =E\psi,\quad L=-\Delta +v(x,y,t),
\end{equation}

where $\Delta= \frac{\partial^2}{\partial x^2}+\frac{\partial^2}{\partial y^2}$ and $E$ is a fixed real constant.

The Novikov-Veselov equation has the Manakov triple representation \cite{Ma2}

\begin{equation}\label{Lax}
\dot{\mathcal L}=\lbrack A,\mathcal L \rbrack -B\mathcal L,
\end{equation}

where 

$$\mathcal L=-\Delta +v+E,$$

$$A=8(\partial^3+\bar{\partial^3})+2(w\partial+\bar{w}\bar{\partial}),$$

$$\bar{B}=-2(\partial w+\bar{\partial}\bar{w}),$$

that is $(v,w)$ solves \eqref{NV} if and only if \eqref{Lax} holds.

\begin{remark}
There is an interesting formal limit of \eqref{NVor} to KP I (resp. KP II) as $E\to+\infty$ (resp. $E\to -\infty)$, under an appropriate scaling, assuming tht the wavelengths in $y$ are much larger than those in $x$.
\end{remark}

\begin{remark}
As far as we know, and contrary to the integrable equations studied above, the Novikov-Veselov equation does not seem  be derived as an asymptotic model from a more general system.
\end{remark}

To make the equation more "PDE like", we first write the equation in a slightly different form.
    
  \vspace{0.3cm}  
Setting $w=w_1+iw_2,$ \eqref {NVor} becomes:

 \begin{equation}
    \label{NVbis}
    \left\lbrace
    \begin{array}{l}
    v_t=2v_{xxx}-6v_{xyy}+2[(vw_1)_x+(vw_2)_y
    ]-2E(\partial_x w_1+\partial_y w_2), \; (x,y)\in \R^2,\; t>0, \\
    \partial_xw_1-\partial_y w_2=-3v_x\\
    \partial_xw_2+\partial_y w_1=3v_y.
    \end{array}\right.
    \end{equation}
    
    Note that the last two equations imply
    
    $$(\partial_{xx}-\partial_{yy})w_2+2\partial_{xy}w_1=0$$
    
    and that $w_1$ et $w_2$ are defined up to an additive constant.
    
    One can express $w=(w_1,w_2)$ in dependence of $v$, in a unique way via the Fourier transform
    
    $$\hat{w_1}=\frac{3(\xi_2^2-\xi_1^2)}{|\xi|^2} \hat{v},\quad \hat{w_2}=\frac{6\xi_1\xi_2}{|\xi|^2} \hat{v}.$$

    We denote $L_1$ and $L_2$ the zero order corresponding operators,
    
    $$w_1=L_1v,\quad w_2=L_2 v.$$
    
    With these notations, the equation reads
    
    \begin{equation}\label{NV}
    v_t-2v_{xxx}+6v_{xyy}+2E(L_1v_x+L_2v_y)-2[(vL_1v)_x+(vL_2v)_y]=0.
    \end{equation}
    
    One remarks that the dispersive part in \eqref{NV} is reminiscent of that of the  Zakharov-Kuznetzov equation \cite{ZR}
    
    $$u_t+u_x+\Delta u_x+uu_x=0,$$
    
    another two-dimensional extension of the KdV equation, which is not integrable though. We refer to \cite{BAKS} for a  systematic study of the dispersive properties of general third order (local) operators in two-dimensions.
    
    \subsubsection{The zero energy case $E=0$.}
    
    We refer to the excellent survey \cite{CMMPSS} for a rather complete account of what is known in this case and which comprises some interesting numerics and a rich bibliography.
    
    We extract from \cite{CMMPSS} the following definition :

    \vspace{0.3cm}
        \begin{definition}The operator $L = -\Delta +v$ is said to be:
    
(i) subcritical if the operator $L$ has a positive Green's function and the equation
$L = 0$ has a strictly positive distributional solution,
    
(ii) critical if $L = 0$ has a bounded strictly positive solution but no positive Green's
function, and
    
(iii) supercritical otherwise.
    
  \end{definition}

  The following conjecture is made in \cite{CMMPSS}.
  
    \vspace{0.3cm}
    {\it Conjecture}. The Novikov-Veselov equation with zero energy has a global solution for
critical and subcritical initial data, but its solution may blow up in finite time for supercritical initial data.\footnote{By the result in \cite{Ang} mentioned below, the Novikov-Veselov equation with zero energy has a {\it  local} solution for a very general class of initial potentials.}
    
    We refer to \cite{CMMPSS} for some partial results toward its resolution.

    \subsubsection{The case $E\neq 0.$}

    For nonzero energy $E$ and potentials $v$ which vanish at infinity, the scattering transform and inverse scattering method was developed by P. Grinevich, R. G.
Novikov, and S.-P. Novikov (see Kazeykina's thesis \cite{Kaz} for an excellent survey and \cite{Gr2, GrMa, GrNo} for the original papers).
    
    We also refer to \cite{Kaz} and to the papers \cite{Kaz1, Kaz2, Kaz3, KaNo1, KaNo2, KaNo3,Nov} for many interesting results on the soliton solutions (absence, decay properties,..) and the long time asymptotics. More precisely it is proven that the rational, non singular solutions introduced in \cite{Gr2, Gr3} at positive energy are multi-solitons. It is also proven that solitons cannot decay faster than $O(|x|^{-3})$ when $E\neq 0$ and than $O(|x|^{-2})$ when $E=0.$ Finally, the evolution solutions corresponding to non singular scattering data are shown to decay asymptotically as $O(t^{-1})$ when $E>0$ and as $O(t^{-3/2})$ when $E<0.$

    \subsubsection{Novikov-Veselov by PDE techniques}
   
We only know of one result on Novikov-Veselov by PDE techniques. Angelopoulos \cite{Ang} has proven the local well-posedness in $H^{1/2}(\R^2)$ of the Cauchy problem in the zero energy case ($E=0$) case following the method used in \cite{MP2}.   It is very likely that the methods in \cite{KoSa} for the Zakharov-Kuznetsov equation \footnote{mainly based on the smoothing properties of the linear group generated by third order partial differential operators. Such a property is very likely to hold for the {\it nonlocal} operator $-2\partial_x^3+6\partial_{xyy}^3+2E(L_1\partial_x+L_2\partial_y)$ which appears in the Novikov-Neselov equation. } would lead to a local well-posedness result when $E\neq 0,$ not in the "optimal" space though. 

Note than one cannot perform simple energy estimates on \eqref{NV} 
(as in the KP or Zakharov-Kuznetsov case) leading "for free" to a local result in $H^s(\R^2), s>2$ (actually even the $L^2$ estimate on \eqref{NV} fails).

Also PDE methods should solve the aforementioned problem of the "long wave" KP limit when $E\to \pm \infty.$ Actually similar questions are already solved for the Gross-Pitaevskii equation (see \cite{CR, BGSS1, BGSS2}).

Finally we do not know of any stability results (by PDE methods) for the Novikov-Veselov solitary waves.

\textbf{Acknowledgement:}\\
The Authors thank Anne de Bouard, Anna Kazeykina, Thanasis Fokas,  Oleg Kiselev and P. Perry for fruitful conversations.  J.-C. S. was partially supported  by the project GEODISP of the Agence Nationale de la Recherche.

\end{document}